\documentclass[11pt]{article}
\usepackage[margin=2.0cm]{geometry}
\usepackage[utf8]{inputenc}
\usepackage[T1]{fontenc}
\usepackage[english]{babel}
\usepackage{amsmath,amsfonts,amssymb,amsthm,mathtools}
\usepackage{bm}
\usepackage{graphicx}
\usepackage{float}
\usepackage{booktabs}
\usepackage[numbers]{natbib}
\usepackage[colorlinks=true,linkcolor=blue,citecolor=blue,urlcolor=blue]{hyperref}
\allowdisplaybreaks\hbadness=10000
\numberwithin{equation}{section}

\theoremstyle{plain}
\newtheorem{theorem}{Theorem}[section]
\newtheorem{proposition}[theorem]{Proposition}
\newtheorem{lemma}[theorem]{Lemma}
\newtheorem{corollary}[theorem]{Corollary}
\theoremstyle{definition}
\newtheorem{remark}[theorem]{Remark}

\newcommand{\E}{\mathbb{E}}
\newcommand{\PP}{\mathbb{P}}
\newcommand{\R}{\mathbb{R}}
\newcommand{\Sph}{\mathbb{S}}
\newcommand{\1}{\mathbf{1}}
\newcommand{\Var}{\mathrm{Var}}
\newcommand{\Cov}{\mathrm{Cov}}
\newcommand{\GOE}{\mathrm{GOE}}
\DeclareMathOperator{\diag}{diag}
\newcommand{\Mb}{\overline{M}}
\newcommand{\dex}{\delta_{\mathrm{exact}}}
\newcommand{\dst}{\delta_{*}}
\newcommand{\Nlm}{N^{\mathrm{lm}}}
\newcommand{\Ncp}{N^{\mathrm{cp}}}
\newcommand{\lmax}{\lambda_{\max}}

\newcommand{\verifrepo}{\href{https://github.com/ydecastro/TENSOR-KSS}{\texttt{github.com/ydecastro/TENSOR-KSS}}}
\newcommand{\doi}[1]{\href{https://doi.org/#1}{\nolinkurl{https://doi.org/#1}}}

\title{Non-asymptotic First- and Second-Order Tail Brackets\\ for the Kostlan--Shub--Smale Field\\ and Spherical $k$-Spin Complexity}
\usepackage{authblk}
\author[1]{Jean-Marc Aza\"{i}s}
\author[2]{Federico Dalmao}
\author[3,4]{Yohann De Castro\thanks{Corresponding author. E-mail: \texttt{yohann.de-castro@ec-lyon.fr}. ORCID: 0000-0002-9008-7474.}}
\affil[1]{{\small Institut de Math\'ematiques de Toulouse\\ Universit\'e de Toulouse, France.}}
\affil[2]{{\small DMEL, CENUR Litoral Norte\\ Universidad de la Rep\'ublica, Salto, Uruguay.}}
\affil[3]{{\small Institut Camille Jordan, CNRS UMR 5208\\ \'Ecole Centrale Lyon, France.}}
\affil[4]{{\small Institut Universitaire de France~(IUF)}}
\date{November 6, 2026}

\begin{document}
\maketitle

\begin{abstract}
For the Kostlan--Shub--Smale field $X$ of degree $k\ge3$ on the sphere $\Sph^{d-1}$, 
we prove the two-sided bound
$(1-\varepsilon(u))\,\dst(u)\le\PP\{\sup X>u\}\le\dst(u)$ at every finite $(k,d)$, where
$\dst(u)$ is the expected number of critical points above $u$, halved when $k$ is even, and $\varepsilon(u)$ is explicit for
$\rho u\ge\min\{d-1,\tfrac32\sqrt{2(d-1)}\}$. The gap is therefore at most of order $e^{-(1+\rho^{2})u^{2}/2}$ with
$\rho^{2}=k/(2(k-1))$.
We also give $\dst$ itself in
closed form, as a finite sum of Hermite polynomials, together with a two-term envelope that is asymptotically
sharp and a single-monomial relaxation of the form $Cu^{d-2}e^{-u^{2}/2}$. As an application we bracket the annealed
complexity of the spherical $k$-spin Hamiltonian at every finite $(k,d)$: the upper bound holds at every
energy, a lower bound of the same form at every energy above $9/4$ times the spectral edge $\sqrt{2(d-1)}/\rho$, above twice
the edge when $d\ge10$, and is consistent with the Auffinger--Ben~Arous--\v{C}ern\'y complexity function as $d\to\infty$.
\end{abstract}



\section{Introduction}\label{sec:intro}

Let $X$ be the Kostlan--Shub--Smale (KSS) field of degree $k\ge3$ on the unit sphere $\Sph^{d-1}\subset\R^{d}$,
$d\ge3$: the centered Gaussian field with covariance $\E[X(s)X(t)]=\langle s,t\rangle^{k}$. It is the canonical
isotropic Gaussian homogeneous polynomial of degree $k$, the Hamiltonian of the spherical $k$-spin model up to
normalization, and the noise part of the spiked tensor model; its paths are polynomials. This paper is about one quantity,
$\PP\{\sup_{\Sph^{d-1}}X>u\}$, and about how well the Kac--Rice first moment describes it at finite $(k,d)$.

Two things are proved. The first is that the Kac--Rice bound is not only an upper bound but a two-sided one,
with an explicit error: Theorem~\ref{thm:target} brackets the supremum probability between
$(1-\varepsilon(u))\,\dst(u)$ and $\dst(u)$, where $\dst(u)$ is the expected number of critical points above $u$ (halved for $k$ even)
and $\varepsilon(u)$ is explicit at every finite $(k,d)$. The second is that $\dst$ itself is computable: it has
a closed form as a finite sum of Hermite polynomials (Theorem~\ref{thm:imf_tail_exact}), together with two
explicit relaxations (Theorems~\ref{thm:imf_tail} and~\ref{thm:smf_tail}) that trade tightness for a form one can
invert. Section~\ref{sec:kspin} applies both to the annealed complexity of the spherical $k$-spin Hamiltonian.

The point of the bracket is that it turns a bound into an approximation. Without it, an explicit upper bound on a
tail can always be suspected of being loose by an unknown factor, and there is no baseline against which to judge
it. Here the baseline is $\dst$ itself, and the bracket says the first moment is right to within
$\varepsilon(u)$, which decays like $e^{-\rho^{2}u^{2}/2}$ with $\rho^{2}=k/(2(k-1))$.

\subsection{The Kac--Rice integral and the conditional Hessian}\label{sub:reduction}\label{sec:main}

Throughout, $\rho:=\sqrt{k/(2(k-1))}$ and $\varsigma^{2}:=1/\rho^{2}=2(k-1)/k$; $\Phi$ and $\varphi$ are the standard Gaussian distribution function and density, $\bar\Phi:=1-\Phi$; $\lmax$ is the largest eigenvalue of a symmetric matrix,~$\|G\|$ the operator norm of a matrix $G$, and $\|Y\|_{p}:=(\E|Y|^{p})^{1/p}$ for a random variable $Y$; $\Ncp_{[u,\infty)}(F)$ and $\Nlm_{[u,\infty)}(F)$ count the critical points and the local maxima of a smooth function $F$ with value $\ge u$, the argument being omitted when $F=X$. The Gaussian Orthogonal Ensemble of size $n$, $\GOE(n)$, is the law of a real symmetric $n\times n$ random matrix $G=(G_{ij})$ whose entries are jointly Gaussian, centered, independent up to symmetry, with the normalization $\Var(G_{ii})=1$ and $\Var(G_{ij})=1/2$ for $i\neq j$; equivalently, the density on the space of real symmetric matrices is proportional to $\exp(-\mathrm{tr}(G^2)/2)$ and the law is invariant under conjugation by any orthogonal matrix~\citep{mehta}. 

Fix $t\in\Sph^{d-1}$ and an orthonormal basis $(e_{a})_{a\le d-1}$ of $t^{\perp}=T_{t}\Sph^{d-1}$. Differentiating $\langle s,t\rangle^{k}$, the covariant derivatives of $X$ at $t$ satisfy
\begin{gather*}
\Cov(\partial_{a}X(t),\partial_{b}X(t))=k\delta_{ab},\qquad \Cov(X(t),\nabla X(t))=0,\qquad \Cov(\nabla X(t),\nabla^{2}X(t))=0,\\
\Cov(X(t),\nabla^{2}_{ab}X(t))=-k\delta_{ab},\qquad
\Cov(\nabla^{2}_{ab}X(t),\nabla^{2}_{cd}X(t))=k(k-1)(\delta_{ac}\delta_{bd}+\delta_{ad}\delta_{bc})+k^{2}\delta_{ab}\delta_{cd}.
\end{gather*}
Hence the matrix
\begin{subequations}\label{eq:kacrice_setup}
\begin{equation}\label{eq:G}
G:=\frac{\nabla^{2}X(t)+k\,X(t)\,I_{d-1}}{\sqrt{2k(k-1)}}\;\sim\;\GOE(d-1)
\end{equation}
is independent of $(X(t),\nabla X(t))$, the gradient $\nabla X(t)\sim N(0,kI_{d-1})$ is independent of $(X(t),\nabla^{2}X(t))$, and conditionally on $\{X(t)=x,\ \nabla X(t)=0\}$,
\begin{equation}\label{eq:condhess}
\nabla^{2}X(t)=\sqrt{2k(k-1)}\,\big(G-\rho x\,I_{d-1}\big),
\end{equation}
since $\sqrt{2k(k-1)}\,\rho=k$; this is the conditional law of \citet[Lemma~3.2(b)]{auffinger2013}. The critical point~$t$ is a nondegenerate local maximum if and only if $\lmax(G)<\rho x$. The Kac--Rice formula for the critical points of $X$ \citep[Theorem~6.1]{armentano2023general}, with $p_{\nabla X(t)}(0)=(2\pi k)^{-(d-1)/2}$ and \eqref{eq:condhess}, gives
\begin{equation}\label{eq:dexact}
\begin{gathered}
\dex(u):=\E\,\Ncp_{[u,\infty)}=C_{k,d}\int_{u}^{\infty}A(\rho x)\,\varphi(x)\,\mathrm{d}x,\\
A(\nu):=\E\big|\det(G-\nu I_{d-1})\big|,
\qquad
C_{k,d}:=\mathrm{Vol}(\Sph^{d-1})\Big(\frac{k-1}{\pi}\Big)^{\frac{d-1}{2}}=\frac{2\sqrt{\pi}}{\Gamma(d/2)}\,(k-1)^{\frac{d-1}2},
\end{gathered}
\end{equation}
an equality for the expected number of critical points above $u$, which Theorem~\ref{thm:imf_tail_exact} computes in closed form. For $k$ even, $X(-s)=X(s)$ and the critical points come in antipodal pairs. 
Since the maximum of~$X$ over the compact $\Sph^{d-1}$ is attained at a critical point, $\1\{\sup X>u\}\le\Ncp_{[u,\infty)}$, and taking expectations 
\begin{equation}\label{eq:KR_main_bound_intro}
    \PP\Big\{\sup_{\Sph^{d-1}}X>u\Big\}\;\le\;\dst(u)\;:=\frac{2}{3+(-1)^k}\;\dex(u)\,.
\end{equation}
\end{subequations}
The bound is loose by the critical points above $u$ that are not the global maximum, saddle points and beaten local maxima alike; Theorem~\ref{thm:target} shows that this slack is a factor $\varepsilon(u)=O(u^{d-2}e^{-\rho^{2}u^{2}/2})$. The bound~\eqref{eq:KR_main_bound_intro} is sharp at second order. For a Gaussian field of constant variance on a parameter set with bounded curvature, \citet[Theorem~8.12]{azais2009level} prove that the Rice bound on the density of the maximum exceeds that density by at most $\exp(-(1+\gamma)u^2/2+o(u^2))$, with $\gamma>0$ given by the conditional variance of the field at two points (and a boundary term, absent here). The scale $\dex$ is therefore the right one for the supremum probability as $u\to\infty$; the constants of that theorem are asymptotic in the level, and Theorem~\ref{thm:target} gives the matching lower bound at finite $(k,d)$ by another route.

Section~\ref{sec:tail_bounds} evaluates $\dex$ in closed form for every $u$ (Theorem~\ref{thm:imf_tail_exact}) and gives $\delta_{\mathrm{IMF}}\ge\dex$ on $[u_{\mathrm{IMF}},\infty)$, with $u_{\mathrm{IMF}}:=\sqrt{2d-1}/\rho$,
and $\delta_{\mathrm{SMF}}\ge\delta_{\mathrm{IMF}}$ on $[u_{\mathrm{SMF}},\infty)$ (Proposition~\ref{thm:uniform_domination}).

\paragraph{Asymptotic baseline as $u\to\infty$, with $d$ fixed.}
\label{sub:baseline}
Replacing $A(\rho x)$ in~\eqref{eq:dexact} by its leading monomial~$(\rho x)^{d-1}$ gives, by Lemma~\ref{lem:exp_det}, the asymptotic equivalent
\begin{equation}\label{eq:asymptotic_baseline_approx}
    \delta_{\mathrm{bl}}(u)\;:=\;\frac{\sqrt{2}}{\Gamma(d/2)}\!\left(\frac{k}{2}\right)^{\frac{d-1}2}u^{d-2}\,e^{-u^2/2}\,,
    \qquad
    \dex(u)\;\sim\;\delta_{\mathrm{bl}}(u)\quad\text{as $u\to\infty$ with $d$ fixed}\,.
\end{equation}
This is the optimal tail rate against which the non-asymptotic bounds are measured; $\delta_{\mathrm{bl}}$ is not itself an upper bound on $\dex$, only the asymptotic equivalent.

\paragraph{Constants.} The constants $\alpha_d,\beta_d$ and the functions $\Phi_d,\Psi_d$ of Table~\ref{tab:tail_bounds} are defined in \eqref{eq:alpha_d}, \eqref{eq:beta_d} and \eqref{eq:phi_def}--\eqref{eq:psi_def}; all are explicit at finite $(k,d)$.

\begin{table}[ht]
\centering
\renewcommand{\arraystretch}{1.5}
\resizebox{\textwidth}{!}{
\begin{tabular}{lll}
\toprule
\textbf{Method} & \textbf{Validity range} & \textbf{Closed-form expression} \\
\midrule
Baseline (not a bound) & $u\to\infty$, $d$ fixed & $\delta_{\mathrm{bl}}(u)\;=\;\dfrac{\sqrt{2}}{\Gamma(d/2)}\!\left(\dfrac{k}{2}\right)^{\frac{d-1}2}u^{d-2}e^{-u^2/2}$ \\
$\delta_{\mathrm{exact}}$ (canonical) & $u\in\R$ (exact) & $\delta_{\mathrm{exact}}(u)\;=\;2(k-1)^{\frac{d-1}2}\bigl[D_1(u)+D_2(u)+D_3(u)+D_4(u)\bigr]$, finite Hermite sum (Thm.~\ref{thm:imf_tail_exact}) \\
$\delta_{\mathrm{IMF}}$ & $u\ge u_{\mathrm{IMF}}$ & $\delta_{\mathrm{IMF}}(u)\;=\;2(k-1)^{\frac{d-1}2}\,\Big(\alpha_d\,\Phi_d(\rho,u)\,e^{-u^2/2}+\Psi_d(\rho,u)\,e^{-(1+\rho^2)u^2/2}\Big)$ \\
$\delta_{\mathrm{SMF}}$ & $u\ge u_{\mathrm{SMF}}$ & $\delta_{\mathrm{SMF}}(u)\;=\;4\,\alpha_d(2k)^{\frac{d-1}2}u^{d-2}e^{-u^2/2}+2^d\beta_d(k-1)^{\frac{d-1}2}u^{2d-3}e^{-3u^2/4}$ \\
\bottomrule
\end{tabular}
}
\vspace{0.5em}
\caption{The Kac--Rice integral $\dex(u)$ in closed form and its two explicit non-asymptotic relaxations, with the asymptotic baseline $\delta_{\mathrm{bl}}$ in the first row, which is the equivalent $\dex\sim\delta_{\mathrm{bl}}$ as $u\to\infty$ at fixed $d$ and not an upper bound. The exact form $\delta_{\mathrm{exact}}$ is the Kac--Rice integral itself, evaluated in closed form for every $u$; the threshold $u_{\mathrm{IMF}}$ marks the range on which $\delta_{\mathrm{exact}}\le\delta_{\mathrm{IMF}}$. The constants $\alpha_{d},\beta_{d}$ and the functions $\Phi_{d},\Psi_{d}$ are defined in \eqref{eq:alpha_d}, \eqref{eq:beta_d} and \eqref{eq:phi_def}--\eqref{eq:psi_def}.}
\label{tab:tail_bounds}
\end{table}

\subsection{The main result}\label{sec:mainresult}

For $t\in\Sph^{d-1}$, let $G$ be the matrix \eqref{eq:G} and $\zeta^{t}$ the helix process of Section~\ref{sec:gap}, the normalized residual of the regression of $X$ on $(X(t),\nabla X(t))$, defined on $M=\Sph^{d-1}\setminus\{\pm t\}$; put
\begin{equation}\label{eq:Xidef}
\Xi(x):=\E\Big[\big|\det(\rho xI-G)\big|\;\PP\Big(\sup_{M}\zeta^{t}>x\;\Big|\;G\Big)\Big],
\end{equation}
which does not depend on $t$ by rotational invariance. Lemma~\ref{lem:beaten} and Proposition~\ref{prop:mult} are stated for $k$ odd; Sections~\ref{sec:layers} to~\ref{sec:rice} hold for both parities, and Section~\ref{sec:final} treats $k$ even.

\newpage

\begin{theorem}[Non-asymptotic second-order bracket]\label{thm:target}
Let $k\ge3$ and $d\ge3$. For every $u>0$,
\[
\big(1-\varepsilon(u)\big)\,\dst(u)\;\le\;\PP\Big\{\sup_{\Sph^{d-1}}X>u\Big\}\;\le\;\dst(u),
\qquad
\varepsilon(u):=\sup_{x\ge u}\frac{\Xi(x)}{A(\rho x)}\in[0,1],
\]
with $\Xi$ as in \eqref{eq:Xidef}, and for $\rho u\ge d-1$
\begin{equation}\label{eq:main_bound}
\begin{gathered}
\varepsilon(u)\;\le\;C^{\mathrm{gap}}_{k,d}\;\frac{(1+\rho u)^{d-1}}{\rho u}\;e^{-\rho^{2}u^{2}/2},\\
C^{\mathrm{gap}}_{k,d}=4\Big(2+\tfrac56\sqrt{d-1}\Big)^{d-1}C^{\mathrm{Rice}}_{k,d},
\qquad
C^{\mathrm{Rice}}_{k,d}=c_{k}\,\frac{\sqrt{\pi}\;k^{3}\;2^{2k}\;(k-1)^{(d-2)/2}}{\Gamma\big(\frac{d-1}{2}\big)},\\
c_{k}=43\ \text{for $k$ odd},\qquad c_{k}=\tfrac92\ \text{for $k$ even}.
\end{gathered}
\end{equation}
\end{theorem}

\begin{remark}
    For $\tfrac32\sqrt{2(d-1)}\le\rho u<d-1$, a non-empty range when $d\ge6$, the bound \eqref{eq:main_bound} holds with the factor $4$ of $C^{\mathrm{gap}}_{k,d}$ replaced by $\tfrac{15}{4}\big(\tfrac95\big)^{(d-1)/2}$.
\end{remark}

Hence $\dst(u)-\PP\{\sup X>u\}$ decays at least like $e^{-(1+\rho^{2})u^{2}/2}$, with $1+\rho^{2}=1+k/(2(k-1))$ the second-order exponent given by \citet[Theorem~8.12]{azais2009level}, here with explicit finite-$u$ constants. 

\paragraph{Strategy.} The first half of the proof (Lemma~\ref{lem:beaten} to Corollary~\ref{cor:cs}) writes the gap as the expected number of critical points above $u$ that are beaten by a larger value of $X$. Conditionally on $X(t)=x$ and $\nabla X(t)=0$, the field at polar angle $\theta$ from $t$ is $X(s)=\cos^{k}\theta\,x+(1-\cos^{k}\theta)\,\zeta^{t}(s)$, where $\zeta^{t}$ is the \emph{helix process}, the normalized residual of the regression of $X$ on $(X(t),\nabla X(t))$; so $t$ is beaten if and only if $\sup\zeta^{t}>x$, the gap equals $C_{k,d}\int_{u}^{\infty}\Xi(x)\varphi(x)\,\mathrm{d}x$ with $\Xi(x)=\E[|\det(\rho xI-G)|\,\PP(\sup\zeta^{t}>x\mid G)]$, and H\"older's inequality with the exponent $p=\rho^{2}x^{2}$ separates the weight, an $L^{p}$ norm of the GOE controlled by the Gaussian tail of $\|G\|$, from the tail of $\sup\zeta^{t}$, whose power $1-1/p$ costs only a constant because the fiber value $\varsigma\lmax(G)$ alone gives the lower bound $\PP\{\sup\zeta^{t}>x\}\ge\bar\Phi(\rho x)$. The second half (Sections~\ref{sec:layers} to~\ref{sec:rice}) bounds that tail by the Rice formula. The helix process is a finite sum of independent KSS fields of lower degree on the equatorial sphere, the layers, with weights depending only on $\theta$; its variance is at most~$\varsigma^{2}$, with equality only at the pole, where it reduces to $\varsigma v^{\top}Gv$. Its supremum is attained at an interior local maximum, never on the pole fibers, so $\PP\{\sup\zeta^{t}>x\}$ is at most the expected number of local maxima above $x$. In the Kac--Rice integral for that number, the Gaussian density at polar angle $\theta$ has exponent at least $y^{2}/(2\varsigma^{2})$, and the diagonal of the conditional Hessian consists of independent Gaussian variables, so Hadamard's inequality bounds the integrand by a product of one-dimensional moments; bounds uniform in $\theta$ give $C^{\mathrm{Rice}}_{k,d}$. Section~\ref{sec:final} assembles the theorem.

\subsection{Related work}

The Kostlan--Shub--Smale field was introduced by \citet{shub1993complexity} and \citet{kostlan1993distribution}
for the roots of random polynomial systems; its covariance $\langle s,t\rangle^{k}$ is that of the noise in the
spiked tensor model. The Kac--Rice framework for the extrema of smooth Gaussian fields is that of
\citet{azais2009level}. The second-order term of the tail of the maximum, the object of Theorem~\ref{thm:target}, goes back to \citet{azais2005distribution} and \citet{azais2008general}, and to \citet{taylor2005validity} by the Euler characteristic route. The conditional law of the Hessian at a critical point, a shifted GOE, is Lemma~3.2(b) of
\citet{auffinger2013}, whose Theorems~2.1 and~2.2 are the exact Kac--Rice counts (Proposition~\ref{prop:ortho_exact} and Remark~\ref{rem:kr_representation} below) and whose Theorems~2.5, 2.8 and~2.17 give their large-$d$ asymptotics
above and below the spectral edge; Theorem~\ref{thm:complexity} is the finite-$(k,d)$ counterpart of their
Theorems~2.5 and~2.8 above the edge, with the same conditional law as input and the Mehta--Fyodorov evaluation of the
expected absolute determinant \citep{mehta,fyodorov2004complexity} in place of the large deviations. The
connection with spherical spin glasses runs through the $p$-spin Hamiltonian, whose landscape is studied by
\citet{fyodorov2004complexity} and \citet{auffinger2013}, and for the spiked model by
\citet{arous2019landscape}. Sharpening the Kac--Rice integral uses the orthogonal-polynomial expansion of the GOE
eigenvalue density \citep{mehta} and Szeg\H{o}'s bounds on the zeros of Hermite polynomials \citep{szego}.
On the statistical side, $\sup X$ is the injective norm of the noise tensor: in the normalization of
\citet{perry2020statistical} its limit is their $\mu_{k}$, and the detection threshold for the spherical prior
sits at $\lambda\asymp\sqrt{d\log k}$ in our variables. The statistical questions attached to the same field are treated in a separate paper and are not used here. The Bayes-optimal
side of the same model, the asymptotic mutual information and minimum mean square error of
\citet{lesieur2017statistical}, is the natural target of a finite-$(k,d)$ bracket for the spiked field; the present
paper does not go there, since it never uses the spike.

\section{Exact and asymptotically sharp first-order bounds}\label{sec:first}
\label{sec:tail_bounds}
This section gives the Kac--Rice integral $\dex$ of~\eqref{eq:dexact} in closed form and two explicit relaxations, tightest first: the exact form $\delta_{\mathrm{exact}}$ of Theorem~\ref{thm:imf_tail_exact}, its relaxation $\delta_{\mathrm{IMF}}$ of Theorem~\ref{thm:imf_tail}, and the single-monomial~$\delta_{\mathrm{SMF}}$ of Theorem~\ref{thm:smf_tail}. They are nested: $\delta_{\mathrm{exact}}\le\delta_{\mathrm{IMF}}$ on $[u_{\mathrm{IMF}},\infty)$ (Theorem~\ref{thm:imf_tail_exact}) and $\delta_{\mathrm{IMF}}\le\delta_{\mathrm{SMF}}$ on $[u_{\mathrm{SMF}},\infty)$ (Proposition~\ref{thm:uniform_domination}). Throughout, $\dex$ is an equality for $\E\Ncp_{[u,\infty)}$ and an upper bound for $\PP\{\sup X>u\}$ by~\eqref{eq:KR_main_bound_intro}.

\subsection{Mehta--Fyodorov representation (MF)}
\label{sub:MF}

The Mehta--Fyodorov representation uses two identities: Fyodorov's formula (Lemma~\ref{lfyod}), which writes the expected absolute characteristic polynomial of a $(d-1)$-GOE at level $\nu$ through the one-point eigenvalue density of a one-size-larger $d$-GOE, and Mehta's Hermite expansion of that density (Lemma~\ref{lem:mehta}); together they give an exact closed algebraic form of the Kac--Rice integrand.

\begin{lemma}[Fyodorov's formula~\citep{fyodorov2004complexity}]
\label{lfyod}
For every $\nu\in\R$,
\begin{equation}\label{eq:fyo}
    A(\nu)\;=\;\E\big|\det(G-\nu I_{d-1})\big|
    \;=\;\frac{2^{\frac{3}2}\,\Gamma((d+2)/2)}{d}\;Q_{d}(\nu)\,,
\end{equation}
where $Q_d(\nu):=e^{\nu^2/2}q_d(\nu)$ and $q_d$ is the eigenvalue intensity (one-point correlation function) of a $d\times d$ GOE matrix: $q_d(\nu)\,\mathrm{d}\nu$ is the expected number of eigenvalues in $\mathrm{d}\nu$, so that $\int_\R q_d(\nu)\,\mathrm{d}\nu = d$.
\end{lemma}

The identity expresses the expected absolute characteristic polynomial of the $(d-1)\times(d-1)$ GOE through the one-point function of the one-size-larger $d\times d$ GOE; see \citet{fyodorov2004complexity}.

\begin{lemma}[Mehta expansion~{\citep[p.~221]{mehta}}]
\label{lem:mehta}
With $c_j:=(2^j j!\sqrt{\pi})^{-1/2}$ and $H_j$ the $j$-th physicist Hermite polynomial,
\begin{align}
    Q_d(\nu) &\;=\;
    \frac{1}{2}\sqrt{\frac d2}\,c_{d-1}c_d\,H_{d-1}(\nu)\big[\mu_d-2\mathcal{I}_d^c(\nu)\big]
    \;+\;\1\{d\text{ odd}\}\,\frac{H_{d-1}(\nu)}{\mu_{d-1}}
    +e^{-\nu^2/2}\sum_{j=0}^{d-1} c_j^2\,H_j^2(\nu)\,, \label{eq:fmehta}
\end{align}
where
\[
    \mu_m:=\int_\R H_m(y)\,e^{-y^2/2}\,\mathrm{d}y\,,\qquad
    \mathcal{I}_d^c(\nu):=\int_\nu^\infty H_d(y)\,e^{-y^2/2}\,\mathrm{d}y\,.
\]
\end{lemma}

By parity, $\mu_m=0$ whenever $m$ is odd. For $d$ even the indicator term is absent and the bracket keeps both of its parts. For $d$ odd the constant $\mu_d$ vanishes, so the bracket reduces to $-2\mathcal{I}_d^c(\nu)$ and the term $H_{d-1}(\nu)/\mu_{d-1}$ is present; the subtractive part $-2\mathcal{I}_d^c(\nu)H_{d-1}(\nu)$ survives in both parities, and it is the part discarded by the IMF bound.

The combination of Lemmas~\ref{lfyod} and~\ref{lem:mehta} produces an exact integral representation of the Kac--Rice probability bound. Fyodorov's formula~\eqref{eq:fyo} gives the expected absolute characteristic polynomial as a constant multiple of $Q_d(\rho x)$, and Mehta's expansion~\eqref{eq:fmehta} supplies an explicit formula for $Q_d$. Substituting into~\eqref{eq:dexact} therefore yields an integrand of the form $Q_d(\rho x)\,e^{-x^2/2}$: a polynomial-and-Hermite quantity weighted by the original Kac--Rice Gaussian, with no further exponential factor entering the integrand. Fyodorov's identity is already stated in the natural variable $Q_d$, so the eigenvalue intensity $q_d=e^{-\nu^2/2}Q_d$ never appears in the substitution: the Gaussian envelope that would otherwise be carried by $q_d$ is absorbed into the definition of $Q_d$. This algebraic identity is what the IMF and SMF bounds of Sections~\ref{sub:IMF} and~\ref{sub:SMF} rest on. We state this in the following proposition.

\begin{proposition}[Exact MF representation]
\label{prop:ortho_exact}
For every $u\in\R$,
\begin{equation}\label{eq:ortho_exact}
    \PP\Big\{\sup_{\Sph^{d-1}}X>u\Big\}
    \;\le\;\dex(u)\;=\;
    2(k-1)^{\frac{d-1}2}\int_u^\infty Q_{d}(\rho x)\,e^{-x^2/2}\,\mathrm{d}x\,.
\end{equation}
\end{proposition}

\begin{proof}[Proof of Proposition~\ref{prop:ortho_exact}]

\noindent
$\bullet$
Starting from~\eqref{eq:dexact}, Lemma~\ref{lfyod} evaluated at the shifted argument $\nu=\rho x$ gives an exact expression for the expected absolute determinant directly in terms of $Q_d$:
\begin{equation}\label{eq:fyo_substitute_app}
    A(\rho x)
    \;=\;\frac{2^{\frac{3}2}\Gamma((d+2)/2)}{d}\,Q_d(\rho x)\,,
\end{equation}
where $Q_d(\nu)=e^{\nu^2/2}q_d(\nu)$ is the unweighted polynomial part of the eigenvalue intensity (one-point correlation function) of a $d\times d$ GOE matrix.

\noindent
$\bullet$ 
We multiply this expression by the standard Gaussian density $\varphi(x)=(2\pi)^{-1/2}e^{-x^2/2}$ and integrate over the interval $[u,\infty)$ to match the Kac--Rice integral:
\[
    \int_u^\infty A(\rho x)\,\varphi(x)\,\mathrm{d}x
    \;=\;\frac{2^{\frac{3}2}\Gamma((d+2)/2)}{d\sqrt{2\pi}}\int_u^\infty Q_d(\rho x)\,e^{-x^2/2}\,\mathrm{d}x\,.
\]
We can simplify the constant prefactor using the Gamma function identity $\Gamma((d+2)/2)=(d/2)\Gamma(d/2)$:
\[
    \frac{2^{\frac{3}2}\Gamma((d+2)/2)}{d\sqrt{2\pi}} 
    \;=\; \frac{2\sqrt{2} \cdot (d/2)\Gamma(d/2)}{d \sqrt{2} \sqrt{\pi}} 
    \;=\; \frac{\Gamma(d/2)}{\sqrt{\pi}}\,.
\]
Finally, with $C_{k,d}=2\sqrt\pi(k-1)^{\frac{d-1}2}/\Gamma(d/2)$ from~\eqref{eq:dexact}, the $\Gamma(d/2)$ and $\sqrt{\pi}$ terms cancel:
\[
    C_{k,d} \frac{\Gamma(d/2)}{\sqrt{\pi}} \;=\; \left( \frac{2\sqrt\pi(k-1)^{\frac{d-1}2}}{\Gamma(d/2)} \right) \frac{\Gamma(d/2)}{\sqrt{\pi}} \;=\; 2(k-1)^{\frac{d-1}2}\,.
\]
This yields the exact representation~\eqref{eq:ortho_exact}.
\end{proof}

\subsection{Strictly exact Mehta--Fyodorov tail bound}
\label{sub:exact_mf}

Substituting the full Mehta expansion~\eqref{eq:fmehta} into Proposition~\ref{prop:ortho_exact} produces the four-piece decomposition
\begin{subequations}\label{eq:Q_d_pieces}
\begin{equation}\label{eq:Q_d_decomp}
    \int_u^\infty Q_d(\rho x)\,e^{-x^2/2}\,\mathrm{d}x \;=\; D_1(u)\;+\;D_2(u)\;+\;D_3(u)\;+\;D_4(u)\,,
\end{equation}
where, with $\nu=\rho x$ and $\Lambda:=2\rho^2-1=1/(k-1)$,
\begin{align*}
    D_1(u) &\;:=\;\tfrac12\sqrt{d/2}\,c_{d-1}c_d\,\mu_d\,\mathcal{L}_d(u)\,,\quad &\text{(zero for $d$ odd, since $\mu_d=0$),}\\
    D_2(u) &\;:=\;-\sqrt{d/2}\,c_{d-1}c_d\,\mathcal{C}_d(u)\,,\quad &\text{(cross integral $\mathcal{C}_d$ defined in~\eqref{eq:Cd_def}),}\\
    D_3(u) &\;:=\;\1\{d\,\mathrm{odd}\}\,\mu_{d-1}^{-1}\,\mathcal{L}_d(u)\,,\quad &\text{(zero for $d$ even),}\\
    D_4(u) &\;:=\;T_d^{\mathrm{exact}}(u)\quad &\text{(squared-Hermite tail of \eqref{eq:Td_exact_intro}),}
\end{align*}
with $\mathcal{L}_d(u):=\int_u^\infty H_{d-1}(\rho x)\,e^{-x^2/2}\,\mathrm{d}x$ and the cross integral
\begin{equation}\label{eq:Cd_def}
    \mathcal{C}_d(u)\;:=\;\int_u^\infty H_{d-1}(\rho x)\,\mathcal{I}_d^c(\rho x)\,e^{-x^2/2}\,\mathrm{d}x\,.
\end{equation}
\end{subequations}

\begin{theorem}[Strictly exact MF tail bound]
\label{thm:imf_tail_exact}
For every $u\in\R$,
\begin{subequations}\label{eq:exact_closed}
\begin{equation}\label{eq:exact_bound}
    \PP\Big\{\sup_{\Sph^{d-1}}X>u\Big\}
    \;\le\;
    \delta_{\mathrm{exact}}(u)
    \;=\;
    2(k-1)^{\frac{d-1}2}\bigl[D_1(u)+D_2(u)+D_3(u)+D_4(u)\bigr]\,,
\end{equation}
where each $D_i$ is given in fully closed form by:
\begin{enumerate}
    \item[$(D_1,D_3)$] the linear Hermite tail $\mathcal{L}_d(u)$ admits the closed form
        \begin{equation}\label{eq:Ld_closed}
            \mathcal{L}_d(u)\;=\;\frac1\rho\,J_{d-1}^{1/\rho^2}(\rho u)\,,
        \end{equation}
        with $J_m^\beta$ given by Lemma~\ref{lem:Jm_beta};
    \item[$(D_4)$] the squared Hermite tail $T_d^{\mathrm{exact}}(u)=\rho^{-1}\sum_{j=0}^{d-1}c_j^2\,K_j^\beta(\rho u)$ of~\eqref{eq:Td_exact_intro}, with $\beta=(3k-2)/k$ and $K_j^\beta$ given by Corollary~\ref{cor:Kj_beta};
    \item[$(D_2)$] the cross integral $\mathcal{C}_d(u)$ admits the two-piece closed form
        \begin{align}
            \mathcal{C}_d(u) &\;=\;\frac1\rho\sum_{\ell=0}^{\lfloor (d-1)/2\rfloor}\frac{2^{\ell+1}\,(d-1)!!}{(d-2\ell-1)!!}\sum_{p=0}^{d-2\ell-1}2^p\,p!\,\binom{d-1}{p}\binom{d-2\ell-1}{p}\,J_{2d-2\ell-2-2p}^\beta(\rho u)\notag\\
            &\quad\;+\;\1\{d\,\mathrm{even}\}\,2^{d/2}\,(d-1)!!\,\mathcal{F}_d(u)\,,\label{eq:Cd_closed}
        \end{align}
    where the Fubini remainder $\mathcal{F}_d(u)$ (only present for $d$ even) reads
        \begin{equation}\label{eq:Fubini_F}
            \mathcal{F}_d(u)\;=\;\sqrt{2\pi}\,\bar\Phi(\rho u)\,\frac{J_{d-1}^{1/\rho^2}(\rho u)}{\rho}
            \;-\;2\rho\sum_{\ell=0}^{(d-2)/2}(2\Lambda)^\ell\,\frac{(d-2)!!}{(d-2\ell-2)!!}\,J_{d-2\ell-2}^\beta(\rho u)\,.
        \end{equation}
\end{enumerate}
\end{subequations}
Moreover, with $u_{\mathrm{IMF}}$, $\alpha_d$, $\Phi_d$ and $\delta_{\mathrm{IMF}}$ as in Proposition~\ref{prop:IMF_decomp} and Theorem~\ref{thm:imf_tail},
\[
    \begin{gathered}
    \delta_{\mathrm{exact}}(u)\;<\;\delta_{\mathrm{IMF}}^\star(u)\;<\;\delta_{\mathrm{IMF}}(u)\qquad\text{for every $u\ge u_{\mathrm{IMF}}$,}\\
    \text{where }\ \delta_{\mathrm{IMF}}^\star(u):=2(k-1)^{\frac{d-1}2}\big[\alpha_d\Phi_d(\rho,u)e^{-u^2/2}+D_4(u)\big].
    \end{gathered}
\]
\end{theorem}

\noindent
Since $\beta>2$, the sums of Lemma~\ref{lem:Jm_beta} alternate. The quantity $2(k-1)^{(d-1)/2}\sum_{i=1}^4 D_i(u)$ is the closed-form evaluation of the Kac--Rice integral $2(k-1)^{(d-1)/2}\int_u^\infty Q_d(\rho x)e^{-x^2/2}\,\mathrm{d}x$ of Proposition~\ref{prop:ortho_exact}, that is of $\dex(u)=\E\Ncp_{[u,\infty)}$ in~\eqref{eq:dexact}; the inequality~\eqref{eq:exact_bound} for $\PP\{\sup X>u\}$ is~\eqref{eq:KR_main_bound_intro}. Both readings are used below: the second in Sections~\ref{sec:gap} to~\ref{sec:final}, where the upper bound is matched by a lower one, and the first in Section~\ref{sec:kspin}, where the count is the annealed complexity. For $k$ even, $X(-s)=X(s)$ and the critical points come in antipodal pairs, so Proposition~\ref{prop:keven} halves every bound of this section on $\PP\{\sup X>u\}$; the bounds on $\E\Ncp_{[u,\infty)}$ are unchanged.

\begin{proof}
Deferred to Appendix~\ref{app:deferred}.
\end{proof}

\begin{figure}[t]
\centering
\includegraphics[width=0.99\linewidth]{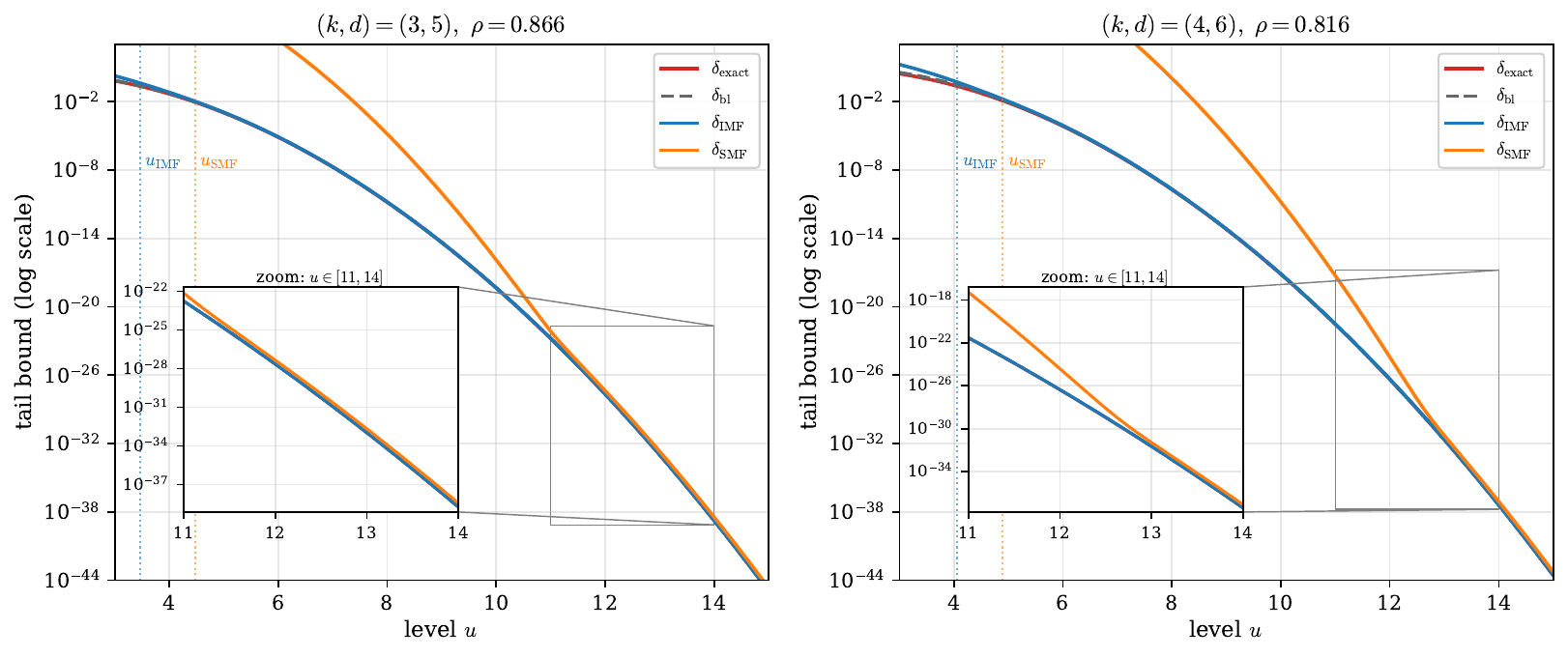}
\caption{The bounds $\delta_{\mathrm{exact}}$ (Theorem~\ref{thm:imf_tail_exact}), $\delta_{\mathrm{IMF}}$
(Theorem~\ref{thm:imf_tail}) and $\delta_{\mathrm{SMF}}$ (Theorem~\ref{thm:smf_tail}), and the asymptotic
baseline $\delta_{\mathrm{bl}}$ of Table~\ref{tab:tail_bounds}, which is not a bound, on logarithmic scale
against the level $u$, for $(k,d)=(3,5)$ (left) and $(4,6)$ (right). The dotted vertical lines mark the
validity thresholds $u_{\mathrm{IMF}}$ and $u_{\mathrm{SMF}}$; each curve is drawn on the whole window. At this
scale $\delta_{\mathrm{IMF}}$ is indistinguishable from $\delta_{\mathrm{exact}}$, the ratio being $1.014$ at
$u=12$ for $(3,5)$ and $1.032$ for $(4,6)$; $\delta_{\mathrm{bl}}$ tracks $\delta_{\mathrm{exact}}$ from
above, the ratio $\delta_{\mathrm{exact}}/\delta_{\mathrm{bl}}$ being $0.993$ at $u=12$ for $(3,5)$ and $0.976$ for $(4,6)$;
$\delta_{\mathrm{SMF}}$ sits above $\delta_{\mathrm{exact}}$ by the factor $2$ of its main term at $(3,5)$, a ratio of
$2.02$ at $u=12$, and by its remainder term $u^{2d-3}e^{-3u^{2}/4}$ at $(4,6)$, a ratio of about $91$ at $u=12$. The inset zooms on
$u\in[11,14]$.}
\label{fig:delta_min}
\end{figure}

Figure~\ref{fig:delta_min} compares the three bounds of Table~\ref{tab:tail_bounds} with the baseline at
$(k,d)=(3,5)$ and $(4,6)$; the figure is produced by the notebook of \verifrepo.

\subsection{Improved Mehta--Fyodorov bound (IMF)}
\label{sub:IMF}

Inspecting the Mehta expansion~\eqref{eq:fmehta}, one sees that the coefficient of the Hermite polynomial $H_{d-1}(\nu)$ is the bracket $[\mu_d-2\mathcal{I}_d^c(\nu)]$, where $\mathcal{I}_d^c(\nu):=\int_\nu^\infty H_d(y)e^{-y^2/2}\,\mathrm{d}y$ is a Hermite tail integral. If a regime can be identified in which $\mathcal{I}_d^c(\nu)$ is \emph{strictly positive}, the subtractive contribution $-2\mathcal{I}_d^c(\nu)H_{d-1}(\nu)$ to $Q_d$ is strictly negative (since $H_{d-1}(\nu)>0$ in the same regime), and discarding it yields a strict upper bound on the Kac--Rice integral.

Pinpointing this regime reduces to locating the largest root of $H_{d-1}$: beyond it, $H_{d-1}(\nu)$ is positive, and by a root-interlacing argument $\mathcal{I}_d^c(\nu)$ is too. Hence the threshold is set by the largest root of $H_{d-1}$. Szeg\H{o}'s classical estimate~\citep[Theorem~6.32]{szego} locates this largest root strictly below $\sqrt{2d-1}$, so that the relevant regime is $\nu\ge\sqrt{2d-1}$, which translates into the threshold $\rho u\ge\sqrt{2d-1}$ in the Kac--Rice integral.

\begin{lemma}[Positivity of the Hermite tail integral]
\label{lem:Idc_positivity}
Let $d\ge 3$ and let $x_1^{(m)}$ denote the largest root of the (physicist) Hermite polynomial $H_m$. For every $\nu\ge x_1^{(d-1)}$, the tail integral $\mathcal{I}_d^c(\nu)$ is strictly positive.
\end{lemma}

\begin{proof}[Proof of Lemma~\ref{lem:Idc_positivity}]

\noindent
$\bullet$
We apply the explicit finite series expansion from Lemma~\ref{lem:hermite_tail_explicit} (specifically Equation~\eqref{eq:hermite_tail_series_1}) with degree $m=d$. This yields:
\[
    \mathcal{I}_d^c(\nu) = e^{-\nu^2/2} \sum_{\ell=0}^{\lfloor (d-1)/2 \rfloor} 2^{\ell+1} \frac{(d-1)!!}{(d-2\ell-1)!!} H_{d-2\ell-1}(\nu) + R_d(\nu)\,,
\]
where the remainder term guarantees $R_d(\nu) \ge 0$.

\noindent
$\bullet$
The roots of consecutive orthogonal polynomials interlace~\citep[Theorem~3.3.2]{szego}; in particular, the largest root $x_1^{(m)}$ of $H_m$ is strictly increasing in the degree $m$. Therefore, for every $\ell \ge 1$, we have the strict inequality:
\[
    x_1^{(d-2\ell-1)} < x_1^{(d-1)}\,.
\]
When we evaluate the sum at any $\nu \ge x_1^{(d-1)}$, the leading polynomial ($\ell=0$) satisfies $H_{d-1}(\nu) \ge 0$. For all subsequent lower-degree polynomials ($\ell \ge 1$; there is at least one since $d\ge3$), since $\nu \ge x_1^{(d-1)} > x_1^{(d-2\ell-1)}$, they evaluate to strictly positive values: $H_{d-2\ell-1}(\nu) > 0$. Because all the constant combinatorial coefficients in the sum are strictly positive and $R_d(\nu) \ge 0$, the entire expression evaluates to a strictly positive number: $\mathcal{I}_d^c(\nu) > 0$.
\end{proof}

\medskip

With the positivity of $\mathcal{I}_d^c$ secured on $\{\nu\ge x_1^{(d-1)}\}$, the subtractive contribution in the Mehta expansion may then be discarded. The two remaining pieces (the linear Hermite tail $\int_u^\infty H_{d-1}(\rho x)e^{-x^2/2}\,\mathrm{d}x$ and the squared Hermite tail $\sum_j c_j^2\int_u^\infty H_j^2(\rho x)e^{-(1+\rho^2)x^2/2}\,\mathrm{d}x$) decay at different Gaussian rates, $e^{-u^2/2}$ and $e^{-(1+\rho^2)u^2/2}$ respectively. This two-scale structure is preserved in the proposition below; the SMF bound of the next subsection collapses it to a single dominant scale. Analytically, we evaluate the linear tail by iterating the three-term Hermite recurrence, and the squared tail by substituting the envelope $H_j(\rho x)\le(2\rho x)^j$ of Lemma~\ref{lem:hermite_bound} into each squared term and invoking Lemma~\ref{lem:gauss_tail}.

\begin{proposition}[Explicit IMF upper bound on $Q_d$]
\label{prop:IMF_decomp}
For every $u>0$ such that $\rho u\ge\sqrt{2d-1}$,
\begin{subequations}\label{eq:IMF_bound}
\begin{equation}\label{eq:IMF_decomp}
    \int_u^\infty Q_d(\rho x)\,e^{-x^2/2}\,\mathrm{d}x
    \;\le\;
    \alpha_d\,\Phi_d(\rho,u)\,e^{-u^2/2}\;+\;\Psi_d(\rho,u)\,e^{-(1+\rho^2)u^2/2}\,,
\end{equation}
where $\alpha_d$ is the dominant Mehta coefficient of Lemma~\ref{lem:dominant_term}, and the explicit polynomial--rational functions are given, with $\Lambda=2\rho^2-1$, by
\begin{align}
    \Phi_d(\rho,u) &\;:=\;\sum_{\ell=0}^{\lfloor (d-2)/2\rfloor}\!2\rho\,(2\Lambda)^\ell\,\frac{(d-2)!!}{(d-2\ell-2)!!}\,(2\rho u)^{d-2\ell-2}
    \;+\;\frac{\1\{d\text{ odd}\}}{u}\,(2\Lambda)^{\frac{d-1}2}\,(d-2)!!\,,\label{eq:phi_def}\\
    \Psi_d(\rho,u) &\;:=\;\frac{c_0^2}{(1+\rho^2)\,u}\;+\;\sum_{j=1}^{d-1}c_j^2\,\frac{(2\rho)^{2j}\,u^{2j-1}}{1+\rho^2-(2j-1)/u^2}\,.\label{eq:psi_def}
\end{align}
\end{subequations}
\end{proposition}

\begin{proof}[Proof of Proposition~\ref{prop:IMF_decomp}]

\noindent
$\bullet$
By Szeg\H{o}'s bound, the largest root of $H_{d-1}$ strictly satisfies $x_1^{(d-1)}<\sqrt{2d-1}$. The proposition assumes $\rho u\ge\sqrt{2d-1}$, ensuring that $\rho x > x_1^{(d-1)}$ for all integration variables $x\ge u$. By Lemma~\ref{lem:Idc_positivity}, this guarantees $\mathcal{I}_d^c(\rho x)>0$ over the entire integration domain. 
Returning to the exact Mehta expansion~\eqref{eq:fmehta}, the contribution of $-2\mathcal{I}_d^c(\rho x) H_{d-1}(\rho x)$ is therefore strictly negative. Discarding it yields a strict upper bound:
\begin{equation}\label{eq:IMF_split}
    \int_u^\infty Q_d(\rho x)e^{-x^2/2}\,\mathrm{d}x \;\le\; \alpha_d\,\mathcal{L}_d(u)\;+\;\mathcal{S}_d(u)\,,
\end{equation}
where $\alpha_d$ is the dominant Mehta coefficient given by~\eqref{eq:alpha_d}, and the components are defined as:
\begin{align*}
    \mathcal{L}_d(u) &:= \int_u^\infty H_{d-1}(\rho x)e^{-x^2/2}\,\mathrm{d}x\,, \\
    \mathcal{S}_d(u) &:= \sum_{j=0}^{d-1}c_j^2\int_u^\infty H_j^2(\rho x)e^{-(1+\rho^2)x^2/2}\,\mathrm{d}x\,.
\end{align*}

\noindent
$\bullet$
We derive a generalized recurrence for $J_m := \int_u^\infty H_m(\rho x)e^{-x^2/2}\,\mathrm{d}x$. Applying the identity $H_m(\rho x)=2\rho x H_{m-1}(\rho x)-2(m-1)H_{m-2}(\rho x)$ and evaluating the first term via integration by parts, we obtain:
\begin{equation*}
    J_m = 2\rho\,H_{m-1}(\rho u)\,e^{-u^2/2} + 2\Lambda(m-1)\,J_{m-2}\,, \qquad \text{where } \Lambda = 2\rho^2-1\,,
\end{equation*}
where we used $\frac{\mathrm{d}}{\mathrm{d}x}[H_{m-1}(\rho x)e^{-x^2/2}]=2\rho(m-1)H_{m-2}(\rho x)e^{-x^2/2}-xH_{m-1}(\rho x)e^{-x^2/2}$. Unrolling this recurrence for $m=d-1$ produces a finite series of boundary evaluations. We bound each boundary evaluation using $H_{d-2\ell-2}(\rho u) \le (2\rho u)^{d-2\ell-2}$ (Lemma~\ref{lem:hermite_bound}). If $d-1$ is even (i.e., $d$ is odd), the recurrence bottoms out with a Gaussian tail, which we bound using Mills' ratio $\int_u^\infty e^{-x^2/2}\,\mathrm{d}x \le u^{-1}e^{-u^2/2}$. Factoring out $e^{-u^2/2}$ yields exactly:
\begin{equation*}
    \mathcal{L}_d(u) \;\le\; \Phi_d(\rho,u)\,e^{-u^2/2}\,,
\end{equation*}
where $\Phi_d$ is the polynomial-rational function defined in~\eqref{eq:phi_def}.

\noindent
$\bullet$
For the squared component, we substitute the bound $H_j(\rho x) \le (2\rho x)^j$ of Lemma~\ref{lem:hermite_bound} directly into the integral:
\begin{equation*}
    \mathcal{S}_d(u) \;\le\; \sum_{j=0}^{d-1} c_j^2 \int_u^\infty (2\rho x)^{2j} e^{-(1+\rho^2)x^2/2}\,\mathrm{d}x\,.
\end{equation*}
We evaluate each integral using Lemma~\ref{lem:gauss_tail} with parameters $a=1+\rho^2$ and $m=2j$; its closed form applies since $(2j-1)/u^2<\rho^2<a$ for $\rho u\ge\sqrt{2d-1}$. For the base case $j=0$, we use the Mills bound $I_0(u;a) \le (au)^{-1} e^{-au^2/2}$ of~\eqref{eq:gauss_tail_mills}. Factoring out the shared exponential weight $e^{-(1+\rho^2)u^2/2}$ leaves the rational function $\Psi_d(\rho,u)$ defined in~\eqref{eq:psi_def}, yielding:
\begin{equation*}
    \mathcal{S}_d(u) \;\le\; \Psi_d(\rho,u)\,e^{-(1+\rho^2)u^2/2}\,.
\end{equation*}

\noindent
$\bullet$
Substituting the bounds for $\mathcal{L}_d(u)$ and $\mathcal{S}_d(u)$ back into~\eqref{eq:IMF_split} completes the proof:
\begin{equation*}
    \int_u^\infty Q_d(\rho x)e^{-x^2/2}\,\mathrm{d}x \;\le\; \alpha_d\,\Phi_d(\rho,u)\,e^{-u^2/2} \;+\; \Psi_d(\rho,u)\,e^{-(1+\rho^2)u^2/2}\,. \qedhere
\end{equation*}
\end{proof}

\medskip

\begin{theorem}[IMF tail bound]
\label{thm:imf_tail}
For all $u\ge u_{\mathrm{IMF}}=\sqrt{2d-1}/\rho$,
\begin{equation*}
    \PP\Big\{\sup_{\Sph^{d-1}}X>u\Big\}
    \;\le\;
    \delta_{\mathrm{IMF}}(u)
    \;:=\;
    2(k-1)^{\frac{d-1}2}\Big(\alpha_d\,\Phi_d(\rho,u)\,e^{-u^2/2}+\Psi_d(\rho,u)\,e^{-(1+\rho^2)u^2/2}\Big)\,.
\end{equation*}
\end{theorem}

\begin{proof}[Proof of Theorem~\ref{thm:imf_tail}]
Apply Proposition~\ref{prop:IMF_decomp} to the right-hand side of~\eqref{eq:ortho_exact} and substitute the constant prefactor $2(k-1)^{\frac{d-1}2}$ from Proposition~\ref{prop:ortho_exact}.
\end{proof}

\paragraph{Sharpness of the IMF bound.} Comparing the leading coefficients in Lemmas~\ref{lfyod} and~\ref{lem:dominant_term} with $A(\nu)\sim\nu^{d-1}$ (Lemma~\ref{lem:exp_det}) gives
\begin{equation}\label{eq:alpha_d_closed}
    \alpha_d\;=\;\big(\sqrt2\,2^{d-1}\Gamma(d/2)\big)^{-1};
\end{equation}
hence $\delta_{\mathrm{IMF}}(u)\sim\delta_{\mathrm{bl}}(u)\sim\delta_{\mathrm{exact}}(u)$ as $u\to\infty$, which is the sense in which the bound is asymptotically sharp.

\subsection{Simplified Mehta--Fyodorov bound (SMF)}
\label{sub:SMF}

The SMF bound applies a second relaxation beyond IMF: the envelopes $H_{d-1}(\nu)<(2\nu)^{d-1}$ and $(\nu^2+1)^{d-1}$ collapse $\Phi_d$ and the squared Hermite sum to a single monomial main term $u^{d-2}e^{-u^2/2}$ plus a faster $e^{-3u^2/4}$ remainder, at the cost of a uniform factor $2$ in the leading constant. Beyond the level $u^\star_d$ the remainder is dominated and the bound collapses to a single monomial (Corollary~\ref{cor:single_term}).

\begin{lemma}[Hermite envelope]
\label{lem:hermite_bound}
For every $m\ge 2$ and every $x>x_1^{(m)}$, $0<H_m(x)<(2x)^m$.
\end{lemma}

\begin{proof}[Proof of Lemma~\ref{lem:hermite_bound}]
Write $H_m(x)=2^m\prod_{j=1}^m(x-x_j^{(m)})$. For $x>x_1^{(m)}$, every factor is strictly positive, so $H_m(x)>0$. For the upper bound, the inequality $H_m(x)<(2x)^m$ is equivalent to $\prod_{j=1}^m(1-x_j^{(m)}/x)<1$. Taking logarithms,
\[
    \sum_{j=1}^m\log(1-x_j^{(m)}/x)\;<\;-\sum_{j=1}^m x_j^{(m)}/x\;=\;0\,,
\]
where the strict inequality uses $\log(1-t)<-t$ for $t\neq 0$ and the equality uses the vanishing sum of Hermite roots, which follows from the absence of the $x^{m-1}$ term in $H_m(x)$.
\end{proof}

\begin{lemma}[Non-asymptotic decomposition of $Q_d$]
\label{lem:dominant_term}
With $c_j$ as in Lemma~\ref{lem:mehta}, $S_d:=\sum_{j=0}^{d-1}c_j^2\,2^j\big[\frac{(2j)!}{j!}\big]^2$, $\widetilde B_d:=\frac{(2d)!\,2^{d+1}}{d!}$ and
\begin{subequations}\label{eq:Qd_dominant}
\begin{equation}\label{eq:beta_d}
    \beta_d := S_d + \sqrt{d/2}\,c_{d-1}c_d\,\frac{(2d-2)!}{(d-1)!}\,2^{d-1}\,\widetilde B_d\,,
\end{equation}
\begin{equation}\label{eq:Qd_decomp}
    Q_d(\nu) \;=\; \alpha_d\,H_{d-1}(\nu)\;+\;\mathcal{R}_d(\nu)\,,
\end{equation}
with the dominant coefficient given by
\begin{equation}\label{eq:alpha_d}
    \alpha_d \;=\;
    \begin{cases}
        \tfrac12\sqrt{d/2}\,c_{d-1}c_d\,\mu_d & \text{if $d$ even},\\[2pt]
        1/\mu_{d-1} & \text{if $d$ odd},
    \end{cases}
\end{equation}
and the remainder bound, for every $\nu\ge 0$,
\begin{equation}\label{eq:remainder_bound}
    |\mathcal{R}_d(\nu)|\;\le\;\beta_d\,(\nu^2+1)^{d-1}\,e^{-\nu^2/2}\,.
\end{equation}
\end{subequations}
\end{lemma}

\begin{proof}
Deferred to Appendix~\ref{app:deferred}.
\end{proof}

The decomposition~\eqref{eq:Qd_decomp} translates into a non-asymptotic upper bound on the Kac--Rice integrand:
\[
    A(\nu)
    \;\le\;\alpha_d\,\frac{2^{\frac{3}2}\Gamma((d+2)/2)}{d}\,(2\nu)^{d-1}
    \;+\;\frac{2^{\frac{3}2}\Gamma((d+2)/2)}{d}\,\beta_d\,(\nu^2+1)^{d-1}\,e^{-\nu^2/2}\,.
\]

\medskip

\noindent
Inserting this inequality into the Kac--Rice integral~\eqref{eq:dexact} and applying the Gaussian-type tail integrals of Lemma~\ref{lem:gauss_tail} (Appendix~\ref{app:gaussian_integrals}) yields the SMF bound.

\begin{theorem}[SMF tail bound]
\label{thm:smf_tail}
Let $u_{\mathrm{SMF}}=2\sqrt{d}$. For all $u\ge u_{\mathrm{SMF}}$,
\begin{equation}\label{eq:smf_bound}
    \PP\Big\{\sup_{\Sph^{d-1}}X>u\Big\}
    \;\le\;
    \delta_{\mathrm{SMF}}(u)
    \;:=\;
    4\,\alpha_d\,(2k)^{\frac{d-1}2}\,u^{d-2}\,e^{-u^2/2}\;+\;2^{d}\,\beta_d\,(k-1)^{\frac{d-1}2}\,u^{2d-3}\,e^{-3u^2/4}\,.
\end{equation}
\end{theorem}

\begin{proof}[Proof of Theorem~\ref{thm:smf_tail}]

\noindent
$\bullet$
By Proposition~\ref{prop:ortho_exact}, the excursion probability is bounded by the integral of the Mehta function:
\[
    \PP\Big\{\sup_{\Sph^{d-1}}X>u\Big\} \;\le\; 2(k-1)^{\frac{d-1}2}\int_u^\infty Q_d(\rho x)\,e^{-x^2/2}\,\mathrm{d}x\,.
\]
Using the non-asymptotic decomposition from Lemma~\ref{lem:dominant_term}, we write $Q_d = \alpha_d H_{d-1} + \mathcal{R}_d$. This splits the corresponding integral into two parts: $I_{\mathrm{main}} + I_{\mathrm{rem}}$.

\noindent
$\bullet$
The main term integral evaluates the dominant Hermite polynomial. By hypothesis, $u \ge u_{\mathrm{SMF}} = 2\sqrt{d}$ and $\rho^2 \ge 1/2$ (since $k \ge 3$). This ensures that $\rho u \ge \sqrt{2d} > \sqrt{2d-1}$, meaning $\rho x > x_1^{(d-1)}$ for all integration variables $x \ge u$. Applying the polynomial envelope from Lemma~\ref{lem:hermite_bound}, we have $0 < H_{d-1}(\rho x) < (2\rho x)^{d-1}$. 
Substituting this strictly positive bound into the main integral gives:
\[
    I_{\mathrm{main}} \;<\; \alpha_d(2\rho)^{d-1}\int_u^\infty x^{d-1}e^{-x^2/2}\,\mathrm{d}x\,.
\]
We evaluate this using Lemma~\ref{lem:gauss_tail} with $a=1$ and $m=d-1$. The required tail hypothesis $(d-2)/u^2 \le 1/2$ holds safely for $u \ge 2\sqrt{d}$, yielding:
\[
    I_{\mathrm{main}} \;\le\; 2\alpha_d(2\rho)^{d-1}u^{d-2}e^{-u^2/2}\,.
\]

\noindent
$\bullet$
Using the remainder envelope~\eqref{eq:remainder_bound}, the absolute value of the remainder integral is bounded by:
\[
    |I_{\mathrm{rem}}| \;\le\; \beta_d\int_u^\infty (\rho^2 x^2+1)^{d-1}e^{-(1+\rho^2)x^2/2}\,\mathrm{d}x\,.
\]
For the polynomial factor, we factor $\rho^2 x^2$ out exactly:
\[
    (\rho^2 x^2+1)^{d-1} \;=\; (\rho^2 x^2)^{d-1}\,\big(1+1/(\rho^2 x^2)\big)^{d-1}
    \;=\; \rho^{2(d-1)}\,x^{2(d-1)}\,\big(1+1/(\rho^2 x^2)\big)^{d-1}\,.
\]
Since $\rho^2 \ge 1/2$ (because $k \ge 3$) and $x \ge u \ge 2\sqrt{d}$, we have $\rho^2 x^2 \ge 2d \ge 1$, so the correction factor satisfies $(1+1/(\rho^2 x^2))^{d-1} \le 2^{d-1}$. Combined with $\rho \le 1$, which gives $\rho^{2(d-1)} \le 1$, this yields $(\rho^2 x^2+1)^{d-1} \le 2^{d-1}\,x^{2(d-1)}$.
For the exponential factor, since $k \ge 3$, we know $1+\rho^2 \ge 3/2$. Substituting these simplifications gives:
\[
    |I_{\mathrm{rem}}| \;\le\; 2^{d-1}\beta_d\int_u^\infty x^{2(d-1)}e^{-3x^2/4}\,\mathrm{d}x\,.
\]
Applying Lemma~\ref{lem:gauss_tail} with $a=3/2$ and $m=2(d-1)$ (its second particular case, valid for $u\ge2\sqrt d$), we obtain:
\[
    |I_{\mathrm{rem}}| \;\le\; 2^{d-1}\beta_d\,u^{2d-3}e^{-3u^2/4}\,.
\]

\noindent
$\bullet$
We combine the bounded integrals and multiply by the global prefactor $2(k-1)^{\frac{d-1}2}$. 
For the main term, the identity $(2\rho)^{d-1}(k-1)^{\frac{d-1}2} = (2k)^{\frac{d-1}2}$ consolidates the constants.
Summing the two appropriately scaled components yields the final bound:
\[
    \PP\Big\{\sup_{\Sph^{d-1}}X>u\Big\} \;\le\; 4\alpha_d(2k)^{\frac{d-1}2}u^{d-2}e^{-u^2/2} \;+\; 2^d\beta_d(k-1)^{\frac{d-1}2}u^{2d-3}e^{-3u^2/4}\,. \qedhere
\]
\end{proof}

\medskip

The two-term structure of $\delta_{\mathrm{SMF}}(u)$ reflects the two exponential scales inherited from the Mehta expansion: the $e^{-u^2/2}$ scale of the linear Hermite tail, which matches the optimal Kac--Rice decay rate, and the $e^{-3u^2/4}$ scale of the squared Hermite tail relaxed through the $(\nu^2+1)^{d-1}$ envelope. Since the latter is strictly faster, the remainder is asymptotically negligible and the bound collapses to its main term for sufficiently large $u$. The following corollary pinpoints the explicit threshold $u^\star_d\ge u_{\mathrm{SMF}}$ beyond which this collapse occurs, at the mild price of a factor of~$2$ in the leading prefactor.

\begin{corollary}[Asymptotic single-term form]
\label{cor:single_term}
There exists an explicit threshold $u^\star_d\ge u_{\mathrm{SMF}}$, depending only on $(k,d)$, such that for every $u\ge u^\star_d$ the remainder term in~\eqref{eq:smf_bound} is dominated by the main term, in which case
\begin{equation*}
    \PP\Big\{\sup_{\Sph^{d-1}}X>u\Big\}
    \;\le\;
    8\,\alpha_d\,(2k)^{\frac{d-1}2}\,u^{d-2}\,e^{-u^2/2}\,.
\end{equation*}
Explicitly, $u^\star_d$ is the smallest $u\ge 2\sqrt d$ satisfying $2^{d-2}\beta_d\,u^{d-1}\,e^{-u^2/4}\le\alpha_d(2\rho)^{d-1}$.
\end{corollary}

\begin{proof}[Proof of Corollary~\ref{cor:single_term}]
Since $\big((k-1)/(2k)\big)^{(d-1)/2}=(2\rho)^{-(d-1)}$, the ratio of the remainder term of~\eqref{eq:smf_bound} to its main term is
$R(u)=2^{d-2}\beta_d\,u^{d-1}e^{-u^2/4}/\big(\alpha_d(2\rho)^{d-1}\big)$, which is decreasing on $[\sqrt{2(d-1)},\infty)\supset[2\sqrt d,\infty)$. Hence $R(u)\le1$ for every $u\ge u^\star_d$, and there the bound of Theorem~\ref{thm:smf_tail} is at most twice its main term, $8\alpha_d(2k)^{(d-1)/2}u^{d-2}e^{-u^2/2}$.
\end{proof}

\begin{proposition}[Domination of the SMF bound by the IMF bound]
\label{thm:uniform_domination}
For every $k\ge 3$ and $d\ge 3$,
\begin{equation}\label{eq:imf_le_smf}
    \delta_{\mathrm{IMF}}(u)\;\le\;\delta_{\mathrm{SMF}}(u)\qquad\text{for every }u\ge u_{\mathrm{SMF}}=2\sqrt d\,.
\end{equation}
\end{proposition}

The proof uses the following envelope of $\Phi_d$.

\begin{lemma}[Geometric envelope on $\Phi_d$]
\label{lem:Phi_d_envelope}
For every $k\ge 3$, $d\ge 3$, and every $u\ge u_{\mathrm{SMF}}=2\sqrt d$,
\begin{equation}\label{eq:Phi_d_envelope}
    \Phi_d(\rho,u)\;\le\;2\,(2\rho)^{d-1}\,u^{d-2}\,.
\end{equation}
\end{lemma}

\begin{proof}
Since $\Lambda=2\rho^2-1=1/(k-1)$ (as in~\eqref{eq:phi_def}), we have $2\Lambda=2/(k-1)$, and we set
\[
    q\;:=\;\frac{2\Lambda}{(2\rho u)^2}\;=\;\frac{2}{(k-1)\,(2\rho u)^2}\,.
\]
On $\{u\ge 2\sqrt d\}$, $\rho^2\ge 1/2$ for $k\ge 3$, hence $(2\rho u)^2\ge 4\rho^2\cdot 4d\ge 8d$, and therefore
\[
    (d-2)\,q\;\le\;\frac{2(d-2)}{8d\,(k-1)}\;=\;\frac{d-2}{4d\,(k-1)}\;\le\;\frac1{4(k-1)}\;\le\;\frac1{8}\,.
\]
The leading term of $\Phi_d$ (the $\ell=0$ summand of~\eqref{eq:phi_def}) equals $2\rho\,(2\rho u)^{d-2}=(2\rho)^{d-1}u^{d-2}$. The $\ell$-th summand for $\ell\ge 1$, divided by the leading term, equals $(d-2)(d-4)\cdots(d-2\ell)\,q^\ell$, a product of $\ell$ factors each at most $(d-2)q\le 1/8$; hence
\[
    \sum_{\ell\ge 1}\frac{\text{$\ell$-th summand}}{\text{leading}}\;\le\;\sum_{\ell\ge 1}(1/8)^\ell\;=\;\frac{1}{7}\,.
\]
Similarly, for $d$ odd, the trailing term $\1\{d\,\mathrm{odd}\}\,u^{-1}\,(2\Lambda)^{(d-1)/2}\,(d-2)!!$ divided by the leading is $q^{(d-1)/2}\,(d-2)!!=\prod_{j=1}^{(d-1)/2}(d-2j)\,q\le(1/8)^{(d-1)/2}\le 1/8$ (the last bound uses $(d-1)/2\ge 1$ for $d\ge 3$). Summing,
\[
    \Phi_d(\rho,u)/\bigl((2\rho)^{d-1}u^{d-2}\bigr)\;\le\;1+\tfrac17+\tfrac18\;=\;\tfrac{71}{56}\;\le\;2\,,
\]
which is~\eqref{eq:Phi_d_envelope}.
\end{proof}

\begin{proof}[Proof of Proposition~\ref{thm:uniform_domination}]
We bound the IMF main and residual terms separately on $[u_{\mathrm{SMF}},\infty)$.

\medskip

\noindent
$\bullet$ \emph{Main term.} On the interval $[u_{\mathrm{SMF}},\infty)$, Lemma~\ref{lem:Phi_d_envelope} gives $\Phi_d(\rho,u)\le 2(2\rho)^{d-1}u^{d-2}$. Multiply by $2(k-1)^{(d-1)/2}\alpha_d e^{-u^2/2}$ and use $\rho^{d-1}(k-1)^{(d-1)/2}=(k/2)^{(d-1)/2}$:
\begin{align*}
    \underbrace{2(k-1)^{\frac{d-1}2}\,\alpha_d\,\Phi_d(\rho,u)\,e^{-u^2/2}}_{\text{IMF main}}
    &\;\le\;4\,\alpha_d\,(2\rho)^{d-1}\,(k-1)^{\frac{d-1}2}\,u^{d-2}\,e^{-u^2/2}\\
    &\;=\;4\,\alpha_d\,(2k)^{\frac{d-1}2}\,u^{d-2}\,e^{-u^2/2}
    \;=\;\underbrace{\delta_{\mathrm{SMF}}^{\mathrm{main}}(u)}_{\text{first term of \eqref{eq:smf_bound}}}\,.
\end{align*}

\noindent
$\bullet$ \emph{Residual.} On $[u_{\mathrm{SMF}},\infty)$, $u^2\ge 4d$ so $(2j-1)/u^2\le(2d-3)/(4d)\le 1/2$ for every $j\le d-1$; combined with $1+\rho^2\ge 3/2$ (since $k\ge 3$), the rational factor in~\eqref{eq:psi_def} satisfies
\[
    \frac{1}{1+\rho^2-(2j-1)/u^2}\;\le\;\frac{1}{1+\rho^2-1/2}\;\le\;1\,.
\]
Hence
\[
    \Psi_d(\rho,u)\;\le\;\frac{c_0^2}{(1+\rho^2)\,u}+\sum_{j=1}^{d-1}c_j^2\,(2\rho)^{2j}\,u^{2j-1}
    \;\le\;u^{2d-3}\sum_{j=0}^{d-1}c_j^2\,(2\rho)^{2j}\,,
\]
using $u^{2j-1}\le u^{2d-3}$ for $j\le d-1$ and $u\ge 1$, the trivial bound $1/[(1+\rho^2)u]\le u^{2d-3}$, and $1/(1+\rho^2)\le 1\le(2\rho)^0=1$ to consolidate the $j=0$ term. Comparing the per-$j$ Szeg\H{o}-style coefficient $c_j^2(2\rho)^{2j}=c_j^2\,4^j\rho^{2j}$ to the entry $c_j^2\,2^j\,((2j)!/j!)^2$ of $S_d$ in the definition~\eqref{eq:beta_d} of $\beta_d$: since $(2j)!/j!=(j+1)(j+2)\cdots(2j)\ge 2^j$ for $j\ge 1$, $((2j)!/j!)^2\ge 4^j$, and using $\rho\le 1$, $(2\rho)^{2j}=4^j\rho^{2j}\le 4^j\le 2^j\cdot 4^j\le 2^j((2j)!/j!)^2$, hence~$c_j^2(2\rho)^{2j}\le c_j^2\,2^j((2j)!/j!)^2$ termwise (for $j=0$ the comparison reads $c_0^2\le c_0^2$ with equality), giving $\sum_j c_j^2(2\rho)^{2j}\le S_d\le\beta_d$. Therefore
\[
    \underbrace{2(k-1)^{\frac{d-1}2}\,\Psi_d(\rho,u)\,e^{-(1+\rho^2)u^2/2}}_{\text{IMF residual}}
    \;\le\;2(k-1)^{\frac{d-1}2}\beta_d\,u^{2d-3}\,e^{-(1+\rho^2)u^2/2}\,.
\]
Comparing to $\delta_{\mathrm{SMF}}^{\mathrm{rem}}(u):=2^d\beta_d(k-1)^{(d-1)/2}u^{2d-3}e^{-3u^2/4}$, the ratio of exponentials is $e^{-(1+\rho^2)u^2/2+3u^2/4}=e^{-\Lambda u^2/4}$, which lies in $(0,1]$ for $\Lambda=2\rho^2-1=1/(k-1)>0$ (i.e., $k\ge 3$). The ratio of polynomial prefactors is $2/2^d=2^{1-d}\le 1$ for $d\ge 1$. Hence the IMF residual is bounded by $\delta_{\mathrm{SMF}}^{\mathrm{rem}}$ on $[u_{\mathrm{SMF}},\infty)$.

\medskip

\noindent
Adding the two pieces, $\delta_{\mathrm{IMF}}(u)\le\delta_{\mathrm{SMF}}^{\mathrm{main}}(u)+\delta_{\mathrm{SMF}}^{\mathrm{rem}}(u)=\delta_{\mathrm{SMF}}(u)$ on $[u_{\mathrm{SMF}},\infty)$, which is~\eqref{eq:imf_le_smf}.
\end{proof}

\section{The gap as a marked Rice integral, and the layer decomposition}\label{sec:gap}

\begin{lemma}\label{lem:beaten}
For $u>0$, $\PP\{\sup X>u\}\le\dex(u)$ and, for $k$ odd,
\begin{equation}\label{eq:gap_count}
\dex(u)-\PP\{\sup X>u\}=\E\,\#\big\{t:\ \nabla X(t)=0,\ X(t)\ge u,\ \exists s:\ X(s)>X(t)\big\}.
\end{equation}
\end{lemma}

\begin{proof}
The first claim is~\eqref{eq:KR_main_bound_intro}. For $k$ odd, $\Var(X(s)-X(t))=2(1-\langle s,t\rangle^{k})$ vanishes only for $s=t$, so the maximum is almost surely attained at a single point $t_{*}$, a consequence of Tsirelson's theorem \citep[Theorem~3]{lifshits1983absolute}. A critical point above $u$ is either $t_{*}$, exactly on $\{\sup X\ge u\}$, or beaten; hence $\Ncp_{[u,\infty)}=\1\{\sup X\ge u\}+\#\{\text{beaten critical points above }u\}$, and $\PP\{\sup X=u\}\le\E\Ncp_{\{u\}}=0$ by \eqref{eq:dexact}, so $\PP\{\sup X\ge u\}=\PP\{\sup X>u\}$; both sides of \eqref{eq:gap_count} are finite by \eqref{eq:dexact}.
\end{proof}

Fix $t\in\Sph^{d-1}$ and write $s=\cos\theta\,t+\sin\theta\,v$ with $\theta\in(0,\pi)$, $v\in\Sph^{d-2}$ the unit sphere of $t^{\perp}$, $r=\cos\theta=\langle s,t\rangle$, and $M:=\Sph^{d-1}\setminus\{\pm t\}$. Standard computations give $\Cov(X(s),X(t))=r^{k}$, $\Cov(X(s),\partial_{a}X(t))=k\,r^{k-1}\sin\theta\,v_{a}$, $\Cov(X(t),\nabla X(t))=0$ and $\Cov(\nabla X(t))=kI_{d-1}$, so the Gaussian regression of $X(s)$ on $(X(t),\nabla X(t))$ reads
\begin{subequations}\label{eq:helix}
\begin{equation}\label{eq:regression}
X(s)=r^{k}X(t)+r^{k-1}\sin\theta\,\langle v,\nabla X(t)\rangle+R^{t}(s),
\end{equation}
with a residual $R^{t}$ independent of $(X(t),\nabla X(t))$. The \emph{helix process} is
\begin{equation}\label{eq:helixdef}
\zeta^{t}(s):=\frac{R^{t}(s)}{1-r^{k}},\qquad s\in M ,
\end{equation}
well defined since $1-r^{k}>0$ on $M$. Substituting $X(t)=x>0$ and $\nabla X(t)=0$ in \eqref{eq:regression} gives $X(s)=r^{k}x+(1-r^{k})\zeta^{t}(s)$, so $X(s)>x$ if and only if $\zeta^{t}(s)>x$, while the excluded points never beat $t$ ($X(-t)=(-1)^{k}x\le x$). Hence, when $(X(t),\nabla X(t))$ is replaced by $(x,0)$ in \eqref{eq:regression}, $x>0$,
\begin{equation}\label{eq:beating}
\big\{\exists s\in\Sph^{d-1}:\ X(s)>x\big\}=\Big\{\sup_{M}\zeta^{t}>x\Big\};
\end{equation}
\end{subequations}
in this form the beating event enters the Rice formula below, since $(G,\zeta^{t})$ is independent of $(X(t),\nabla X(t))$.

\begin{proposition}[The gap as a marked Rice integral]\label{prop:mult}
Let $k$ be odd and $u>0$. Then
\begin{equation}\label{eq:Xi}
\dex(u)-\PP\{\sup X>u\}=C_{k,d}\int_{u}^{\infty}\Xi(x)\,\varphi(x)\,\mathrm{d}x ,
\end{equation}
with $\Xi$ as in \eqref{eq:Xidef}. In particular $0\le\Xi(x)\le A(\rho x)$ and $\dex(u)-\PP\{\sup X>u\}\le\varepsilon(u)\,\dex(u)$ with $\varepsilon(u)=\sup_{x\ge u}\Xi(x)/A(\rho x)$.
\end{proposition}

\begin{proof}
We use the Kac--Rice formula with a lower semi-continuous weight \citep[Theorem~6.1]{armentano2023general}, which, for a Gaussian field, requires only $C^{1}$ paths and a continuous density, bounded on compacts; it counts all zeros, degenerate ones included. We apply it on the charts of a finite atlas of $\Sph^{d-1}$ to the $C^{\infty}$ Gaussian field $\nabla X$, with $\nabla X(t)\sim N(0,kI_{d-1})$, and to the weight $g(t,X(\cdot))=\1\{X(t)>u\}\1\{\sup_{s}X(s)>X(t)\}$, the indicator of an open condition on the continuous functions $(t,X(\cdot))\mapsto X(t)$ and $\sup_{s}X(s)-X(t)$, hence lower semi-continuous in $t$ and in $X(\cdot)$ for uniform convergence; the integrand $|\det\nabla^{2}X|\,p_{\nabla X}(0)\,\mathrm{d}t$ is chart independent at critical points. This writes the expected number of beaten critical points above $u$ as $\int_{\Sph^{d-1}}\E[|\det\nabla^{2}X(t)|\1\{X(t)>u\}\1\{\sup_{s}X(s)>X(t)\}\mid\nabla X(t)=0]\,p_{\nabla X(t)}(0)\,\mathrm{d}t$, where $X(t)>u$ may replace $X(t)\ge u$ in \eqref{eq:gap_count} since $\E\Ncp_{\{u\}}=0$ by \eqref{eq:dexact}. Condition further on $X(t)=x$, of density $\varphi$, independent of $\nabla X(t)$. By \eqref{eq:condhess}, $|\det\nabla^{2}X(t)|=(2k(k-1))^{(d-1)/2}|\det(\rho xI-G)|$; by \eqref{eq:beating} the beating indicator is $\1\{\sup_{M}\zeta^{t}>x\}$; and $(G,\zeta^{t})$ is independent of $(X(t),\nabla X(t))$, so the conditional expectation is an unconditional one over $(G,\zeta^{t})$. Conditioning on $G$ inside gives $\Xi(x)$, and $p_{\nabla X(t)}(0)=(2\pi k)^{-(d-1)/2}$ gives $C_{k,d}$. The bounds follow from $0\le\PP(\cdot\mid G)\le1$ and \eqref{eq:dexact}.
\end{proof}

\begin{corollary}[H\"older decoupling]\label{cor:cs}
For every $x>0$ and $p>1$, with $P(x):=\PP\{\sup_{M}\zeta^{t}>x\}$,
\begin{equation}\label{eq:Xi_cs}
\Xi(x)\;\le\;\big\|\det(\rho xI-G)\big\|_{p}\;P(x)^{1-1/p},\qquad\text{and}\qquad P(x)\ge\bar\Phi(\rho x).
\end{equation}
\end{corollary}

\begin{proof}
$\Xi(x)=\E[|\det(\rho xI-G)|\,\1\{\sup_{M}\zeta^{t}>x\}]$, the mark in \eqref{eq:Xidef} being a conditional probability given $G$; apply H\"older's inequality. The lower bound is $\sup_{M}\zeta^{t}\ge\varsigma\lmax(G)\ge\varsigma G_{11}$ (Lemma~\ref{lem:poles}) and $\varsigma\rho=1$.
\end{proof}

The first factor is an $L^{p}$ norm of the GOE (Lemma~\ref{lem:goe}). The second is the tail of $\sup\zeta^{t}$, a field of variance at most $\varsigma^{2}$ (Lemma~\ref{lem:variance}); by the lower bound in \eqref{eq:Xi_cs} it decays at best like $e^{-x^{2}/(2\varsigma^{2})}=e^{-\rho^{2}x^{2}/2}$, so the rate of Theorem~\ref{thm:target} is the best this route can give, and with $p=\rho^{2}x^{2}$ the power $1-1/p$ loses only a constant factor. The Rice formula gives the tail without the mean shift $\E\sup\zeta^{t}\approx\varsigma\sqrt{2(d-1)}$ that Borell--TIS would put into the exponent.

\subsection{The layer decomposition}\label{sec:layers}

Let $\Mb:=[0,\pi]\times\Sph^{d-2}$ be the compact blow-up of the sphere at $\pm t$, $M\cong(0,\pi)\times\Sph^{d-2}$; the sets $\{0\}\times\Sph^{d-2}$ and $\{\pi\}\times\Sph^{d-2}$ are the \emph{pole fibers}.

\begin{proposition}[Layers]\label{prop:layers}
There exist independent centered Gaussian fields $Z_{2},\dots,Z_{k}$ on $\Sph^{d-2}$, $Z_{j}$ a KSS field of degree $j$ (covariance $\langle v,v'\rangle^{j}$), such that
\begin{subequations}\label{eq:layer_decomp}
\begin{equation}\label{eq:layers}
\zeta^{t}(s)=\sum_{j=2}^{k}w_{j}(\theta)\,Z_{j}(v),
\qquad
w_{j}(\theta)=\binom kj^{1/2}\frac{r^{\,k-j}(1-r^{2})^{j/2}}{1-r^{k}},\qquad r=\cos\theta .
\end{equation}
Write $w(\theta):=(w_{j}(\theta))_{2\le j\le k}$. The $w_{j}$ are linearly independent on $(0,\pi)$ and $Z_{j}(v)=L_{j}[\zeta^{t}(\cdot,v)]$ for the explicit linear functionals
\begin{equation}\label{eq:dualbasis}
L_{j}[f]:=\sum_{i=2}^{k}(\Gamma^{-1})_{ji}\int_{0}^{\pi}w_{i}(\theta)f(\theta)\,\mathrm{d}\theta,\qquad \Gamma=\Big(\int_{0}^{\pi}w_{i}w_{j}\,\mathrm{d}\theta\Big)_{2\le i,j\le k},\qquad L_{j}[w_{i}]=\delta_{ij}.
\end{equation}
\end{subequations}
\end{proposition}

\begin{proof}
Let $s_{i}=r_{i}t+\sqrt{1-r_{i}^{2}}\,v_{i}$, $c=\langle v_{1},v_{2}\rangle$, $q=\sqrt{(1-r_{1}^{2})(1-r_{2}^{2})}$. By the binomial formula $\langle s_{1},s_{2}\rangle^{k}=\sum_{j=0}^{k}\binom kj(r_{1}r_{2})^{k-j}q^{j}c^{j}$, and by \eqref{eq:regression} the covariance of the regression parts is $r_{1}^{k}r_{2}^{k}+k(r_{1}r_{2})^{k-1}qc$, exactly the terms $j=0,1$. Hence $\E[R^{t}(s_{1})R^{t}(s_{2})]=\sum_{j\ge2}\binom kj(r_{1}r_{2})^{k-j}q^{j}c^{j}$, which divided by $(1-r_{1}^{k})(1-r_{2}^{k})$ is the covariance of the right side of \eqref{eq:layers}; a centered Gaussian field is determined by its covariance. Linear independence: if $\sum_{j}a_{j}w_{j}\equiv0$, multiplying by $1-\cos^{k}\theta$ and dividing by $\cos^{k}\theta$ on $(0,\pi/2)$ gives $\sum_{j}a_{j}\binom kj^{1/2}\tan^{j}\theta\equiv0$, a polynomial in $\tan\theta$ vanishing on $(0,\infty)$, so all $a_{j}=0$. The $w_{j}$ are bounded on $(0,\pi)$ (Lemma~\ref{lem:poles}(i), which concerns the deterministic $w_{j}$ only), so the Gram matrix $\Gamma$ in $L^{2}(0,\pi)$ is positive definite and \eqref{eq:dualbasis} satisfies $L_{j}[w_{i}]=\delta_{ij}$. Define $Z_{j}(v):=L_{j}[\zeta^{t}(\cdot,v)]$: the integrals exist almost surely and commute with covariances, since $\E\int_{0}^{\pi}|w_{i}(\theta)\zeta^{t}(\theta,v)|\,\mathrm{d}\theta\le\int_{0}^{\pi}|w_{i}|\sqrt{V}\,\mathrm{d}\theta<\infty$, where $V=\Var\zeta^{t}=\sum_{l}w_{l}^{2}$ is bounded by the covariance identity above. From $\E[\zeta^{t}(\theta,v)\zeta^{t}(\theta',v')]=\sum_{l}w_{l}(\theta)w_{l}(\theta')\langle v,v'\rangle^{l}$ and $L_{j}[w_{l}]=\delta_{jl}$ one gets $\Cov(Z_{i}(v),Z_{j}(v'))=\delta_{ij}\langle v,v'\rangle^{i}$ and $\Cov(\zeta^{t}(s),Z_{j}(v))=w_{j}(\theta)$, hence $\Var\big(\zeta^{t}(s)-\sum_{j}w_{j}(\theta)Z_{j}(v)\big)=V-2\sum_{j}w_{j}^{2}+\sum_{j}w_{j}^{2}=0$. So the $Z_{j}$ are independent centered Gaussian fields with covariance $\langle v,v'\rangle^{j}$, in the Gaussian space of $X$; taking their continuous versions, \eqref{eq:layers} holds almost surely for each $s$, hence for all $s\in M$ by continuity of both sides.
\end{proof}

\begin{lemma}[The layer $j=2$ is the Hessian]\label{lem:global}
Almost surely $Z_{2}(v)=v^{\top}Gv$ for all $v\in\Sph^{d-2}$, with $G$ the matrix \eqref{eq:G} at the same $t$; the layers $Z_{3},\dots,Z_{k}$ are independent of $G$.
\end{lemma}

\begin{proof}
From $\langle s,t\rangle^{k}$, $\Cov(X(s),\nabla^{2}_{ab}X(t))=k(k-1)r^{k-2}\sin^{2}\theta\,v_{a}v_{b}-kr^{k}\delta_{ab}$, $\Cov(X(t),\nabla^{2}_{ab}X(t))=-k\delta_{ab}$ and $\Cov(\nabla X(t),\nabla^{2}X(t))=0$, so in \eqref{eq:regression} the $\delta_{ab}$ terms cancel: $\Cov(R^{t}(s),\nabla^{2}_{ab}X(t))=k(k-1)r^{k-2}\sin^{2}\theta\,v_{a}v_{b}$. Dividing by $(1-r^{k})\sqrt{2k(k-1)}$ and noting that the $kX(t)I$ part of $G$ is uncorrelated with $R^{t}$, $\Cov(\zeta^{t}(s),G_{ab})=w_{2}(\theta)v_{a}v_{b}$. Applying $L_{j}$ of \eqref{eq:dualbasis} in $\theta$ gives $\Cov(Z_{j}(v),G_{ab})=\delta_{j2}v_{a}v_{b}$; hence $\Cov(Z_{2}(v),v^{\top}Gv)=\sum_{a,b}v_{a}^{2}v_{b}^{2}=1$, and both variables have unit variance ($\Var(v^{\top}Gv)=\sum_{a}v_{a}^{4}+2\sum_{a<b}v_{a}^{2}v_{b}^{2}=1$), so they are almost surely equal, for each $v$ and then, by continuity, for all $v$; the layers $j\ge3$ are uncorrelated with $G$ and, lying with $G$ in the Gaussian space of $X$, independent of it.
\end{proof}

\begin{lemma}[Variance]\label{lem:variance}
$V(\theta):=\Var\zeta^{t}(s)=\|w(\theta)\|^{2}=\dfrac{1-r^{2k}-k\,r^{2k-2}(1-r^{2})}{(1-r^{k})^{2}}\le\varsigma^{2}=\dfrac{2(k-1)}{k}$ on $(0,\pi)$, with equality only in the limit $\theta\to0$ (and $\theta\to\pi$ for $k$ even).
\end{lemma}

\begin{proof}
The closed form is the covariance identity of Proposition~\ref{prop:layers} at $s_{1}=s_{2}$. The numerator is even in $r$ and $(1-r^{k})^{2}\ge(1-|r|^{k})^{2}$, so it suffices to take $r\in[0,1)$, where the claim reads $\psi(r):=2(k-1)(1-r^{k})^{2}-k[1-r^{2k}-kr^{2k-2}(1-r^{2})]\ge0$. Expanding, $\psi=(k-2)(1-r^{k})^{2}+kr^{k}\chi(r)$ with $\chi(r)=kr^{k-2}-(k-2)r^{k}-2$,~$\chi(1)=0$, $\chi'(\varrho)=k(k-2)\varrho^{k-3}(1-\varrho^{2})\le2k(k-2)(1-\varrho)$ on $[0,1]$, so $-\chi(r)=\int_{r}^{1}\chi'\le k(k-2)(1-r)^{2}$. The arithmetic--geometric mean inequality gives $1-r^{k}=(1-r)\sum_{i<k}r^{i}\ge kr^{(k-1)/2}(1-r)$, hence $kr^{k}(-\chi)\le k^{2}(k-2)r^{k}(1-r)^{2}\le(k-2)(1-r^{k})^{2}$ and $\psi\ge0$, strictly for $r\in(0,1)$; the expansion of Lemma~\ref{lem:poles} gives $V\to\varsigma^{2}$ at the pole.
\end{proof}

By~\eqref{eq:regression}, $V(\theta)$ equals $\Var(X(s)\mid X(t),\nabla X(t))/(1-r^{k})^{2}$, so $\varsigma^{2}=\sup_{\theta}V$ is the constant $\sigma_{t}^{2}$ of \citet[Theorem~8.12]{azais2009level} for the KSS field with derivatives along the sphere, where their $\kappa_{t}$ vanishes because the cone of $\Sph^{d-1}$ at $t$ is the whole tangent space, the sphere having no boundary; their exponent $1+1/\sigma_{t}^{2}$ is $1+\rho^{2}$.

\begin{lemma}[Poles and parities]\label{lem:poles}
(i) Each $w_{j}$ is $C^{\infty}$ on $[0,\pi)$, and as $\theta\to0$
\begin{gather*}
w_{j}(\theta)=\tfrac2k\tbinom kj^{1/2}\theta^{\,j-2}\big(1+O(\theta^{2})\big),\qquad
w_{2}(\theta)=\varsigma\big[1-\tfrac{k-2}{4}\theta^{2}+O(\theta^{4})\big],\\
w_{3}(\theta)=b_{3}\theta+O(\theta^{3}),\qquad b_{3}:=\tfrac2k\tbinom k3^{1/2}>0 .
\end{gather*}
For $k$ odd, $w_{j}(\pi-\eta)=(-1)^{k-j}\tfrac12\binom kj^{1/2}\eta^{j}(1+O(\eta^{2}))$ as $\eta\to0$; for $k$ even, $w_{j}(\pi-\theta)=(-1)^{k-j}w_{j}(\theta)$.
(ii) $\zeta^{t}$ extends to an almost surely continuous field $\bar\zeta(\theta,v)=\sum_{j}w_{j}(\theta)Z_{j}(v)$ on $\Mb$, with $\bar\zeta(0,v)=\varsigma v^{\top}Gv$, and $\bar\zeta(\pi,\cdot)\equiv0$ for $k$ odd, $\bar\zeta(\pi,v)=\varsigma v^{\top}Gv$ for $k$ even. In particular $\sup_{M}\zeta^{t}\ge\varsigma\lmax(G)$.
(iii) Almost surely $Z_{j}(-v)=(-1)^{j}Z_{j}(v)$ for all $v$; for $k$ even, $\zeta^{t}(-s)=\zeta^{t}(s)$.
\end{lemma}

\begin{proof}
(i) Write $1-\cos^{k}\theta=\theta^{2}h(\theta)$ and $\sin\theta=\theta g(\theta)$ with $h,g$ smooth and even, $h(0)=k/2$, $g(0)=1$; then $w_{j}=\binom kj^{1/2}\cos^{k-j}\theta\,g^{j}\theta^{j-2}/h$. The second-order term of $w_{2}$ comes from $\cos^{k-2}\theta=1-\frac{k-2}{2}\theta^{2}+O(\theta^{4})$, $\sin^{2}\theta=\theta^{2}(1-\frac13\theta^{2}+O(\theta^{4}))$ and $1-\cos^{k}\theta=\frac k2\theta^{2}(1-\frac{3k-2}{12}\theta^{2}+O(\theta^{4}))$. At $\theta=\pi-\eta$, $r=-\cos\eta$, $1-r^{k}=1+\cos^{k}\eta\to2$ for $k$ odd, and for $k$ even $r\mapsto-r$ leaves $1-r^{k}$ and $1-r^{2}$ invariant. (ii) The $Z_{j}$ are continuous and the $w_{j}$ extend continuously to $[0,\pi]$ by (i); the boundary values follow from (i) and Lemma~\ref{lem:global}. (iii) $\Var(Z_{j}(-v)-(-1)^{j}Z_{j}(v))=2-2(-1)^{j}\langle-v,v\rangle^{j}=0$, and continuity; for $k$ even, $-s=(\pi-\theta,-v)$ and $w_{j}(\pi-\theta)Z_{j}(-v)=w_{j}(\theta)Z_{j}(v)$.
\end{proof}

\begin{figure}[t]
\centering
\includegraphics[width=0.58\linewidth]{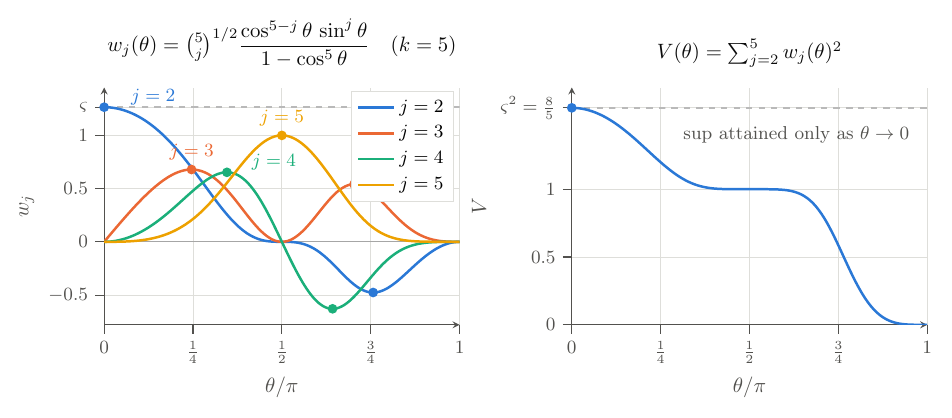}
\caption{Layer weights $w_{j}$ of \eqref{eq:layers} for $k=5$ against $\theta/\pi$ (left) and $V=\sum_{j}w_{j}^{2}$ (right). Only $w_{2}$ is maximal at the pole, where it tends to $\varsigma=\sqrt{8/5}$ (dashed); the higher layers vanish there and peak in the interior; $V$ decreases from $\varsigma^{2}$ to $0$.}
\label{fig:weights}
\end{figure}

\section{Interiority, the Rice bound, and the proof of Theorem~\ref{thm:target}}\label{sec:second}

\subsection{The supremum is attained in the interior}\label{sec:interior}

\begin{proposition}\label{prop:interior}
Let $k\ge3$. Almost surely the supremum of $\bar\zeta$ over $\Mb$ is not attained on a pole fiber. Consequently $\sup_{M}\zeta^{t}$ is almost surely attained at a local maximum of $\zeta^{t}$ on $M$, and for every $x\in\R$
\begin{equation}\label{eq:interior}
\1\Big\{\sup_{M}\zeta^{t}>x\Big\}\;\le\;\Nlm_{(x,\infty)}(\zeta^{t})\quad\text{a.s.},\qquad\text{hence}\qquad
\PP\Big\{\sup_{M}\zeta^{t}>x\Big\}\le\E\,\Nlm_{(x,\infty)}(\zeta^{t}).
\end{equation}
\end{proposition}

\begin{proof}
On the fiber at $t$, $\bar\zeta(0,v)=\varsigma v^{\top}Gv$ has maximum $\varsigma\lmax(G)$, attained exactly at $v=\pm v_{1}$, $v_{1}$ a unit top eigenvector: the top eigenvalue is almost surely simple, the GOE having a density and the matrices with a repeated eigenvalue forming a proper algebraic subvariety. By Lemma~\ref{lem:poles}(i) there are constants $a_{j}$ with $|w_{2}-\varsigma|\le a_{2}\theta^{2}$, $|w_{3}-b_{3}\theta|\le a_{3}\theta^{2}$ and $|w_{j}|\le a_{j}\theta^{2}$ for $j\ge4$ on $[0,\pi]$, so, on an event of probability one, for all $v$,
\begin{equation}\label{eq:meridian}
\bar\zeta(\theta,v)-\varsigma v^{\top}Gv=b_{3}\,\theta\,Z_{3}(v)+\theta^{2}R(\theta,v),\qquad |R(\theta,v)|\le C(v):=1+\sum_{j}a_{j}|Z_{j}(v)|<\infty .
\end{equation}
The layer $Z_{3}$ is odd (Lemma~\ref{lem:poles}(iii)) and independent of $G$, hence of $v_{1}$; given $G$, $Z_{3}(v_{1})\sim N(0,1)$, so $Z_{3}(v_{1})\ne0$ almost surely. Let $v^{*}\in\{\pm v_{1}\}$ with $Z_{3}(v^{*})>0$: by \eqref{eq:meridian}, $\bar\zeta(\theta,v^{*})>\varsigma\lmax(G)=\max_{v}\bar\zeta(0,v)$ for~$0<\theta<b_{3}Z_{3}(v^{*})/C(v^{*})$, so the supremum over $\Mb$ is not attained on the fiber at $t$. For $k$ even the fiber at $-t$ carries the same field and $\bar\zeta(\pi-\theta,-v)=\bar\zeta(\theta,v)$, and the same argument applies. For $k$ odd, $\bar\zeta(\pi,\cdot)\equiv0$ while $\sup\bar\zeta>0$ almost surely: $\lmax(G)\ne0$ almost surely since $G$ has a density; if $\lmax(G)>0$ this is the fiber at $t$; if $\lmax(G)<0$, Lemma~\ref{lem:poles}(i) gives $\bar\zeta(\pi-\eta,v)=-\frac12\binom k2^{1/2}\eta^{2}v^{\top}Gv+O(\eta^{3})>0$ for small $\eta$. The supremum over the compact $\Mb$ is attained, hence at a point of $M$, which is a local maximum of $\zeta^{t}$; on $\{\sup_{M}\zeta^{t}>x\}$ its value exceeds $x$, whence \eqref{eq:interior}.
\end{proof}

\subsection{Kac--Rice for the helix process}\label{sec:kr}

From now on $n=d-2$, $s=(\theta,v)\in M$, $r=\cos\theta$, and $(e_{r},e_{1},\dots,e_{n})$ is the orthonormal frame at $s$ made of the meridian direction $e_{r}=\partial_{\theta}$ and of $n$ directions tangent to the orbit sphere $\{\theta\}\times\Sph^{n}$, of radius $\sin\theta$. Primes are derivatives in $\theta$, $w=(w_{j})_{2\le j\le k}$, and
\begin{equation}\label{eq:scalars}
V=\|w\|^{2},\qquad C=\langle w,w'\rangle=\tfrac12V',\qquad g_{r}=\|w'\|^{2},\qquad \Delta=Vg_{r}-C^{2},\qquad g_{t}=\frac{1}{\sin^{2}\theta}\sum_{j}j\,w_{j}^{2}.
\end{equation}

\begin{lemma}[One-point law of a layer]\label{lem:layerlaw}
Let $Z$ be a KSS field of degree $j\ge2$ on $\Sph^{m}$ and $v\in\Sph^{m}$. Then $Z(v)\sim N(0,1)$, $\nabla Z(v)\sim N(0,jI_{m})$ and $\widehat G:=(\nabla^{2}Z(v)+jZ(v)I_{m})/\sqrt{2j(j-1)}\sim\GOE(m)$ are independent.
\end{lemma}

\begin{proof}
The covariances of Section~\ref{sec:main} with $(k,d)$ replaced by $(j,m+1)$: $\widehat G$ is uncorrelated with $Z(v)$ and $\nabla Z(v)$ and has the GOE covariances.
\end{proof}

\begin{lemma}[Scalars of the helix process]\label{lem:scalars}
For $\theta\in(0,\pi)$: 

\noindent
$(i)$ $(\zeta^{t}(s),\partial_{\theta}\zeta^{t}(s))$ has covariance $\Sigma=\begin{psmallmatrix}V&C\\ C&g_{r}\end{psmallmatrix}$, and the tangential derivatives $\partial_{a}\zeta^{t}(s)$ are independent of it and of each other, each $N(0,g_{t})$. 

\noindent
$(ii)$ $g_{t}=k(1-r^{2k-2})/(1-r^{k})^{2}$ and $\sum_{j}j(j-1)w_{j}^{2}=k(k-1)\sin^{4}\theta/(1-r^{k})^{2}$. 

\noindent
$(iii)$ $\Delta=(\sin\theta\cos\theta)^{-2}\sum_{j<l}(l-j)^{2}w_{j}^{2}w_{l}^{2}$ (by continuity at $\pi/2$) and
\[
\Delta\ge\frac{2^{6-2k}}{k^{4}}\quad\text{if }\cos\theta\ge0\text{ or }k\text{ even},\qquad
\Delta\ge4^{1-k}\sin^{8}\theta\quad\text{if }\cos\theta\le0\text{ and }k\text{ odd};
\]
in particular $\Delta>0$, $g_{r}\ge\Delta/V>0$ and $g_{t}>0$.
\end{lemma}

\begin{proof}
(i) $\zeta^{t}(s)=\sum_{j}w_{j}Z_{j}(v)$, $\partial_{\theta}\zeta^{t}(s)=\sum_{j}w_{j}'Z_{j}(v)$, $\partial_{a}\zeta^{t}(s)=\sin^{-1}\theta\sum_{j}w_{j}\partial_{a}Z_{j}(v)$, and Lemma~\ref{lem:layerlaw}. (ii) The identities $\sum_{j}j\binom kja^{k-j}b^{j}=kb(a+b)^{k-1}$ and $\sum_{j}j(j-1)\binom kja^{k-j}b^{j}=k(k-1)b^{2}(a+b)^{k-2}$ at $a=r^{2}$, $b=\sin^{2}\theta$; the terms $j=0,1$ of the first contribute $kr^{2k-2}\sin^{2}\theta$, those of the second vanish. (iii) Lagrange's identity $\Delta=\sum_{j<l}(w_{j}w_{l}'-w_{l}w_{j}')^{2}$ and the logarithmic derivative of \eqref{eq:layers},
\begin{equation}\label{eq:logder}
\frac{w_{j}'}{w_{j}}=j\cot\theta-(k-j)\tan\theta-g,\qquad g:=\frac{k\,r^{k-1}\sin\theta}{1-r^{k}},\qquad\theta\ne\tfrac{\pi}{2},
\end{equation}
give $w_{j}w_{l}'-w_{l}w_{j}'=(l-j)w_{j}w_{l}/(\sin\theta\cos\theta)$; the identities derived from \eqref{eq:logder} are polynomial in $(\cos\theta,\sin\theta)$ and extend to $\pi/2$ by continuity. Keeping the pairs $(2,3)$ and $(k-1,k)$,
\[
\Delta\ge\max\Big(\tbinom k2\tbinom k3\frac{r^{4k-12}\sin^{8}\theta}{(1-r^{k})^{4}},\ k\,\frac{\sin^{4k-4}\theta}{(1-r^{k})^{4}}\Big).
\]
If $r\ge0$, $1-r^{k}=(1-r)\sum_{i<k}r^{i}\le k\sin^{2}\theta$, so the two terms are at least $r^{4k-12}/k^{4}$ and $\sin^{4k-12}\theta/k^{3}$, and $\max(r^{2},\sin^{2}\theta)\ge\frac12$ gives $\Delta\ge k^{-4}2^{6-2k}$. If $r\le0$ and $k$ is odd, $1-r^{k}\le2$, the terms are at least $r^{4k-12}\sin^{8}\theta/16$ and $\sin^{4k-4}\theta/16$, and $\Delta\ge\frac{\sin^{8}\theta}{16}\max(r^{2},\sin^{2}\theta)^{2k-6}\ge4^{1-k}\sin^{8}\theta$. For $k$ even, $\Delta$ is symmetric about $\pi/2$ by Lemma~\ref{lem:poles}(i).
\end{proof}

\begin{lemma}[Kac--Rice for the local maxima]\label{lem:kr}
Let $B=-\nabla^{2}\zeta^{t}(s)$ and write ``given $(y,0)$'' for the conditioning on $\zeta^{t}(s)=y$, $\nabla\zeta^{t}(s)=0$. For $x>0$,
\begin{equation}\label{eq:KR}
\begin{gathered}
\E\,\Nlm_{(x,\infty)}(\zeta^{t})=\omega_{n}\int_{0}^{\pi}\sin^{n}\theta\int_{x}^{\infty}p_{\theta}(y)\,\E\big[\det B\,\1\{B\succ0\}\,\big|\,(y,0)\big]\,\mathrm{d}y\,\mathrm{d}\theta,\\
p_{\theta}(y)=\frac{\Delta^{-1/2}\,g_{t}^{-n/2}}{(2\pi)^{d/2}}\,e^{-y^{2}g_{r}/(2\Delta)},
\end{gathered}
\end{equation}
with $\omega_{n}=2\pi^{(n+1)/2}/\Gamma(\frac{n+1}{2})$ the volume of $\Sph^{n}$ and $p_{\theta}(y)$ the density of $(\zeta^{t}(s),\nabla\zeta^{t}(s))$ at $(y,0)$.
\end{lemma}

\begin{proof}
The paths of $\zeta^{t}$ are $C^{\infty}$ on $M$. By Lemma~\ref{lem:scalars}, $\nabla\zeta^{t}(s)$ has the nondegenerate covariance $\diag(g_{r},g_{t}I_{n})$, continuous in $s$, so its density is continuous and bounded on compact subsets of $M$, and the Kac--Rice formula with lower semi-continuous weights of the proof of Proposition~\ref{prop:mult} applies to $\nabla\zeta^{t}$ on the charts of $M$, with weights that are indicators of open conditions on $(s,\zeta^{t}(\cdot),\nabla^{2}\zeta^{t}(\cdot))$. First, for a compact $K\subset M$ and $\epsilon>0$, the weight $\1\{|\det\nabla^{2}\zeta^{t}(s)|<\epsilon\}$ gives $\E\#\{s\in K:\nabla\zeta^{t}(s)=0,\ |\det\nabla^{2}\zeta^{t}(s)|<\epsilon\}\le\epsilon\int_{K}p_{\nabla\zeta^{t}(s)}(0)\,\mathrm{d}s\to0$ as $\epsilon\to0$; exhausting $M$ by compacts, almost surely no critical point of $\zeta^{t}$ in $M$ is degenerate, so every local maximum has $B\succ0$. Second, the weight $\1\{\nabla^{2}\zeta^{t}(s)\prec0\}\1\{\zeta^{t}(s)>x\}$, on compacts exhausting $M$, gives $\E\Nlm_{(x,\infty)}(\zeta^{t})=\int_{M}\E[|\det\nabla^{2}\zeta^{t}|\1\{\nabla^{2}\zeta^{t}\prec0\}\1\{\zeta^{t}>x\}\mid\nabla\zeta^{t}=0]\,p_{\nabla\zeta^{t}(s)}(0)\,\mathrm{d}s$. The law of the field around $s$ depends on $s$ only through $\theta$ (the layers $j\ge3$ are isotropic and~$Z_{2}(Rv)=v^{\top}(R^{\top}GR)v$ with $R^{\top}GR\sim\GOE(d-1)$), so $\int_{M}\mathrm{d}s=\omega_{n}\int_{0}^{\pi}\sin^{n}\theta\,\mathrm{d}\theta$; conditioning on $\zeta^{t}(s)=y$ and Lemma~\ref{lem:scalars}(i) give the density, with $\det\Sigma=\Delta$ and $(\Sigma^{-1})_{11}=g_{r}/\Delta$.
\end{proof}

\subsection{The diagonal of the conditional Hessian}\label{sec:hessian}

\begin{lemma}\label{lem:diag}
Given $(y,0)$, the entries $B_{11},\dots,B_{nn},B_{rr}$ of $B=-\nabla^{2}\zeta^{t}(s)$ are independent Gaussian variables,
\[
B_{aa}\sim N\big(\mu_{t}\,y,\ 2\beta_{tt}\big),\qquad \mu_{t}=\frac{k}{1-r^{k}},\qquad \beta_{tt}=\frac{k(k-1)}{(1-r^{k})^{2}},\qquad \frac{\mu_{t}}{\sqrt{2\beta_{tt}}}=\rho,
\]
and $\E[B_{rr}^{2}\mid(y,0)]\le\|w''\|^{2}\big(1+y^{2}g_{r}/\Delta\big)$.
\end{lemma}

\begin{proof}
In the coordinates $(\theta,v)$ the round metric is the warped product $\mathrm{d}\theta^{2}+\sin^{2}\theta\,h$, and for a unit field $\widetilde U$ of $(\Sph^{n},h)$, $\nabla_{\widetilde U}\widetilde U=\nabla^{h}_{\widetilde U}\widetilde U-\sin\theta\cos\theta\,\partial_{\theta}$, $\nabla_{\partial_{\theta}}\partial_{\theta}=0$, $e_{a}=\widetilde U/\sin\theta$. Hence, for $f=\sum_{j}w_{j}Z_{j}$,
\[
(\nabla^{2}f)_{rr}=\sum_{j}w_{j}''Z_{j}(v),\qquad
(\nabla^{2}f)_{aa}=\frac{1}{\sin^{2}\theta}\sum_{j}w_{j}\,(\nabla^{2}_{\Sph^{n}}Z_{j})_{aa}(v)+\cot\theta\sum_{j}w_{j}'Z_{j}(v).
\]
For $f=\cos\theta$ and $f=\sin\theta\,v_{e}$, restrictions of linear forms, whose Hessian on the sphere is $-f$ times the metric, both formulas give $-f$. By Lemma~\ref{lem:layerlaw}, $(\nabla^{2}_{\Sph^{n}}Z_{j})_{aa}=\sqrt{2j(j-1)}\widehat G^{(j)}_{aa}-jZ_{j}(v)$ with independent $\widehat G^{(j)}\sim\GOE(n)$, independent of all $Z_{j}(v)$, $\nabla Z_{j}(v)$. So
\[
B_{aa}=\sum_{j}\Big[\frac{j\,w_{j}}{\sin^{2}\theta}-\cot\theta\,w_{j}'\Big]Z_{j}(v)-\sum_{j}\frac{\sqrt{2j(j-1)}\,w_{j}}{\sin^{2}\theta}\,\widehat G^{(j)}_{aa}.
\]
Multiplying \eqref{eq:logder} by $\sin\theta\cos\theta\,w_{j}$ gives $jw_{j}=\sin\theta\cos\theta\,w_{j}'+(k\sin^{2}\theta+g\sin\theta\cos\theta)w_{j}$, so the bracket equals $(k+g\cot\theta)w_{j}=\frac{k}{1-r^{k}}w_{j}$ ($g\cot\theta=kr^{k}/(1-r^{k})$): the first sum is $\mu_{t}\zeta^{t}(s)=\mu_{t}y$. The second is centered Gaussian, independent of the local data, of variance $\sin^{-4}\theta\sum_{j}2j(j-1)w_{j}^{2}=2\beta_{tt}$ (Lemma~\ref{lem:scalars}(ii)), and independent across $a$. Given $(y,0)$, each $B_{aa}$ is therefore $\mu_{t}y$ plus a centered Gaussian variable built on $(\widehat G^{(j)}_{aa})_{2\le j\le k}$, while $B_{rr}$ below is built on the $Z_{j}(v)$ alone; the two families are independent, which gives the joint independence of $B_{11},\dots,B_{nn},B_{rr}$ given $(y,0)$. Next $B_{rr}=-\sum_{j}w_{j}''Z_{j}(v)$ is a function of the $Z_{j}(v)$, hence independent of the $\widehat G^{(j)}$, and unaffected by the conditioning on $\partial_{a}\zeta^{t}(s)=0$. Given $(\zeta^{t}(s),\partial_{\theta}\zeta^{t}(s))=(y,0)$ it is Gaussian with variance at most $\Var B_{rr}=\|w''\|^{2}$ and mean $\gamma^{\top}\Sigma^{-1}(y,0)^{\top}$, $\gamma=\Cov(B_{rr},(\zeta^{t}(s),\partial_{\theta}\zeta^{t}(s)))$, bounded by $\sqrt{\gamma^{\top}\Sigma^{-1}\gamma}\sqrt{(\Sigma^{-1})_{11}}\,y\le\|w''\|\sqrt{g_{r}/\Delta}\,y$, since $\gamma^{\top}\Sigma^{-1}\gamma\le\Var B_{rr}$ is the explained variance.
\end{proof}

\begin{lemma}[Bounds on the weights]\label{lem:wbounds}
For $\theta\in(0,\pi)$: (i) $\|w''(\theta)\|\le W_{k}$, where $W_{k}:=9k^{3}$ for $k$ even and $W_{k}:=43k^{3}$ for $k$ odd; (ii) if $k$ is odd and $\cos\theta\le0$, $V(\theta)\le2^{k}\sin^{4}\theta$.
\end{lemma}

\begin{proof}
With $S_{k}(r)=\sum_{i<k}r^{i}=(1-r^{k})/(1-r)$ and $\sin^{2}\theta=(1-r)(1+r)$, for $j\ge2$
\begin{equation}\label{eq:wregular}
w_{j}(\theta)=s_{k}(\theta)\,T_{j}(\theta),\qquad s_{k}(\theta):=\frac{1}{S_{k}(\cos\theta)},\qquad T_{j}(\theta):=\binom kj^{1/2}\cos^{k-j}\theta\,(1+\cos\theta)\,\sin^{j-2}\theta ,
\end{equation}
a form with no singularity at $\theta=0$. Write $T=(T_{j})_{2\le j\le k}$, so that $w=s_{k}T$ and
\[
\|w''\|\le s_{k}\|T''\|+2|s_{k}'|\,\|T'\|+|s_{k}''|\,\|T\| .
\]
\emph{The vector $T$.} With $p=\sin^{2}\theta$, $\sum_{j\ge2}\binom kj\cos^{2(k-j)}\theta\sin^{2j}\theta=\PP\{\mathrm{Bin}(k,p)\ge2\}\le\binom k2p^{2}$, so $\|T(\theta)\|^{2}=(1+\cos\theta)^{2}\sin^{-4}\theta\sum_{j\ge2}\binom kj\cos^{2(k-j)}\theta\sin^{2j}\theta\le4\binom k2$ for every $\theta\in\R$, hence $\sup_{\R}\|T\|\le\sqrt{2k(k-1)}$. Each $T_{j}$ is a trigonometric polynomial of degree $k-1$, and Bernstein's inequality \citep[Chapter~4]{devore1993} holds for $\R^{k-1}$-valued trigonometric polynomials with the same constant (for fixed $\theta_{0}$, apply it to $\langle T(\cdot),u\rangle$ with~$u=T'(\theta_{0})/\|T'(\theta_{0})\|$); hence $\|T'\|\le(k-1)\sqrt{2k(k-1)}$ and $\|T''\|\le(k-1)^{2}\sqrt{2k(k-1)}$.

\emph{The scalar $s_{k}$ for $\cos\theta\ge0$.} On $[0,1]$, $S_{k}\ge1$, $S_{k}'=\sum_{i<k}i\,r^{i-1}\le(k-1)S_{k}$ and $S_{k}''\le(k-1)(k-2)S_{k}$. From $s_{k}'=S_{k}'\sin\theta/S_{k}^{2}$ and $s_{k}''=(S_{k}'\cos\theta-S_{k}''\sin^{2}\theta)/S_{k}^{2}+2S_{k}'^{2}\sin^{2}\theta/S_{k}^{3}$ we get $s_{k}\le1$, $|s_{k}'|\le k-1$, $|s_{k}''|\le(k-1)+(k-1)(k-2)+2(k-1)^{2}=3(k-1)^{2}$, and
\[
\|w''\|\le\sqrt{2k(k-1)}\,\big[(k-1)^{2}+2(k-1)^{2}+3(k-1)^{2}\big]=6(k-1)^{2}\sqrt{2k(k-1)}\le 9k^{3}.
\]
For $k$ even, $\|w''(\pi-\theta)\|=\|w''(\theta)\|$ by Lemma~\ref{lem:poles}(i), which proves (i) in that case.

\emph{The scalar $s_{k}$ for $\cos\theta<0$, $k$ odd.} Put $r=-x$, $x\in(0,1]$. Then $S_{k}(r)=(1+x^{k})/(1+x)\in[\frac12,1]$, $S_{k}'(r)=n(x)/(1+x)^{2}$ with $n(x)=1-kx^{k-1}-(k-1)x^{k}$ decreasing from $1$ to $2-2k$, and $S_{k}''(r)=k(k-1)x^{k-2}/(1+x)+2S_{k}'(r)/(1+x)$. If $n(x)<0$ then $(2k-1)x^{k-1}>1$, so $x\ge5^{-1/2}$ and $(1+x)^{2}>2$; hence $|S_{k}'|\le k-1$ and $|S_{k}''|\le k(k-1)+2(k-1)\le k(k+1)$. Therefore $s_{k}\le2$, $|s_{k}'|\le4(k-1)$, $|s_{k}''|\le4\big(k(k+1)+k-1\big)+16(k-1)^{2}\le20k^{2}$, and
\[
\|w''\|\le\sqrt{2k(k-1)}\,\big[2(k-1)^{2}+8(k-1)^{2}+20k^{2}\big]\le30k^{2}\sqrt{2k(k-1)}\le43k^{3}.
\]
(ii) For $k$ odd and $r\le0$, $1-r^{k}\ge1$, so $w_{j}^{2}\le\binom kj\sin^{2j}\theta\le\binom kj\sin^{4}\theta$.
\end{proof}

\subsection{The Rice bound}\label{sec:rice}

\begin{lemma}[Gaussian tail integrals]\label{lem:tail}
For $m\ge0$ and $a>0$ with $a(1+a)\ge2m$, $\int_{a}^{\infty}(1+z)^{m}e^{-z^{2}/2}\,\mathrm{d}z\le2\,(1+a)^{m}a^{-1}e^{-a^{2}/2}$, and $a\mapsto(1+a)^{m}e^{-a^{2}/2}$ is decreasing on $\{a(1+a)\ge m\}$. For $m=0$ the factor $2$ can be dropped: $\bar\Phi(a)\le\varphi(a)/a$.
\end{lemma}

\begin{proof}
Let $J$ be the integral. Since $1\le z/a$ on $[a,\infty)$, integrating by parts,
\begin{align*}
J\le\frac1a\int_{a}^{\infty}z(1+z)^{m}e^{-z^{2}/2}\,\mathrm{d}z&=\frac{(1+a)^{m}e^{-a^{2}/2}}{a}+\frac{m}{a}\int_{a}^{\infty}(1+z)^{m-1}e^{-z^{2}/2}\,\mathrm{d}z\\
&\le\frac{(1+a)^{m}e^{-a^{2}/2}}{a}+\frac{m}{a(1+a)}J,
\end{align*}
whence the claim; the logarithmic derivative of $(1+a)^{m}e^{-a^{2}/2}$ is $\frac{m}{1+a}-a$.
\end{proof}

\begin{proposition}[Rice bound with explicit constant]\label{prop:rice}
For $\rho x\ge\sqrt{2(d-1)}$,
\begin{equation}\label{eq:ricebound}
\PP\Big\{\sup_{M}\zeta^{t}>x\Big\}\;\le\;\E\,\Nlm_{(x,\infty)}(\zeta^{t})\;\le\;\bar P(x):=\frac{\sqrt{\pi}\,W_{k}\,2^{2k}\,(k-1)^{(d-2)/2}}{\Gamma\big(\frac{d-1}{2}\big)}\;\frac{(1+\rho x)^{d-1}}{\rho x}\;e^{-\rho^{2}x^{2}/2},
\end{equation}
with $W_{k}$ as in Lemma~\ref{lem:wbounds}; the constant is $C^{\mathrm{Rice}}_{k,d}$ of \eqref{eq:main_bound} for $k$ odd and $2C^{\mathrm{Rice}}_{k,d}$ for $k$ even.
\end{proposition}

\begin{proof}
The first inequality is \eqref{eq:interior}. On $\{B\succ0\}$ Hadamard's inequality gives $0<\det B\le\prod_{i}B_{ii}$, so $\det B\,\1\{B\succ0\}\le\prod_{i}(B_{ii})_{+}$, and by Lemma~\ref{lem:diag}
\[
\E\big[\det B\,\1\{B\succ0\}\mid(y,0)\big]\le\E[(B_{rr})_{+}\mid(y,0)]\prod_{a=1}^{n}\E[(B_{aa})_{+}\mid(y,0)]=:H_{\theta}(y),
\]
with $\E[(B_{aa})_{+}\mid(y,0)]=\sqrt{2\beta_{tt}}\,\varpi(\rho y)$, $\varpi(u)=u\Phi(u)+\varphi(u)\le u+1$ for $u\ge0$, and $\E[(B_{rr})_{+}\mid(y,0)]\le\|w''\|(1+y^{2}g_{r}/\Delta)^{1/2}$. Put $\varkappa:=\sqrt{\Delta/g_{r}}$; then $\varkappa^{2}=V-C^{2}/g_{r}\le V\le\varsigma^{2}$ (Lemma~\ref{lem:variance}), so $\rho\varkappa\le1$. Substituting $y=\varkappa z$ in the inner integral of \eqref{eq:KR}: $\Delta^{-1/2}\varkappa=g_{r}^{-1/2}$, $y^{2}g_{r}/\Delta=z^{2}$, $\rho y\le z$, $(1+z^{2})^{1/2}\le1+z$, and the lower limit is $x/\varkappa\ge\rho x$; hence
\[
\int_{x}^{\infty}p_{\theta}(y)H_{\theta}(y)\,\mathrm{d}y\le\frac{(2\beta_{tt})^{n/2}}{(2\pi)^{d/2}g_{t}^{n/2}}\cdot\frac{\|w''\|}{\sqrt{g_{r}}}\;J\Big(\frac{x}{\varkappa}\Big),
\qquad J(a):=\int_{a}^{\infty}(1+z)^{d-1}e^{-z^{2}/2}\,\mathrm{d}z .
\]
By Lemma~\ref{lem:scalars}(ii), $\sin^{2}\theta\,\beta_{tt}/g_{t}=(k-1)/\sum_{i=0}^{k-2}r^{2i}\le k-1$, so $\sin^{n}\theta(2\beta_{tt}/g_{t})^{n/2}\le(2(k-1))^{n/2}$ and
\begin{equation}\label{eq:thetaint}
\E\,\Nlm_{(x,\infty)}(\zeta^{t})\le\frac{\omega_{n}(2(k-1))^{n/2}}{(2\pi)^{d/2}}\;W_{k}\int_{0}^{\pi}\frac{J(x/\varkappa(\theta))}{\sqrt{g_{r}(\theta)}}\,\mathrm{d}\theta .
\end{equation}
Away from the antipode (all $\theta$ if $k$ is even, $\cos\theta\ge0$ if $k$ is odd), Lemmas~\ref{lem:scalars}(iii) and~\ref{lem:variance} give $g_{r}\ge\Delta/V\ge2^{5-2k}/k^{4}$, and Lemma~\ref{lem:tail} with $m=d-1$ (its hypothesis $\rho x(1+\rho x)\ge2(d-1)$ holds) gives $J(x/\varkappa)\le J(\rho x)\le2(1+\rho x)^{d-1}(\rho x)^{-1}e^{-\rho^{2}x^{2}/2}$: the integrand of \eqref{eq:thetaint} is at most $k^{2}2^{k-3/2}\,(1+\rho x)^{d-1}(\rho x)^{-1}e^{-\rho^{2}x^{2}/2}$. Near the antipode ($k$ odd, $\cos\theta<0$), Lemmas~\ref{lem:scalars}(iii) and~\ref{lem:wbounds}(ii) give $1/\sqrt{g_{r}}\le\sqrt{V/\Delta}\le2^{3k/2-1}/\sin^{2}\theta$ and $\varkappa\le\sqrt V\le2^{k/2}\sin^{2}\theta$, hence $1/\sqrt{g_{r}}\le2^{2k-1}/\varkappa$; with $a=x/\varkappa\ge\rho x$, Lemma~\ref{lem:tail} gives $\varkappa^{-1}J(x/\varkappa)\le\frac{2}{x}(1+a)^{d-1}e^{-a^{2}/2}\le\frac{2}{x}(1+\rho x)^{d-1}e^{-\rho^{2}x^{2}/2}$, so the integrand is at most $2^{2k}\rho\,(1+\rho x)^{d-1}(\rho x)^{-1}e^{-\rho^{2}x^{2}/2}$. As~$k^{2}2^{k-3/2}\le2^{2k}$ and $\rho<1$, the integral in \eqref{eq:thetaint} is at most $\pi2^{2k}(1+\rho x)^{d-1}(\rho x)^{-1}e^{-\rho^{2}x^{2}/2}$, and $\omega_{n}(2(k-1))^{n/2}(2\pi)^{-d/2}=(k-1)^{(d-2)/2}/(\sqrt\pi\,\Gamma(\frac{d-1}{2}))$ give \eqref{eq:ricebound}.
\end{proof}

\subsection{Two moments of the GOE and the proof of Theorem~\ref{thm:target}}\label{sec:final}

\begin{lemma}[Two facts on the GOE]\label{lem:goe}
Let $G\sim\GOE(N)$, $N\ge2$, $\nu>0$.
\begin{itemize}
\item[(i)] $\displaystyle\E\det(\nu I-G)=\sum_{i=0}^{\lfloor N/2\rfloor}\frac{(-1)^{i}N!}{i!\,(N-2i)!\,4^{i}}\,\nu^{N-2i}$, and $\E|\det(\nu I-G)|\ge\E\det(\nu I-G)\ge\tfrac34\nu^{N}$ for $\nu\ge N$.
\item[(ii)] $\PP\{\|G\|>\sqrt{2N}+t\}\le2e^{-t^{2}/2}$ for $t\ge0$; hence $\|\,\|G\|\,\|_{q}\le\sqrt{2N}+\tfrac56\sqrt q$ for $q\ge8$, and for $p\ge8/N$
\[
\big\|\det(\nu I-G)\big\|_{p}\;\le\;\Big(\nu+\sqrt{2N}+\tfrac56\sqrt{Np}\Big)^{N}.
\]
\end{itemize}
\end{lemma}

\begin{proof}
(i)~$\det(\nu I-G)=\sum_{j}(-1)^{j}e_{j}(G)\nu^{N-j}$, $e_{j}(G)$ the sum of the principal $j\times j$ minors, each a $\GOE(j)$ matrix $G'$. In $\E\det G'=\sum_{\pi}\mathrm{sgn}(\pi)\E\prod_{i}G'_{i\pi(i)}$ the entries are centered and independent up to symmetry, so a term survives only if every entry appears an even number of times, i.e.\ $\pi$ is a product of $j/2$ disjoint transpositions; each such $\pi$ contributes $(-1)^{j/2}2^{-j/2}$ and there are $(j-1)!!$ of them. Hence $\E e_{j}(G)=\binom Nj(-1)^{j/2}(j-1)!!2^{-j/2}$ for $j$ even, $0$ for $j$ odd, and $\binom N{2i}(2i-1)!!2^{-i}=N!/(i!(N-2i)!4^{i})$. The terms $t_{i}=\frac{N!}{i!(N-2i)!}(4\nu^{2})^{-i}\nu^{N}$ satisfy $t_{i+1}/t_{i}\le N^{2}/(4\nu^{2})\le\frac14$ for $\nu\ge N$, so the alternating sum is at least $t_{0}-t_{1}\ge\frac34\nu^{N}$.
(ii) $\lmax(G)=\sup_{|v|=1}v^{\top}Gv$ is the supremum of a centered Gaussian process of unit variance, and $\Var(v^{\top}Gv-u^{\top}Gu)=2(1-\langle u,v\rangle^{2})\le4(1-\langle u,v\rangle)=\Var(\sqrt2\langle\xi,v\rangle-\sqrt2\langle\xi,u\rangle)$ for $\xi\sim N(0,I_{N})$, so the Sudakov--Fernique inequality gives $\E\lmax(G)\le\sqrt2\,\E\|\xi\|\le\sqrt{2N}$ and Borell--TIS gives $\PP\{\lmax(G)>\sqrt{2N}+t\}\le e^{-t^{2}/2}$; the factor $2$ comes from $\|G\|=\max(\lmax(G),\lmax(-G))$ and $-G\sim G$. Then $\E(\|G\|-\sqrt{2N})_{+}^{q}=\int_{0}^{\infty}qt^{q-1}\PP\{\|G\|>\sqrt{2N}+t\}\,\mathrm{d}t\le2q\int_{0}^{\infty}t^{q-1}e^{-t^{2}/2}\,\mathrm{d}t=q\,2^{q/2}\Gamma(q/2)\le2\sqrt{\pi q}\,(q/e)^{q/2}e^{1/(6q)}$ by Stirling's bound $\Gamma(z)\le\sqrt{2\pi}\,z^{z-1/2}e^{-z}e^{1/(12z)}$; the $q$-th root is $(2\sqrt{\pi q})^{1/q}e^{1/(6q^{2})}\sqrt{q/e}\le0.82\sqrt q\le\frac56\sqrt q$ for $q\ge8$, the prefactor being decreasing in $q$; Minkowski gives $\|\,\|G\|\,\|_{q}\le\sqrt{2N}+\frac56\sqrt q$. Finally $|\det(\nu I-G)|\le(\nu+\|G\|)^{N}$, so $\|\det(\nu I-G)\|_{p}\le\|\nu+\|G\|\,\|_{Np}^{N}\le(\nu+\|\,\|G\|\,\|_{Np})^{N}$ by Minkowski in~$L^{Np}$, with $q=Np\ge8$.
\end{proof}

\begin{proof}[Proof of Theorem~\ref{thm:target}, $k$ odd]
The upper bound is Lemma~\ref{lem:beaten}, and $\dex(u)-\PP\{\sup X>u\}\le\varepsilon(u)\dex(u)$ is Proposition~\ref{prop:mult}. Let $\rho u\ge\min\{d-1,\tfrac32\sqrt{2(d-1)}\}$ and $x\ge u$; put $N:=d-1\ge2$ and $a:=\rho x$. Then $a\ge2$, $\sqrt{2N}\le a$, $Na^{2}\ge8$ and $a(1+a)\ge N$. Take $p=a^{2}$ in Corollary~\ref{cor:cs}. By Lemma~\ref{lem:goe}(ii), $\|\det(\rho xI-G)\|_{p}\le(a+\sqrt{2N}+\frac56\sqrt N\,a)^{N}\le a^{N}(2+\frac56\sqrt N)^{N}$, since $\sqrt{2N}\le a$. By the lower bound in \eqref{eq:Xi_cs} and $\bar\Phi(a)\ge\frac{a}{1+a^{2}}\varphi(a)$,
\[
P(x)^{-1/p}\le\Big(\frac{\sqrt{2\pi}\,(1+a^{2})}{a}\Big)^{1/a^{2}}e^{1/2}\le3\qquad(a\ge2),
\]
the left side being decreasing in $a$ and equal to $2.61$ at $a=2$; hence $P(x)^{1-1/p}\le3P(x)\le3\bar P(x)$ by Proposition~\ref{prop:rice}. If $\rho u\ge d-1$, then $a\ge N$ and $A(\rho x)\ge\frac34(\rho x)^{N}$ (Lemma~\ref{lem:goe}(i)), so
\[
\frac{\Xi(x)}{A(\rho x)}\;\le\;3\cdot\tfrac43\Big(2+\tfrac56\sqrt{d-1}\Big)^{d-1}\bar P(x)
=C^{\mathrm{gap}}_{k,d}\,\frac{(1+\rho x)^{d-1}}{\rho x}\,e^{-\rho^{2}x^{2}/2}.
\]
If $\tfrac32\sqrt{2N}\le\rho u<N$, write $a=e_{x}\sqrt{2N}$ with $e_{x}\ge\frac32$; then \eqref{eq:denominator_bound} of Appendix~\ref{app:deferred} gives $A(\rho x)\ge(a^{2}-2N-1)^{N/2}=a^{N}(1-e_{x}^{-2})^{N/2}\big(1-\tfrac{1}{2N(e_{x}^{2}-1)}\big)^{N/2}\ge\tfrac45\big(\tfrac59\big)^{N/2}a^{N}$ by Bernoulli's inequality, and the factor $3\cdot\frac43$ becomes $3\cdot\frac54(\frac95)^{N/2}=\frac{15}{4}(\frac95)^{N/2}$. In both cases the bound is decreasing in $x$ on $[u,\infty)$ (the logarithmic derivative of $(1+a)^{d-1}a^{-1}e^{-a^{2}/2}$ is $\frac{d-1}{1+a}-\frac1a-a<0$ for $a(1+a)\ge d-1$); the supremum over $x\ge u$ is at $x=u$, which gives the two bounds of the theorem.
\end{proof}

\begin{proposition}[$k$ even]\label{prop:keven}
Let $k$ be even. Then $\PP\{\sup X>u\}\le\dst(u)=\frac12\dex(u)$ and $\dst(u)-\PP\{\sup X>u\}=\frac12C_{k,d}\int_{u}^{\infty}\Xi(x)\varphi(x)\,\mathrm{d}x$ for $u>0$. Moreover $\1\{\sup_{M}\zeta^{t}>x\}\le\frac12\Nlm_{(x,\infty)}(\zeta^{t})$ almost surely, so $\PP\{\sup_{M}\zeta^{t}>x\}\le\frac12\bar P(x)$ for $\rho x\ge\sqrt{2(d-1)}$, and Theorem~\ref{thm:target} holds as stated.
\end{proposition}

\begin{proof}
$X(-s)=X(s)$ and the antipodal map is an isometry, so critical points come in antipodal pairs with equal values and Hessians and $\Ncp_{[u,\infty)}$ is even. On $\mathbb{RP}^{d-1}$ the field $\bar X([s])=X(s)$ has $\Var(\bar X([s])-\bar X([t]))=2(1-\langle s,t\rangle^{k})>0$ for $[s]\ne[t]$, so its maximum is almost surely attained at a single point \citep[Theorem~3]{lifshits1983absolute}, and $\PP\{\sup X=u\}\le\E\Ncp_{\{u\}}=0$ by \eqref{eq:dexact}. Hence $\Ncp_{[u,\infty)}=2\1\{\sup X\ge u\}+\#\{\text{beaten critical points above }u\}$, and the marked Rice formula of Proposition~\ref{prop:mult}, whose proof uses the parity of $k$ only through \eqref{eq:beating} ($X(-t)=x$ never beats $t$ strictly), gives the identity. If $s_{*}\in M$ is the interior local maximum at which $\sup_{M}\zeta^{t}$ is attained (Proposition~\ref{prop:interior}), then $-s_{*}\in M$, $\zeta^{t}(-s)=\zeta^{t}(s)$ for all $s\in M$ (Lemma~\ref{lem:poles}(iii)), so $-s_{*}\ne s_{*}$ is a local maximum with the same value: on $\{\sup_{M}\zeta^{t}>x\}$, $\Nlm_{(x,\infty)}(\zeta^{t})\ge2$. Sections~\ref{sec:layers} to~\ref{sec:rice} hold for $k$ even without change (Proposition~\ref{prop:interior} treats both fibers), so $\PP\{\sup_{M}\zeta^{t}>x\}\le\frac12\E\Nlm_{(x,\infty)}(\zeta^{t})\le\frac12\bar P(x)$ by Proposition~\ref{prop:rice}, and the proof of Theorem~\ref{thm:target} applies with $\bar P$ replaced by $\frac12\bar P$, which is the factor $\frac12$ in $C^{\mathrm{Rice}}_{k,d}$ for $k$ even.
\end{proof}

\paragraph{Numerical illustration.} Table~\ref{tab:final} lists, in the range $\rho u\ge d-1$ of the theorem, the Rice integral \eqref{eq:KR} with Hadamard's bound $H_{\theta}$ evaluated numerically, its explicit bound $\bar P(u)$, and the two resulting values of the bound on $\Xi(u)/A(\rho u)$ from the proof of Theorem~\ref{thm:target}: with the numerical integral, and with the explicit bound \eqref{eq:main_bound} (scripts in \verifrepo). The explicit $C^{\mathrm{Rice}}_{k,d}$ is crude by a factor between $3\cdot10^{6}$ and $9\cdot10^{6}$, which comes from the uniform bounds $W_{k}$ and the lower bounds on $\Delta$ ($W_{k}$ alone is loose by a factor of order $k^{2}$: numerically $\sup_{\theta}\|w''(\theta)\|\le2k$ for $3\le k\le10$); the factor $3(2+\frac56\sqrt{d-1})^{d-1}$ of the proof is $5.4\cdot10^{2}$ for $d=5$ and $2.6\cdot10^{3}$ for $d=6$. With the numerical integral the bracket is non-vacuous from $1.51$ to $1.56\,u_{\mathrm{IMF}}$ on, $u_{\mathrm{IMF}}=\sqrt{2d-1}/\rho$, with the explicit constant from $2.33$ to $2.40\,u_{\mathrm{IMF}}$ on.

\begin{table}[H]
\centering\footnotesize\setlength\tabcolsep{6pt}
\begin{tabular}{lcccccc}
\toprule
$(k,d)$ & $u$ & $u/u_{\mathrm{IMF}}$ & Rice integral \eqref{eq:KR} & $\bar P(u)$ ($\tfrac12\bar P$ for $k=4$) & $\varepsilon$, numerical & $\varepsilon$, \eqref{eq:main_bound} \\
\midrule
$(3,5)$ & $8.66$ & $2.5$ & $4.8\cdot10^{-11}$ & $1.6\cdot10^{-4}$ & $2.8\cdot10^{-8}$ & $0.12$\\
        & $13.86$ & $4.0$ & $1.0\cdot10^{-29}$ & $4.8\cdot10^{-23}$ & $5.6\cdot10^{-27}$ & $3.5\cdot10^{-20}$\\
$(3,6)$ & $9.57$ & $2.5$ & $7.1\cdot10^{-13}$ & $3.9\cdot10^{-6}$ & $2.0\cdot10^{-9}$ & $0.014$\\
        & $15.32$ & $4.0$ & $1.4\cdot10^{-35}$ & $1.1\cdot10^{-28}$ & $3.7\cdot10^{-32}$ & $3.7\cdot10^{-25}$\\
$(4,6)$ & $10.16$ & $2.5$ & $1.4\cdot10^{-12}$ & $8.8\cdot10^{-6}$ & $2.0\cdot10^{-9}$ & $0.030$\\
        & $16.25$ & $4.0$ & $2.8\cdot10^{-35}$ & $2.4\cdot10^{-28}$ & $3.7\cdot10^{-32}$ & $8.3\cdot10^{-25}$\\
\bottomrule
\end{tabular}
\caption{Numerical against explicit; ``$\varepsilon$, numerical'' is the bound of the proof of Theorem~\ref{thm:target} on $\Xi(u)/A(\rho u)$ with the Rice integral evaluated numerically. For $k=4$, $\bar P$ and $\varepsilon$ include the factor $\frac12$ of Proposition~\ref{prop:keven}. All entries are rounded up.}
\label{tab:final}
\end{table}

\normalsize
\section{A non-asymptotic annealed complexity}\label{sec:kspin}

The Kostlan--Shub--Smale field $X$ coincides with the Hamiltonian of the (centered, isotropic) spherical $k$-spin model from spin-glass physics \citep{auffinger2013,arous2019landscape}. Its critical-point structure, the count of local maxima above an energy level $E$, is a central object of study. Theorem~\ref{thm:complexity} concerns the annealed complexity, the first moment of the count; the quenched complexity is studied through the second moment by \citet{subag2017complexity}, and the critical points near the maximum by \citet{subag2017extremal}. \citet{auffinger2013} write $\Theta_{0,p}(u)$ for the growth rate of the number of local minima of the order-$p$ model below the level $Nu$, and $\Theta_p$ for the rate of all critical points. By the symmetry $H\stackrel d=-H$, the rate of local maxima above $Nu$ is~$\Theta_{0,p}(-u)$, and~$\Theta_{0,p}=\Theta_p$ on $u\le-E_\infty$ \citep[Eqs.~(2.15) and~(2.16)]{auffinger2013}. We write $\Theta_k:=\Theta_{0,p}\big|_{p=k}$. Here $E_\infty:=2\sqrt{(k-1)/k}$ is the spectral-edge energy \citep[Eq.~(2.13)]{auffinger2013}, with $\rho E_\infty=\sqrt2$, and $E_0>E_\infty$ is the ground-state energy above which no critical point survives \citep[Eq.~(2.22) and Theorem~2.12]{auffinger2013}. In our normalization the level $E$ corresponds to $u_{\mathrm{ABC}}=E/\sqrt d$ (Remark~\ref{rem:abc_normalization}).
Theorem~\ref{thm:imf_tail_exact} and Lemma~\ref{lem:goe} turn the asymptotic statement $\lim_{d\to\infty}d^{-1}\log\E[\Nlm_{[E,\infty)}]=\Theta_k(-u_{\mathrm{ABC}})$ into a non-asymptotic two-sided bracket at every finite $(k,d)$.

\begin{theorem}[Non-asymptotic annealed complexity]
\label{thm:complexity}
Let $k,d\ge 3$. For $E>0$ put
\begin{subequations}\label{eq:complexity}
\begin{equation}\label{eq:complexity_thresholds}
    e\;:=\;\frac{\rho E}{\sqrt{2(d-1)}}\,,\qquad \kappa(e)\;:=\;\frac{2e\,\exp(-(e-1)^{2})}{\sqrt{e^{2}-1}}\quad(e>1)\,.
\end{equation}
\begin{enumerate}
    \item[\textup{(i)}] For every $E\in\R$,
    \begin{equation}\label{eq:complexity_upper}
        \E\big[\Nlm_{[E,\infty)}\big]\;\le\;\E\big[N_{[E,\infty)}^{\mathrm{cp}}\big]\;=\;\delta_{\mathrm{exact}}(E)\,,
    \end{equation}
    where $N_{[E,\infty)}^{\mathrm{cp}}$ counts all critical points with $X(s)\ge E$.
    \item[\textup{(ii)}] For every $E$ with $e\ge 3/2$,
    \begin{equation}\label{eq:complexity_bracket}
        \big(1-4\,\kappa(e)^{d-1}\big)\,\delta_{\mathrm{exact}}(E)
        \;\le\;\E\big[\Nlm_{[E,\infty)}\big]
        \;\le\;\delta_{\mathrm{exact}}(E)\,.
    \end{equation}
\end{enumerate}
\end{subequations}
\end{theorem}
The function $\kappa$ is decreasing on $[3/2,\infty)$ and $\kappa(e)<1$ if and only if $e>e_{*}=1.9224\ldots$; thus the lower bound in \eqref{eq:complexity_bracket} is non-trivial for $e\ge 9/4$ at every $d\ge3$, for $e\ge2$ when $d\ge10$ ($\kappa(2)=4\exp(-1)/\sqrt3=0.850$), and for every fixed $e>e_{*}$ when $d$ is large. By Step~3b of the proof in Appendix~\ref{app:deferred}, $E/\sqrt d\to eE_\infty$ as $d\to\infty$: $e$ is the energy in units of the spectral edge.

\begin{corollary}[Annealed complexity of deep local minima, by sign symmetry]
\label{cor:minima_dual}
Let $N_{(-\infty,E']}^{\mathrm{lm,min}}$ denote the number of local minima of $X$ with value at most $E'$.
\begin{subequations}\label{eq:complexity_min}
\begin{enumerate}
    \item[\textup{(i$'$)}] For every $E'\in\R$,
    \begin{equation}\label{eq:complexity_min_upper}
        \E\big[N_{(-\infty,E']}^{\mathrm{lm,min}}\big]\;\le\;\delta_{\mathrm{exact}}(-E')\,.
    \end{equation}
    \item[\textup{(ii$'$)}] For every $E'<0$ with $e':=-\rho E'/\sqrt{2(d-1)}\ge3/2$,
    \begin{equation}\label{eq:complexity_min_bracket}
        \big(1-4\,\kappa(e')^{d-1}\big)\,\delta_{\mathrm{exact}}(-E')
        \;\le\;\E\big[N_{(-\infty,E']}^{\mathrm{lm,min}}\big]
        \;\le\;\delta_{\mathrm{exact}}(-E')\,.
    \end{equation}
\end{enumerate}
\end{subequations}
\end{corollary}
\begin{proof}
The Gaussian symmetry $X\stackrel d=-X$ gives $N_{(-\infty,E']}^{\mathrm{lm,min}}\stackrel d= N_{[-E',\infty)}^{\mathrm{lm}}$ for every $E'\in\R$, and Theorem~\ref{thm:complexity} applies to the right-hand side.
\end{proof}
The bracket (ii$'$) is informative when $4\kappa(e')^{d-1}<1$, that is for minima deeper than about twice the spectral edge ($e'\ge9/4$ at every $d\ge3$, and every fixed $e'>e_{*}$ for $d$ large). The band $1\le e'\le e_{*}$ is covered by (i$'$) only; in the limit $d\to\infty$ it contains the reduced ground-state energy $E_0/E_\infty$ of \citet[Eq.~(2.22) and Theorem~2.12]{auffinger2013} ($1.015$ for $k=3$, $1.13$ for $k=10$, below $e_{*}$ for every $k\le7\cdot10^{4}$).

\begin{remark}[Exact Kac--Rice representation]
\label{rem:kr_representation}
The Kac--Rice formula with the conditional law~\eqref{eq:condhess} of the Hessian yields, for every $E\in\R$,
\begin{equation}\label{eq:complexity_exact}
    \E\big[N_{[E,\infty)}^{\mathrm{lm}}\big]
    \;=\;C_{k,d}\int_E^\infty\E\big[|\det(G-\rho x\,I_{d-1})|\,\1\{G-\rho x\,I_{d-1}\prec 0\}\big]\,\varphi(x)\,\mathrm{d}x\,,
\end{equation}
with $C_{k,d}$ as in~\eqref{eq:dexact}; dropping the indicator gives~\eqref{eq:complexity_upper}.
\end{remark}

\begin{corollary}[Asymptotic match with \citet{auffinger2013} in the high-energy regime]
\label{cor:abc_match}
For $u\le -E_\infty$,
\begin{subequations}\label{eq:abc_match}
\begin{equation}\label{eq:Thetap_definition}
    \Theta_k(u)\;=\;\tfrac12\log(k-1)\,-\,\tfrac{k-2}{4(k-1)}u^2\,-\,I_1(u)\,,\qquad I_1(u):=\tfrac{2}{E_\infty^2}\int_u^{-E_\infty}\sqrt{z^2-E_\infty^2}\,\mathrm{d}z\,.
\end{equation}
As $d\to\infty$ with $E=e\sqrt{2(d-1)}/\rho$ for fixed reduced energy $e\ge 1$,
\begin{equation}\label{eq:complexity_asymp}
    \lim_{d\to\infty}\frac{1}{d}\log\E\big[N_{[E,\infty)}^{\mathrm{lm}}\big]
    \;=\;\Theta_k(-eE_\infty)\,,
\end{equation}
and
\begin{equation}\label{eq:Thetap_explicit}
    \Theta_k(-eE_\infty)\;=\;\tfrac12\log(k-1)\,-\,\tfrac{(k-2)e^2}{k}\,-\,2\!\int_1^e\sqrt{w^2-1}\,\mathrm{d}w\,.
\end{equation}
\end{subequations}
\end{corollary}
\begin{proof}
The function $\Theta_k$ is the complexity function of \citet[Eqs.~(2.14) and~(2.16), Theorem~2.5]{auffinger2013} at $p=k$; \eqref{eq:Thetap_explicit} is \eqref{eq:Thetap_definition} after the substitution $z=-E_\infty w$ in $I_1$. The limit \eqref{eq:complexity_asymp} follows from the equality $\E[N^{\mathrm{cp}}_{[E,\infty)}]=\delta_{\mathrm{exact}}(E)$ of part~\textup{(i)} of Theorem~\ref{thm:complexity}, the level correspondence $E/\sqrt d\to eE_\infty$, and \citet[Theorems~2.5 and~2.8]{auffinger2013}: Step~3 of the proof of Theorem~\ref{thm:complexity} in Appendix~\ref{app:deferred}.
\end{proof}
Exponentially many local maxima persist throughout the band $E_\infty<E/\sqrt d<E_0$, where $\Theta_k(-eE_\infty)>0$.
The restriction $e\ge3/2$ concerns only the non-asymptotic bracket~\eqref{eq:complexity_bracket}, which is informative when $4\kappa(e)^{d-1}<1$. The limit~\eqref{eq:complexity_asymp} holds on the whole range $e\ge1$, since it follows from part~\textup{(i)} and \citet[Theorems~2.5 and~2.8]{auffinger2013}. A finite-$d$ lower bracket on $1<e\le e_{*}$ would need the large deviations of~$\lmax(G)$ near the edge \citep[Theorem~6.2]{ben-guionnet} in place of the Gaussian concentration of Lemma~\ref{lem:goe}(ii); we leave it open.
\begin{remark}[Normalization reconciliation with \citet{auffinger2013}]
\label{rem:abc_normalization}
\citet{auffinger2013} use $H_{N,p}$ on the sphere of radius $\sqrt N$ ($\E[H^2]=N$) with energy $u=H/N$, so that our level $E$ is their $u_{\mathrm{ABC}}=E/\sqrt d$ with ambient dimension $N=d$ (Step~3b of the proof in Appendix~\ref{app:deferred}); above the spectral edge $E_\infty$, critical points are local maxima with overwhelming probability, $E_\infty$ being the right edge of the limiting GOE spectral support; we use $f=H(\sqrt N\,\cdot)/\sqrt N$ on the unit sphere ($\E[f^2]=1$) with energy $E$ on the scale $\sqrt N\,u$. By~\eqref{eq:uABC_to_e}, $e=\sqrt{d/(d-1)}\,u_{\mathrm{ABC}}/E_\infty$, so the two conventions agree as $d\to\infty$.
\end{remark}

\begin{remark}[Three nested non-asymptotic complexity bounds]
\label{rem:complexity_hierarchy}
Substituting $\delta_{\mathrm{IMF}}$ (for $E\ge u_{\mathrm{IMF}}$) or~$\delta_{\mathrm{SMF}}$ (for $E\ge u_{\mathrm{SMF}}$) for $\delta_{\mathrm{exact}}$ in~\eqref{eq:complexity_bracket} yields three nested non-asymptotic upper bounds on the annealed complexity: $\delta_{\mathrm{exact}}$ tightest (Theorem~\ref{thm:imf_tail_exact}), $\delta_{\mathrm{IMF}}$ asymptotically sharp (Theorem~\ref{thm:imf_tail}), $\delta_{\mathrm{SMF}}$ a single monomial past $u^\star_d$ (Corollary~\ref{cor:single_term}); the lower bound of~\eqref{eq:complexity_bracket} is unchanged.
\end{remark}

\appendix

\section{Deferred proofs}\label{app:deferred}

For readability the three longest proofs are collected here.

\begin{proof}[Proof of Theorem~\ref{thm:imf_tail_exact}]
We check each $D_i$ in turn; only $D_2$ requires care.

\medskip

\noindent
$\bullet$ \emph{$D_4$ (squared Hermite).} By Lemma~\ref{lem:mehta} the squared-Hermite part of $Q_d(\nu)$ is $e^{-\nu^2/2}\sum_{j=0}^{d-1}c_j^2H_j(\nu)^2$, so with $\nu=\rho x$ and the change of variable $y=\rho x$,
\[
    \int_u^\infty H_j(\rho x)^2\,e^{-(1+\rho^2)x^2/2}\,\mathrm{d}x\;=\;\frac1\rho\int_{\rho u}^\infty H_j(y)^2\,e^{-\beta y^2/2}\,\mathrm{d}y\;=\;\frac1\rho\,K_j^\beta(\rho u)\,,
    \qquad\beta\;=\;\frac{1+\rho^2}{\rho^2}\;=\;\frac{3k-2}{k}\,,
\]
and summing over $j$,
\begin{equation}\label{eq:Td_exact_intro}
    D_4(u)\;=\;T_d^{\mathrm{exact}}(u)\;:=\;\frac1\rho\sum_{j=0}^{d-1}c_j^2\,K_j^\beta(\rho u)\,,
\end{equation}
each $K_j^\beta$ being given in closed form by Corollary~\ref{cor:Kj_beta}; no inequality is used in this step.

\medskip

\noindent
$\bullet$ \emph{$D_1,D_3$ (linear Hermite).} The change of variable $y=\rho x$ in $\mathcal{L}_d(u)=\int_u^\infty H_{d-1}(\rho x)\,e^{-x^2/2}\,\mathrm{d}x$ gives $\mathcal{L}_d(u)=\rho^{-1}J_{d-1}^{1/\rho^2}(\rho u)$, and Lemma~\ref{lem:Jm_beta} with $\beta=1/\rho^2$, $\gamma=2\rho^2$, $\vartheta=2\rho^2-1=\Lambda$ supplies the closed form. For $d$ odd, $D_3$ exhibits a Mills-ratio-type remainder $\bar\Phi(\rho u\sqrt\beta)=\bar\Phi(u)$ (since $d-1$ is even and $J_0^{1/\rho^2}$ enters); for $d$ even, $D_1$ has no such remainder (since $d-1$ is odd).

\medskip

\noindent
$\bullet$ \emph{$D_2$, cross-integral, polynomial part.} Apply Lemma~\ref{lem:hermite_tail_explicit} to $\mathcal{I}_d^c(\rho x)$:
\[
    \mathcal{I}_d^c(\rho x)\;=\;e^{-(\rho x)^2/2}\sum_{\ell=0}^{\lfloor(d-1)/2\rfloor}\!\frac{2^{\ell+1}(d-1)!!}{(d-2\ell-1)!!}\,H_{d-2\ell-1}(\rho x)\;+\;\1\{d\,\mathrm{even}\}\,2^{d/2}(d-1)!!\!\!\int_{\rho x}^\infty\!e^{-y^2/2}\,\mathrm{d}y\,.
\]
Substituting into~\eqref{eq:Cd_def}, the Hermite-series part becomes
\[
    \sum_\ell\frac{2^{\ell+1}(d-1)!!}{(d-2\ell-1)!!}\!\!\int_u^\infty\!H_{d-1}(\rho x)\,H_{d-2\ell-1}(\rho x)\,e^{-(1+\rho^2)x^2/2}\,\mathrm{d}x\,.
\]
The classical product linearization $H_a(\rho x)\,H_b(\rho x)=\sum_{p=0}^{\min(a,b)}2^p\,p!\,\binom{a}{p}\binom{b}{p}\,H_{a+b-2p}(\rho x)$, combined with~$y=\rho x$ and Lemma~\ref{lem:Jm_beta} with $\beta=(3k-2)/k$, evaluates each integral to $\rho^{-1}J_{a+b-2p}^\beta(\rho u)$, producing the first line of~\eqref{eq:Cd_closed}.

\medskip

\noindent
$\bullet$ \emph{$D_2$, cross-integral, Fubini remainder ($d$ even only).} The remainder contribution to $\mathcal{C}_d(u)$ is, after substituting $\mathcal{I}_d^c$'s Gaussian-tail piece,
\[
    2^{d/2}(d-1)!!\,\mathcal{F}_d(u)\,,\qquad \mathcal{F}_d(u)\;:=\;\int_u^\infty H_{d-1}(\rho x)\,\Big[\!\!\int_{\rho x}^\infty e^{-y^2/2}\,\mathrm{d}y\Big]\,e^{-x^2/2}\,\mathrm{d}x\,.
\]
The double integral is over $\{(x,y):x\ge u,\ y\ge\rho x\}$. By Fubini, swapping the order, $y$ ranges over $[\rho u,\infty)$ and for each such $y$, $x$ ranges over $[u,y/\rho]$:
\[
    \mathcal{F}_d(u)\;=\;\int_{\rho u}^\infty\!e^{-y^2/2}\!\int_u^{y/\rho}\!H_{d-1}(\rho x)\,e^{-x^2/2}\,\mathrm{d}x\,\mathrm{d}y
    \;=\;\frac1\rho\!\int_{\rho u}^\infty\!e^{-y^2/2}\bigl[J_{d-1}^{1/\rho^2}(\rho u)-J_{d-1}^{1/\rho^2}(y)\bigr]\,\mathrm{d}y\,,
\]
where the second equality applies the change of variable $z=\rho x$ to the inner integral and writes the difference of $J^{1/\rho^2}_{d-1}$ at the two endpoints. Splitting at the bracket gives
\begin{equation}\label{eq:Fd_split}
    \mathcal{F}_d(u)\;=\;\underbrace{\frac{J_{d-1}^{1/\rho^2}(\rho u)}{\rho}\!\int_{\rho u}^\infty\!e^{-y^2/2}\,\mathrm{d}y}_{=:\,\mathcal{F}_d^{(1)}(u)}\;-\;\underbrace{\frac{1}{\rho}\!\int_{\rho u}^\infty\!e^{-y^2/2}\,J_{d-1}^{1/\rho^2}(y)\,\mathrm{d}y}_{=:\,\mathcal{F}_d^{(2)}(u)}\,.
\end{equation}

\emph{First piece $\mathcal F_d^{(1)}$.} Direct evaluation: $\int_{\rho u}^\infty e^{-y^2/2}\,\mathrm{d}y = \sqrt{2\pi}\,\bar\Phi(\rho u)$, so
\[
    \mathcal{F}_d^{(1)}(u)\;=\;\sqrt{2\pi}\,\bar\Phi(\rho u)\,\frac{J_{d-1}^{1/\rho^2}(\rho u)}{\rho}\,,
\]
which is the first term of~\eqref{eq:Fubini_F}.

\emph{Second piece $\mathcal F_d^{(2)}$.} For $d$ even, $d-1$ is odd, so the iterated form~\eqref{eq:Jm_beta_iter} of Lemma~\ref{lem:Jm_beta} applied with~$\beta=1/\rho^2$, $\gamma=2/\beta=2\rho^2$, $\vartheta=(2-\beta)/\beta=2\rho^2-1=\Lambda$ has $R^{1/\rho^2}_{d-1}=0$ and reads
\[
    J_{d-1}^{1/\rho^2}(y)\;=\;2\rho^2\,e^{-y^2/(2\rho^2)}\!\sum_{\ell=0}^{(d-2)/2}(2\Lambda)^\ell\,\frac{(d-2)!!}{(d-2\ell-2)!!}\,H_{d-2\ell-2}(y)\,.
\]
Substituting into $\mathcal{F}_d^{(2)}(u)$ and pulling the constants $2\rho^2$ and the outer $1/\rho$ outside the sum,
\begin{align*}
    \mathcal{F}_d^{(2)}(u)
    &\;=\;\frac{1}{\rho}\!\int_{\rho u}^\infty\!e^{-y^2/2}\cdot 2\rho^2\,e^{-y^2/(2\rho^2)}\sum_{\ell=0}^{(d-2)/2}(2\Lambda)^\ell\,\frac{(d-2)!!}{(d-2\ell-2)!!}\,H_{d-2\ell-2}(y)\,\mathrm{d}y\\
    &\;=\;\underbrace{\frac{2\rho^2}{\rho}}_{=\,2\rho}\;\sum_{\ell=0}^{(d-2)/2}(2\Lambda)^\ell\,\frac{(d-2)!!}{(d-2\ell-2)!!}\!\int_{\rho u}^\infty\!H_{d-2\ell-2}(y)\,e^{-(1+1/\rho^2)y^2/2}\,\mathrm{d}y\\
    &\;=\;2\rho\,\sum_{\ell=0}^{(d-2)/2}(2\Lambda)^\ell\,\frac{(d-2)!!}{(d-2\ell-2)!!}\,J_{d-2\ell-2}^{\beta}(\rho u)\,,
\end{align*}
where the combined exponent $-(1/\rho^2+1)y^2/2=-\beta y^2/2$ identifies $\beta=1+1/\rho^2=(1+\rho^2)/\rho^2=(3k-2)/k$, the same $\beta$ appearing in $D_4$ (Corollary~\ref{cor:Kj_beta}), and the last line applies the definition $J_m^\beta(\rho u)=\int_{\rho u}^\infty H_m(y)\,e^{-\beta y^2/2}\,\mathrm{d}y$.

\emph{Assembly.} Substituting $\mathcal{F}_d^{(1)}$ and $\mathcal{F}_d^{(2)}$ into~\eqref{eq:Fd_split} and noting the explicit minus sign in front of $\mathcal{F}_d^{(2)}$, the second-piece coefficient becomes $-2\rho$, yielding the second line of~\eqref{eq:Fubini_F}.

\medskip

\noindent
$\bullet$ \emph{Pointwise tightness.} Let $u\ge u_{\mathrm{IMF}}$. On $\{\rho x\ge\sqrt{2d-1}\}$ both $\mathcal{I}_d^c(\rho x)>0$ (Lemma~\ref{lem:Idc_positivity}) and $H_{d-1}(\rho x)>0$, so $D_2(u)<0$; with $D_1+D_3=\alpha_d\,\mathcal{L}_d(u)\le\alpha_d\,\Phi_d(\rho,u)\,e^{-u^2/2}$ (Proposition~\ref{prop:IMF_decomp}) this gives $\delta_{\mathrm{exact}}(u)<\delta_{\mathrm{IMF}}^\star(u)$. For the second inequality, $\sqrt{2d-1}>x_1^{(j)}$ for every $j\le d-1$, so Lemma~\ref{lem:hermite_bound} gives $0\le H_j(\rho x)\le(2\rho x)^j$ on $\{\rho x\ge\sqrt{2d-1}\}$, with equality for $j=0,1$ and strict inequality for $j\ge2$; hence
\[
    \int_u^\infty H_j(\rho x)^2\,e^{-(1+\rho^2)x^2/2}\,\mathrm{d}x\;\le\;\int_u^\infty(2\rho x)^{2j}\,e^{-(1+\rho^2)x^2/2}\,\mathrm{d}x\,,
\]
strictly for $j\ge2$, an index present since $d\ge3$. Summing with the weights $c_j^2$ and bounding the right-hand side as in Proposition~\ref{prop:IMF_decomp} gives $D_4(u)<\Psi_d(\rho,u)\,e^{-(1+\rho^2)u^2/2}$, that is $\delta_{\mathrm{IMF}}^\star(u)<\delta_{\mathrm{IMF}}(u)$.
\end{proof}

\medskip

\begin{proof}[Proof of Lemma~\ref{lem:dominant_term}]

\noindent
$\bullet$
We separate the Mehta expansion~\eqref{eq:fmehta} into three distinct pieces, $Q_d = T_1 + T_2 + T_3$, defined as follows:
\begin{align*}
    T_1(\nu) &= e^{-\nu^2/2}\sum_{j=0}^{d-1}c_j^2 H_j^2(\nu)\,, \\
    T_2(\nu) &= \tfrac{1}{2}\sqrt{d/2}\,c_{d-1}c_d\,H_{d-1}(\nu)\,\big[\mu_d-2\mathcal{I}_d^c(\nu)\big]\,, \\
    T_3(\nu) &= \1\{d\text{ odd}\}\,\frac{H_{d-1}(\nu)}{\mu_{d-1}}\,.
\end{align*}

\noindent
$\bullet$
We extract the constant prefactor multiplying $H_{d-1}(\nu)$, which depends on the parity of the dimension $d$:
\begin{itemize}
    \item \textbf{For $d$ even:} $T_3 = 0$, and the constant part of the bracket in $T_2$ isolates $\alpha_d = \tfrac{1}{2}\sqrt{d/2}\,c_{d-1}c_d\,\mu_d$.
    \item \textbf{For $d$ odd:} $\mu_d = 0$ (by parity), leaving only $T_3$ to contribute. This gives $\alpha_d = 1/\mu_{d-1}$ (here $d-1$ is even, so $\mu_{d-1} \neq 0$).
\end{itemize}
The exact values of the integrals are obtained via direct evaluation of the Gaussian moments
\[
    \int_{\R} y^{2j}\,e^{-y^2/2}\,\mathrm{d}y \;=\; \sqrt{2\pi}\,(2j-1)!!\,,
\]
yielding the closed form $\mu_{2p}=\sqrt{2\pi}\,(2p)!/p!$.

\noindent
$\bullet$
The total remainder $\mathcal R_d(\nu):=Q_d(\nu)-\alpha_d H_{d-1}(\nu)$ has the same closed form in both parities:
\begin{itemize}
    \item \textbf{For $d$ even}, $T_3=0$ and the constant part of $T_2$'s bracket is $\alpha_d H_{d-1}(\nu)$, so $\mathcal R_d = T_1 + (T_2-\alpha_d H_{d-1}) = T_1 - \sqrt{d/2}\,c_{d-1}c_d\,\mathcal I_d^c(\nu)\,H_{d-1}(\nu)$.
    \item \textbf{For $d$ odd}, $\mu_d=0$ so~$T_2 = -\sqrt{d/2}\,c_{d-1}c_d\,\mathcal I_d^c(\nu)\,H_{d-1}(\nu)$, and $T_3 = H_{d-1}(\nu)/\mu_{d-1} = \alpha_d H_{d-1}(\nu)$, so $\mathcal R_d = T_1 + T_2 = T_1 - \sqrt{d/2}\,c_{d-1}c_d\,\mathcal I_d^c(\nu)\,H_{d-1}(\nu)$.
\end{itemize}
In both cases
\begin{equation}\label{eq:Rd_unified}
    \mathcal R_d(\nu) \;=\; T_1(\nu) \;-\; \sqrt{d/2}\,c_{d-1}c_d\,\mathcal I_d^c(\nu)\,H_{d-1}(\nu)\,,
\end{equation}
and we bound the two terms separately.

We first bound $T_1$. The physicist Hermite polynomial admits the explicit expansion
\[
    H_j(\nu)\;=\;\sum_{m=0}^{\lfloor j/2\rfloor}\frac{(-1)^m\,j!}{m!\,(j-2m)!}(2\nu)^{j-2m}\,,
\]
and we use the elementary coefficient inequality $j!\,2^{j-2m}/m!\le(2j)!/(2m)!$, valid for every $j\ge m\ge 0$. The latter is proved by induction: the base case $j=m$ gives $2^{-m}\le 1$, and the step $j\to j+1$ multiplies the LHS by $2(j+1)$ and the RHS by $2(2j+1)(j+1)$, giving a ratio of $2j+1\ge 1$. The triangle inequality then yields
\[
    |H_j(\nu)| \;\le\; \sum_{m=0}^{\lfloor j/2\rfloor}\frac{(2j)!}{(2m)!(j-2m)!}|\nu|^{j-2m}
    \;\le\;\sum_{p=0}^{j}\frac{(2j)!}{j!}\binom{j}{p}|\nu|^p
    \;=\;\frac{(2j)!}{j!}(1+|\nu|)^j\,,
\]
where the second inequality adds the (non-negative) odd-power monomials $|\nu|^p$ with $p\not\equiv j\pmod 2$, using $(2j)!/((2m)!(j-2m)!)=(2j)!/j!\cdot j!/((2m)!(j-2m)!)\le(2j)!/j!\cdot\binom{j}{j-2m}$. Squaring and applying $(1+|\nu|)^2\le 2(1+\nu^2)$,
\[
    T_1(\nu) \;\le\; S_d(1+\nu^2)^{d-1}e^{-\nu^2/2}\,, \qquad \text{with $S_d$ as in Lemma~\ref{lem:dominant_term}}.
\]

\noindent
$\bullet$ \emph{A uniform envelope for $\mathcal{I}_d^c$.} We bound $\mathcal{I}_d^c$ on the whole half-line $\{\nu\ge 0\}$ by a single envelope carrying the same polynomial--Gaussian form $(1+|\nu|)^{d-1}e^{-\nu^2/2}$ as $H_{d-1}$. The explicit Hermite-tail series of Lemma~\ref{lem:hermite_tail_explicit} (applied with degree $m=d$) reads
\[
    \mathcal{I}_d^c(\nu)\;=\;e^{-\nu^2/2}\!\!\sum_{\ell=0}^{\lfloor(d-1)/2\rfloor}\!2^{\ell+1}\frac{(d-1)!!}{(d-2\ell-1)!!}\,H_{d-2\ell-1}(\nu)\;+\;R_d(\nu)\,,
\]
with $R_d(\nu)=0$ for $d$ odd and $R_d(\nu)=2^{d/2}(d-1)!!\int_\nu^\infty e^{-y^2/2}\,\mathrm{d}y$ for $d$ even. By the coefficient bound $|H_m(\nu)|\le\frac{(2m)!}{m!}(1+|\nu|)^m$ established above for $T_1$ (whence $|H_{d-2\ell-1}(\nu)|\le\frac{(2(d-2\ell-1))!}{(d-2\ell-1)!}(1+|\nu|)^{d-2\ell-1}\le\frac{(2(d-2\ell-1))!}{(d-2\ell-1)!}(1+|\nu|)^{d-1}$ since $1+|\nu|\ge 1$) and the monotone-Mills bound $\int_\nu^\infty e^{-y^2/2}\,\mathrm{d}y\le\sqrt{\pi/2}\,e^{-\nu^2/2}$ (valid for every $\nu\ge 0$, because $\nu\mapsto e^{\nu^2/2}\int_\nu^\infty e^{-y^2/2}\,\mathrm{d}y$ is non-increasing, with value $\sqrt{\pi/2}$ at $\nu=0$, since its derivative $\nu\,e^{\nu^2/2}\!\int_\nu^\infty e^{-y^2/2}\mathrm{d}y-1$ is negative by the Mills inequality $\nu\,e^{\nu^2/2}\!\int_\nu^\infty e^{-y^2/2}\mathrm{d}y<1$), we obtain
\begin{equation}\label{eq:Idc_uniform_env}
\begin{aligned}
    |\mathcal{I}_d^c(\nu)|&\;\le\;E_d\,(1+|\nu|)^{d-1}\,e^{-\nu^2/2}\,,\\
    E_d&\;:=\;\sum_{\ell=0}^{\lfloor(d-1)/2\rfloor}2^{\ell+1}\frac{(d-1)!!}{(d-2\ell-1)!!}\frac{(2(d-2\ell-1))!}{(d-2\ell-1)!}\;+\;2^{d/2}(d-1)!!\sqrt{\tfrac\pi2}\,.
\end{aligned}
\end{equation}

\noindent
$\bullet$ \emph{$E_d$ is dominated by $\widetilde B_d$.} The summands $a_\ell:=2^{\ell+1}\frac{(d-1)!!}{(d-2\ell-1)!!}\frac{(2(d-2\ell-1))!}{(d-2\ell-1)!}$ are strictly decreasing in $\ell$: with $p:=d-2\ell-1\ge 2$, a direct computation gives the ratio $a_{\ell+1}/a_\ell=\frac{p}{2(2p-1)(2p-3)}\le\frac{1}{2(2p-3)}<1$. Hence, using $a_0=2\frac{(2d-2)!}{(d-1)!}$ and that the sum has at most $\lfloor(d-1)/2\rfloor+1\le\frac{d+1}{2}$ terms,
\[
    \sum_{\ell}a_\ell\;\le\;\frac{d+1}{2}\,a_0\;=\;(d+1)\frac{(2d-2)!}{(d-1)!}\,,\qquad
    2^{d/2}(d-1)!!\sqrt{\tfrac\pi2}\;\le\;2^{d}\frac{(2d-2)!}{(d-1)!}\,,
\]
the second bound using $(d-1)!!\le(d-1)!\le\frac{(2d-2)!}{(d-1)!}$ and $\sqrt{\pi/2}\le 2^{d/2}$. Therefore
\[
    E_d\;\le\;(d+1+2^d)\frac{(2d-2)!}{(d-1)!}\;\le\;(2d-1)\,2^{d+2}\,\frac{(2d-2)!}{(d-1)!}\;=\;\frac{(2d)!\,2^{d+1}}{d!}\;=\;\widetilde B_d\,,
\]
so that, by~\eqref{eq:Idc_uniform_env}, $|\mathcal{I}_d^c(\nu)|\le\widetilde B_d\,(1+|\nu|)^{d-1}\,e^{-\nu^2/2}$ for every $\nu\ge 0$.

\noindent
$\bullet$ \emph{Assembly.} Combining the last bound with $|H_{d-1}(\nu)|\le\frac{(2d-2)!}{(d-1)!}(1+|\nu|)^{d-1}$ and the elementary inequality $(1+|\nu|)^{2(d-1)}\le 2^{d-1}(1+\nu^2)^{d-1}$ (which follows from $(1+|\nu|)^2\le 2(1+\nu^2)$), the cross term in~\eqref{eq:Rd_unified} obeys
\[
    \sqrt{d/2}\,c_{d-1}c_d\,|\mathcal{I}_d^c(\nu)|\,|H_{d-1}(\nu)|\;\le\;\sqrt{d/2}\,c_{d-1}c_d\,\frac{(2d-2)!}{(d-1)!}\,2^{d-1}\,\widetilde B_d\,(1+\nu^2)^{d-1}e^{-\nu^2/2}\,.
\]
Adding the bound $T_1(\nu)\le S_d\,(1+\nu^2)^{d-1}e^{-\nu^2/2}$ established above yields the remainder bound~\eqref{eq:remainder_bound} with the constant $\beta_d$ of \eqref{eq:beta_d}.
\end{proof}

\medskip

\begin{proof}[Proof of Theorem~\ref{thm:complexity} and Corollary~\ref{cor:abc_match}]
The inputs are the Kac--Rice formula of Section~\ref{sub:reduction}, the conditional GOE law of the Hessian \citep[Lemma~3.2(b)]{auffinger2013}, and the two facts on the GOE of Lemma~\ref{lem:goe}.

\medskip

\noindent
$\bullet$ \emph{Step~1: Exact Kac--Rice representation and upper bound.}
The Kac--Rice computation is the one of Section~\ref{sub:reduction}, displays~\eqref{eq:condhess}--\eqref{eq:dexact}, with one change: the target is the expected count of local maxima with value in $[E,\infty)$, for which the Kac--Rice formula of Section~\ref{sub:reduction} is an equality, so the indicator $\1\{\nabla^2 X(t)\prec 0\}$ is kept inside the conditional expectation. By~\eqref{eq:condhess}, $\nabla^2 X(t)\prec 0$ if and only if $G-\rho x\,I_{d-1}\prec 0$, that is $\lmax(G)<\rho x$. The constants gather as in~\eqref{eq:dexact} and give the exact representation~\eqref{eq:complexity_exact}. Dropping the indicator $\1\{G-\rho x\,I_{d-1}\prec 0\}\le 1$ gives $\E[N_{[E,\infty)}^{\mathrm{lm}}]\le C_{k,d}\int_E^\infty A(\rho x)\,\varphi(x)\,\mathrm{d}x=\delta_{\mathrm{exact}}(E)$ by~\eqref{eq:dexact}, which proves part~(i).

\smallskip

\noindent
$\bullet$ \emph{Step~2: Lower bound for $e\ge3/2$.}
Write $N:=d-1$, $\lmax:=\lmax(G)$, and fix $x\ge E$; put $\nu:=\rho x$ and $e_{x}:=\nu/\sqrt{2N}\ge e\ge3/2$. Since $G-\nu I\prec0$ if and only if $\lmax<\nu$, and $|\det(G-\nu I)|=\det(\nu I-G)$ on that event,
\begin{equation}\label{eq:lm_split}
    \E\big[|\det(G-\nu I)|\,\1\{G-\nu I\prec0\}\big]\;=\;A(\nu)\;-\;\E\big[|\det(\nu I-G)|\,\1\{\lmax\ge\nu\}\big]\,,
    \qquad A(\nu)=\E|\det(G-\nu I)|\,.
\end{equation}

\emph{(2a) Lower bound on $A(\nu)$.} By Lemma~\ref{lem:goe}(i), $\E\det(\nu I-G)=\sum_{i}\frac{(-1)^{i}N!}{i!\,(N-2i)!\,4^{i}}\nu^{N-2i}=2^{-N}H_{N}(\nu)$, the physicists' Hermite polynomial, whose zeros $x_{1}>x_{2}>\dots>x_{N}$ are symmetric about $0$ and satisfy~\mbox{$x_{1}<\sqrt{2N+1}$} (Szeg\H{o}, as in Lemma~\ref{lem:Idc_positivity}). Since $\nu\ge\frac32\sqrt{2N}>\sqrt{2N+1}$, pairing the zeros $\pm x_{j}$ (for $N$ odd the zero at $0$ contributes a factor $\nu\ge(\nu^{2}-x_{1}^{2})^{1/2}$),
\begin{equation}\label{eq:denominator_bound}
    A(\nu)\;\ge\;\E\det(\nu I-G)\;=\;\prod_{j=1}^{N}(\nu-x_{j})\;\ge\;(\nu^{2}-x_{1}^{2})^{N/2}\;\ge\;(\nu^{2}-2N-1)^{N/2}\,.
\end{equation}

\emph{(2b) Upper bound on the remainder.} Since $|\det(\nu I-G)|\le(\nu+\|G\|)^{N}$ and $\{\lmax\ge\nu\}\subset\{\|G\|\ge\nu\}$, the layer-cake formula for $s\mapsto(\nu+s)^{N}$ gives
\[
    \E\big[|\det(\nu I-G)|\,\1\{\lmax\ge\nu\}\big]\;\le\;(2\nu)^{N}\,\PP\{\|G\|\ge\nu\}\;+\;N\int_{\nu}^{\infty}(\nu+s)^{N-1}\,\PP\{\|G\|\ge s\}\,\mathrm{d}s\,.
\]
By Lemma~\ref{lem:goe}(ii), with $t_{0}:=\nu-\sqrt{2N}=(e_{x}-1)\sqrt{2N}$ and the substitution $s=\sqrt{2N}+t$,
\[
    \E\big[|\det(\nu I-G)|\,\1\{\lmax\ge\nu\}\big]\;\le\;2(2\nu)^{N}\exp(-t_{0}^{2}/2)\;+\;2N\int_{t_{0}}^{\infty}(\nu+\sqrt{2N}+t)^{N-1}\exp(-t^{2}/2)\,\mathrm{d}t\,.
\]
For $t\ge t_{0}$, $-\frac{\mathrm{d}}{\mathrm{d}t}\big[(\nu+\sqrt{2N}+t)^{N-1}\exp(-t^{2}/2)\big]=(\nu+\sqrt{2N}+t)^{N-1}\exp(-t^{2}/2)\big(t-\tfrac{N-1}{\nu+\sqrt{2N}+t}\big)$ and $t-\frac{N-1}{\nu+\sqrt{2N}+t}\ge t_{0}-\frac{N}{2\nu}=\sqrt{2N}\big((e_{x}-1)-\tfrac{1}{4e_{x}}\big)=:c>0$, so, as $\nu+\sqrt{2N}+t_{0}=2\nu$ and $2\nu c=N\big(4e_{x}(e_{x}-1)-1\big)\ge2N$ for $e_{x}\ge3/2$,
\begin{align*}
    N\int_{t_{0}}^{\infty}(\nu+\sqrt{2N}+t)^{N-1}\exp(-t^{2}/2)\,\mathrm{d}t&\;\le\;\frac{N(2\nu)^{N-1}\exp(-t_{0}^{2}/2)}{c}\\
    &\;=\;\frac{(2\nu)^{N}\exp(-t_{0}^{2}/2)}{4e_{x}(e_{x}-1)-1}\;\le\;\tfrac12(2\nu)^{N}\exp(-t_{0}^{2}/2)\,.
\end{align*}
Hence
\begin{equation}\label{eq:numerator_bound}
    \E\big[|\det(\nu I-G)|\,\1\{\lmax\ge\nu\}\big]\;\le\;3\,(2\nu)^{N}\,\exp\big(-(e_{x}-1)^{2}N\big)\,.
\end{equation}

\emph{(2c) The amplifier.} Dividing \eqref{eq:numerator_bound} by \eqref{eq:denominator_bound} and writing $\nu=e_{x}\sqrt{2N}$,
\[
    \frac{\E\big[|\det(\nu I-G)|\,\1\{\lmax\ge\nu\}\big]}{A(\nu)}\;\le\;3\bigg[\frac{2e_{x}\exp(-(e_{x}-1)^{2})}{\sqrt{e_{x}^{2}-1-\frac{1}{2N}}}\bigg]^{N}\;\le\;3\bigg[\frac{2e\,\exp(-(e-1)^{2})}{\sqrt{e^{2}-1-\frac{1}{2N}}}\bigg]^{N}\,,
\]
the bracket being decreasing in $e_{x}$ on $[3/2,\infty)$ (its logarithmic derivative is $\frac1{e_{x}}-2(e_{x}-1)-\frac{e_{x}}{e_{x}^{2}-1-1/(2N)}$, and $\frac1{e_{x}}-2(e_{x}-1)\le-\frac13$ there). Finally $e^{2}-1\ge\frac54$ and Bernoulli's inequality give $\big(1-\frac{1}{2N(e^{2}-1)}\big)^{-N/2}\le\big(1-\frac{1}{4(e^{2}-1)}\big)^{-1}\le\frac54$, whence
\begin{equation}\label{eq:amplifier_bound}
    \frac{\E\big[|\det(\nu I-G)|\,\1\{\lmax\ge\nu\}\big]}{A(\nu)}\;\le\;\tfrac{15}{4}\,\kappa(e)^{N}\;\le\;4\,\kappa(e)^{d-1}\qquad\text{for every }x\ge E\,.
\end{equation}

\emph{(2d) Integration.} By \eqref{eq:lm_split} and \eqref{eq:amplifier_bound}, $\E[|\det(G-\rho xI)|\1\{G-\rho xI\prec0\}]\ge(1-4\kappa(e)^{d-1})A(\rho x)$ for all $x\ge E$; multiplying by $C_{k,d}\varphi(x)$, integrating over $[E,\infty)$ and using \eqref{eq:complexity_exact} together with $C_{k,d}\int_{E}^{\infty}A(\rho x)\varphi(x)\,\mathrm{d}x=\delta_{\mathrm{exact}}(E)$ from Step~1 gives the lower bound in \eqref{eq:complexity_bracket}.

\smallskip

\noindent
$\bullet$ \emph{Step~3: High-dimensional limit via the ABC complexity theorem.}
We avoid redoing the asymptotic of $\delta_{\mathrm{exact}}$ from scratch via the Mehta--Fyodorov decomposition, which would require fine control of the partial Hermite sum $\Phi_d(\rho,E)$ at the natural scale $\rho E=e\sqrt{2(d-1)}\asymp\sqrt d$, where the leading-monomial approximation $\Phi_d\sim 2\rho(2\rho E)^{d-2}$ is invalid (every term in the sum contributes at comparable order). Instead, we reduce~\eqref{eq:complexity_asymp} directly to \citet[Theorem~2.8]{auffinger2013} via an explicit field/level correspondence between our paper's KSS field $X$ on $\Sph^{d-1}$ and ABC's spherical $k$-spin Hamiltonian $H_{N,p}$ on $S^{N-1}(\sqrt N)$.

\smallskip

\emph{(3a) Field and Hessian correspondence.} \citet{auffinger2013} work, throughout their Section~3, with the rescaled field
\begin{subequations}\label{eq:abc_setting}
\begin{equation}\label{eq:abc_rescaled_field}
    f_{N,p}(s)\;:=\;\frac{1}{\sqrt N}\,H_{N,p}(\sqrt N\,s)\,,\qquad s\in\Sph^{N-1}\,,
\end{equation}
of unit variance on the unit sphere, with covariance $\E[f(s)f(s')]=\langle s,s'\rangle^p$ \citep[eqs.~(3.1)--(3.2)]{auffinger2013}. Identifying their $N$ with $d$ and $p$ with $k$, the law of $f$ on $\Sph^{d-1}$ coincides with that of our paper's $X$. ABC's conditional-Hessian formula \citep[Lemma~3.2, part~(b)]{auffinger2013} is stated directly in this unit-sphere setting: for every $s\in\Sph^{N-1}$ and every $x\in\R$,
\begin{equation}\label{eq:abc_conditional_hessian}
    \nabla^2 f(s)\,\big|\,f(s)=x\;\stackrel d=\;M^{N-1}\,\sqrt{2(N-1)\,p(p-1)}\;-\;x\,p\,I_{N-1}\,,
\end{equation}
\end{subequations}
where $M^{N-1}$ is a GOE in the ABC normalization $\E[(M^{N-1}_{ij})^2]=(1+\delta_{ij})/(2(N-1))$. Setting $G_{N-1}:=\sqrt{N-1}\,M^{N-1}$, which has the normalization $\E[(G_{N-1})_{ij}^2]=(1+\delta_{ij})/2$, and using $p/\sqrt{2p(p-1)}=\sqrt{p/(2(p-1))}=:\rho$, the right-hand side rearranges to $\sqrt{2p(p-1)}\,(G_{N-1}-\rho\,x\,I_{N-1})$. With $(N,p)=(d,k)$ this is exactly the conditional Hessian law~\eqref{eq:condhess}
\begin{equation*}
    \nabla^2 X(s)\,\big|\,X(s)=x\;\stackrel d=\;\sqrt{2k(k-1)}\,(G-\rho\,x\,I_{d-1})\,,\qquad\rho=\sqrt{k/(2(k-1))}\,,
\end{equation*}
the law used throughout the paper. No metric rescaling of the Hessian under a sphere-radius change is needed: ABC perform the rescaling once in passing from $H_{N,p}$ on $S^{N-1}(\sqrt N)$ to $f_{N,p}$ on $\Sph^{N-1}$ via~\eqref{eq:abc_rescaled_field}, and their Hessian formula~\eqref{eq:abc_conditional_hessian} is already on the unit sphere.

\smallskip

\emph{(3b) Level correspondence.} The threshold $X(s)\ge E$ in our paper corresponds to $H(\sqrt N\,s)\ge E\sqrt N$, i.e., $H/N\ge E/\sqrt N=:u_{\mathrm{ABC}}$. On the scale $E=e\sqrt{2(d-1)}/\rho$, using $\rho E_\infty=\sqrt 2$ (which follows from $\rho^2=k/(2(k-1))$ and $E_\infty^2=4(k-1)/k$),
\begin{equation}\label{eq:uABC_to_e}
    u_{\mathrm{ABC}}\;=\;\frac{E}{\sqrt d}\;=\;\frac{e\sqrt{2(d-1)}}{\rho\sqrt d}\;\xrightarrow[d\to\infty]{}\;\frac{e\sqrt 2}{\rho}\;=\;\frac{e\sqrt 2}{\rho E_\infty}\cdot E_\infty\;=\;e\,E_\infty\,.
\end{equation}
For $e\ge 1$ the limiting level satisfies $-eE_\infty\le-E_\infty$, so the case-1 branch of~\eqref{eq:Thetap_definition} applies at $eE_\infty$; at finite~$d$,~$u_{\mathrm{ABC}}=eE_\infty\sqrt{(d-1)/d}<eE_\infty$, and Step~(3c) passes to the limit by monotonicity and continuity.

\smallskip

\emph{(3c) Asymptotic rate for the total critical-point count.} By \citet[Theorem~2.8]{auffinger2013} (log-rate for the total critical-point count below level $Nu$) combined with the upper-tail symmetry $H\stackrel d=-H$ (an immediate consequence of $J\stackrel d=-J$ in the disorder), for every fixed $u\ge E_\infty$ (for smaller $u$ the rate is ABC's $\Theta_p(-u)$, whose formula differs from~\eqref{eq:Thetap_definition} but which is continuous on $\R$):
\[
    \lim_{N\to\infty}\frac1N\log\E\big[\#\{\sigma:\nabla H(\sigma)=0,\,H(\sigma)\in(Nu,\infty)\}\big]\;=\;\Theta_k(-u)\,.
\]
Since $u_{\mathrm{ABC}}\to eE_\infty$, the count is non-increasing in the level and $\Theta_p$ is continuous, the limit persists along the $d$-dependent level $u_{\mathrm{ABC}}$: bracket $u_{\mathrm{ABC}}$ between $eE_\infty\mp\eta$, apply the display at both fixed levels, and let $\eta\to0$.
By the Kac--Rice equality~\eqref{eq:dexact} (recall from Step~1 that $\delta_{\mathrm{exact}}(E)$ equals the Kac--Rice expected critical-point count above level $E$), the left-hand side equals $\lim_d \tfrac1d\log\delta_{\mathrm{exact}}(E)$ under the field/level correspondence of (3a)--(3b). Hence
\begin{equation}\label{eq:delta_exact_asymptotic}
    \lim_{d\to\infty}\frac1d\log\delta_{\mathrm{exact}}(E)\;=\;\Theta_k(-eE_\infty)
    \;=\;\tfrac12\log(k-1)-\tfrac{(k-2)e^2}{k}-2\!\int_1^e\sqrt{w^2-1}\,\mathrm{d}w
\end{equation}
for $e\ge 1$, using~\eqref{eq:Thetap_explicit}.

\smallskip

\emph{(3d) Local-maxima specialization and the bracketing.} ABC's Theorem~2.8 covers the total critical-point count; for local maxima (Morse index $N-1$ on the sphere) specifically, Remark~2.6 following \citet[Theorem~2.5]{auffinger2013} gives
$\lim_N\frac1N\log\E[\#\{\text{local maxima of }H\text{ above }Nu\}]=\Theta_{0,p}(-u)$,
and the case-1 branch $\Theta_{0,p}(u)=\Theta_p(u)$ on $u\le-E_\infty$~\citep[Eqs.~(2.15) and~(2.16)]{auffinger2013} states that for $e\ge 1$, the local-maxima rate coincides with the total-critical-point rate: every critical point above the spectral edge is, with overwhelming probability, a local maximum. Therefore, by the bracketing of (3c), $\Theta_{0,p}$ being continuous on $\R$,
\[
    \lim_{d\to\infty}\frac1d\log\E[N_{[E,\infty)}^{\mathrm{lm}}]\;=\;\Theta_k(-eE_\infty).
\]
Equivalently, and consistently with our non-asymptotic bracketing of~\eqref{eq:complexity_bracket}, combining the upper bound $\E[N^{\mathrm{lm}}]\le\delta_{\mathrm{exact}}$ with the lower bound $\E[N^{\mathrm{lm}}]\ge(1-4\kappa(e)^{d-1})\delta_{\mathrm{exact}}$ (valid for $e\ge3/2$), and observing that the amplifier $4\kappa(e)^{d-1}$ vanishes exponentially as $d\to\infty$ for $e>e_{*}$,
\[
    \tfrac1d\log\delta_{\mathrm{exact}}+\tfrac1d\log\big(1-4\kappa(e)^{d-1}\big)\;\le\;\tfrac1d\log\E[N^{\mathrm{lm}}]\;\le\;\tfrac1d\log\delta_{\mathrm{exact}}\,,
\]
both bracketing terms converge to $\Theta_k(-eE_\infty)$ by~\eqref{eq:delta_exact_asymptotic}. The limit~\eqref{eq:complexity_asymp} is thus established for every~$e\ge 1$ by the exact Kac--Rice equality of Step~1 and \citet[Theorems~2.5 and~2.8]{auffinger2013} (Steps~3c--3d); on~$e>e_{*}$ the non-asymptotic bracket above re-derives the same limit with explicit finite-$d$ error bars. The band $1\le e\le e_{*}$ is covered by the asymptotic statement~\eqref{eq:complexity_asymp}, while only the matching finite-$d$ lower bracket, not the limit itself, lies outside the scope of the present argument.\qedhere
\end{proof}

\section{Technical lemmas}

\subsection{Gaussian-type tail integrals}
\label{app:gaussian_integrals}
The proofs of the three tail bounds use a common tool: integration by parts on tail integrals of the form $\int_u^\infty x^m e^{-ax^2/2}\,\mathrm{d}x$, iterated in steps of two. We isolate it here.

\begin{lemma}[Gaussian-type tail integrals]
\label{lem:gauss_tail}
For $a>0$, integer $m\ge 0$, and $u>0$, define $I_m(u;a):=\int_u^\infty x^m e^{-ax^2/2}\,\mathrm{d}x$. The base case $m=0$ obeys the Mills-type bound
\begin{subequations}\label{eq:gauss_tail}
\begin{equation}\label{eq:gauss_tail_mills}
    I_0(u;a)\;=\;\int_u^\infty e^{-ax^2/2}\,\mathrm{d}x\;\le\;\frac{e^{-au^2/2}}{a\,u}\,,
\end{equation}
and $I_1(u;a)=a^{-1}e^{-au^2/2}$ exactly. For every integer $m\ge 2$,
\begin{equation}\label{eq:gauss_tail_rec}
    I_m(u;a) \;=\; \frac{u^{m-1}\,e^{-au^2/2}}{a}\;+\;\frac{m-1}{a}\,I_{m-2}(u;a)\,,
\end{equation}
and, provided $(m-1)/(au^2)<1$,
\begin{equation}\label{eq:gauss_tail_closed}
    I_m(u;a)\;\le\;\frac{u^{m-1}\,e^{-au^2/2}}{a-(m-1)/u^2}\,.
\end{equation}
\end{subequations}
In particular, 
\begin{align*}
    I_{d-1}(u;1)    &\le 2u^{d-2}e^{-u^2/2}\,,
    &\text{for }u^2\ge 2(d-2)\,,\\
    I_{2d-2}(u;3/2) &\le u^{2d-3}e^{-3u^2/4}\,,
    &\text{for }u\ge 2\sqrt d\,.
\end{align*}
\end{lemma}

\begin{proof}
On $[u,\infty)$, $1\le x/u$, so $I_0(u;a)\le u^{-1}\int_u^\infty x\,e^{-ax^2/2}\,\mathrm{d}x=e^{-au^2/2}/(au)$, which is~\eqref{eq:gauss_tail_mills}; $I_1$ is explicit. For $m\ge2$, integrating $x^{m-1}\cdot x\,e^{-ax^2/2}$ by parts gives~\eqref{eq:gauss_tail_rec}. Iterating it down to $I_1$ or $I_0$, bounding $I_0$ by~\eqref{eq:gauss_tail_mills} and each product $(m-1)(m-3)\cdots(m-2j+1)$ by $(m-1)^j$,
\[
    I_m(u;a)\;\le\;\frac{u^{m-1}e^{-au^2/2}}{a}\sum_{j\ge0}\Big(\frac{m-1}{au^2}\Big)^{j},
\]
a geometric series whose sum is~\eqref{eq:gauss_tail_closed} when $(m-1)/(au^2)<1$. The two particular cases follow from~\eqref{eq:gauss_tail_closed}: in the respective ranges the denominators $1-(d-2)/u^2$ and $\tfrac32-(2d-3)/u^2$ are at least $\tfrac12$ and $1$.
\end{proof}

\subsection{Asymptotic leading term}
\begin{lemma}[Leading term]\label{lem:exp_det}
As $\nu\to\infty$, $A(\nu)\sim\nu^{d-1}$; consequently
$\int_u^\infty A(\rho x)\varphi(x)\,\mathrm{d}x\sim\int_u^\infty(\rho x)^{d-1}\varphi(x)\,\mathrm{d}x$ as $u\to\infty$.
\end{lemma}
\begin{proof}
For $\nu\ge1$, $\nu^{-(d-1)}|\det(G-\nu I)|=\prod_i|1-\mu_i/\nu|\le(1+\max_i|\mu_i|)^{d-1}$, which is integrable, and the product tends to $1$ almost surely; dominated convergence gives the first claim. For the second, given~$\varepsilon\in(0,1)$ pick $x_\varepsilon$ with $|A(\rho x)/(\rho x)^{d-1}-1|\le\varepsilon$ for $x\ge x_\varepsilon$; for $u\ge x_\varepsilon$ the ratio of the two integrals lies in $[1-\varepsilon,1+\varepsilon]$. Finally, \eqref{eq:gauss_tail_rec} with $a=1$, $m=d-1$ gives $I_{d-1}(u;1)\ge u^{d-2}e^{-u^2/2}$ and \eqref{eq:gauss_tail_closed} gives $I_{d-1}(u;1)\le u^{d-2}e^{-u^2/2}/(1-(d-2)/u^2)$, so $\int_u^\infty x^{d-1}\varphi(x)\,\mathrm{d}x\sim u^{d-2}\varphi(u)$, and $C_{k,d}\,\rho^{d-1}/\sqrt{2\pi}=\sqrt2\,(k/2)^{(d-1)/2}/\Gamma(d/2)$ yields~\eqref{eq:asymptotic_baseline_approx}.
\end{proof}

\subsection{Explicit Hermite tail integrals}

\begin{lemma}[Explicit Hermite tail integrals]
\label{lem:hermite_tail_explicit}
For every $u\in\R$ and integer $m\ge 1$,
\begin{equation}\label{eq:hermite_tail_series_1}
    \int_u^\infty H_m(x)\,e^{-x^2/2}\,\mathrm{d}x
    \;=\; e^{-u^2/2}\sum_{\ell=0}^{\lfloor (m-1)/2\rfloor}2^{\ell+1}\frac{(m-1)!!}{(m-2\ell-1)!!}H_{m-2\ell-1}(u)\;+\;R_m(u)\,,
\end{equation}
where $R_m(u)=0$ if $m$ is odd and $R_m(u)=2^{m/2}(m-1)!!\int_u^\infty e^{-x^2/2}\,\mathrm{d}x$ if $m$ is even.
\end{lemma}

\begin{proof}

\medskip

\noindent
$\bullet$
We evaluate $J_m = \int_u^\infty H_m(x) e^{-x^2/2} \,\mathrm{d}x$. The integrations by parts below use only the finiteness of the boundary term $H_{m-1}(u)\,e^{-u^2/2}$ at $u$ and the vanishing of the integrand at $+\infty$, so no positivity of $u$ is required.
We use two standard identities for the physicist Hermite polynomials:
\begin{subequations}\label{eq:hermite_identities}
\begin{align}
    H_m'(x) &= 2m H_{m-1}(x) \label{eq:hermite_deriv} \\
    H_m(x) &= 2x H_{m-1}(x) - 2(m-1)H_{m-2}(x) \label{eq:hermite_recurrence}
\end{align}
\end{subequations}
We begin by evaluating the integral of $2x H_{m-1}(x) e^{-x^2/2}$ using integration by parts. Let $f(x) = H_{m-1}(x)$ and $g'(x) = x e^{-x^2/2}$, which implies $g(x) = -e^{-x^2/2}$.
\begin{align*}
    \int_u^\infty x H_{m-1}(x) e^{-x^2/2} \,\mathrm{d}x &= \left[ -H_{m-1}(x) e^{-x^2/2} \right]_u^\infty + \int_u^\infty H_{m-1}'(x) e^{-x^2/2} \,\mathrm{d}x \\
    &= H_{m-1}(u) e^{-u^2/2} + 2(m-1) \int_u^\infty H_{m-2}(x) e^{-x^2/2} \,\mathrm{d}x
\end{align*}
where we substituted \eqref{eq:hermite_deriv} in the second term. We then substitute $2x H_{m-1}(x)$ using a rearrangement of \eqref{eq:hermite_recurrence}, $2x H_{m-1}(x) = H_m(x) + 2(m-1)H_{m-2}(x)$, into the left-hand side:
\begin{equation*}
    \frac{1}{2} \int_u^\infty \Big( H_m(x) + 2(m-1)H_{m-2}(x) \Big) e^{-x^2/2} \,\mathrm{d}x = H_{m-1}(u) e^{-u^2/2} + 2(m-1) \int_u^\infty H_{m-2}(x) e^{-x^2/2} \,\mathrm{d}x
\end{equation*}
Distributing the integral and isolating $J_m = \int_u^\infty H_m(x) e^{-x^2/2} \,\mathrm{d}x$ yields the exact recurrence relation:
\begin{equation}\label{eq:J_recurrence}
    J_m = 2H_{m-1}(u) e^{-u^2/2} + 2(m-1) J_{m-2}
\end{equation}
By iteratively applying \eqref{eq:J_recurrence}, we can expand $J_m$:
\begin{align*}
    J_m &= 2H_{m-1}(u) e^{-u^2/2} + 2(m-1) \Big[ 2H_{m-3}(u) e^{-u^2/2} + 2(m-3) J_{m-4} \Big] \\
    &= e^{-u^2/2} \Big[ 2H_{m-1}(u) + 2^2(m-1)H_{m-3}(u) + 2^3(m-1)(m-3)H_{m-5}(u) + \dots \Big] + \text{Remainder}
\end{align*}
The $\ell$-th term of this expansion involves the cascading product $(m-1)(m-3)\cdots(m-2\ell+1)$, which evaluates precisely to the double factorial ratio $\frac{(m-1)!!}{(m-2\ell-1)!!}$. Incorporating the accumulated factor of $2^{\ell+1}$, the summation becomes:
\begin{equation*}
    e^{-u^2/2} \sum_{\ell=0}^{\lfloor (m-1)/2 \rfloor} 2^{\ell+1} \frac{(m-1)!!}{(m-2\ell-1)!!} H_{m-2\ell-1}(u)
\end{equation*}
The recursion stops differently depending on parity:
\begin{itemize}
    \item If $m$ is odd, the final evaluation step corresponds to $\ell = (m-1)/2$, meaning $H_0(u)$ is evaluated and no subsequent integrals are required. Thus, $R_m(u) = 0$.
    \item If $m$ is even, the final evaluated boundary term corresponds to $\ell = (m-2)/2$. The remaining integral falls back to $J_0 = \int_u^\infty H_0(x) e^{-x^2/2}\,\mathrm{d}x = \int_u^\infty e^{-x^2/2}\,\mathrm{d}x$. The accumulated multiplier in front of $J_0$ is $2^{m/2} (m-1)!!$, yielding exactly the prescribed $R_m(u)$.
\end{itemize}

\end{proof}

\medskip

The next lemma extends Lemma~\ref{lem:hermite_tail_explicit} from the canonical Gaussian weight $e^{-y^2/2}$ to an arbitrary weight $e^{-\beta y^2/2}$ with $\beta>0$. The general $\beta$-version gives the exact squared-Hermite tail $T_d^{\mathrm{exact}}(u)$ of the $D_4$ step in the proof of Theorem~\ref{thm:imf_tail_exact}: when integrating $H_j(\rho x)^2$ against $e^{-(1+\rho^2)x^2/2}$ on $[u,\infty)$ and substituting~$y=\rho x$, the relevant exponent is $\beta=(1+\rho^2)/\rho^2=(3k-2)/k>2$ for $k\ge 3$, a regime to which Lemma~\ref{lem:hermite_tail_explicit} ($\beta=1$) does not apply.

\begin{lemma}[Generalized Hermite--Gaussian tail recurrence]
\label{lem:Jm_beta}
For every $\beta>0$, every integer $m\ge 0$, and every $u\in\R$, set
\begin{subequations}\label{eq:Jm_beta}
\begin{equation}\label{eq:Jm_beta_def}
    J_m^\beta(u)\;:=\;\int_u^\infty H_m(y)\,e^{-\beta y^2/2}\,\mathrm{d}y\,.
\end{equation}
Write $\gamma:=2/\beta$ and $\vartheta:=(2-\beta)/\beta$. The base cases are
\begin{equation}\label{eq:Jm_beta_init}
    J_0^\beta(u)\;=\;\sqrt{\tfrac{2\pi}\beta}\,\bar\Phi(u\sqrt\beta)\,,\qquad
    J_1^\beta(u)\;=\;\gamma\,e^{-\beta u^2/2}\,,
\end{equation}
and for every $m\ge 2$, the family $\{J_m^\beta\}$ satisfies the three-term recurrence
\begin{equation}\label{eq:Jm_beta_rec}
    J_m^\beta(u)\;=\;\gamma\,H_{m-1}(u)\,e^{-\beta u^2/2}\;+\;2(m-1)\,\vartheta\;J_{m-2}^\beta(u)\,.
\end{equation}
Iterating~\eqref{eq:Jm_beta_rec} yields the explicit closed form
\begin{equation}\label{eq:Jm_beta_iter}
    J_m^\beta(u)\;=\;\gamma\,e^{-\beta u^2/2}\!\!\sum_{\ell=0}^{\lfloor (m-1)/2\rfloor}(2\vartheta)^\ell\,\frac{(m-1)!!}{(m-2\ell-1)!!}\,H_{m-2\ell-1}(u)\;+\;R_m^\beta(u)\,,
\end{equation}
\end{subequations}
where $R_m^\beta(u)=0$ if $m$ is odd, and $R_m^\beta(u)=(2\vartheta)^{m/2}(m-1)!!\,J_0^\beta(u)$ if $m$ is even.
\end{lemma}

\noindent
For $\beta=1$, $\gamma=2$ and $\vartheta=1$, so $\gamma(2\vartheta)^\ell=2^{\ell+1}$, and Lemma~\ref{lem:Jm_beta} reduces exactly to~\eqref{eq:hermite_tail_series_1}; the remainder $R_m^\beta$ specializes to $2^{m/2}(m-1)!!\int_u^\infty e^{-y^2/2}\,\mathrm{d}y$, matching $R_m$. For $\beta>2$, $\vartheta<0$, and the iterated sum~\eqref{eq:Jm_beta_iter} carries alternating signs in $\ell$ but remains exact pointwise; in subsequent estimates one must not pass to absolute values of individual summands.

\begin{proof}[Proof of Lemma~\ref{lem:Jm_beta}]
\noindent
$\bullet$
Apply integration by parts with $f(y)=H_{m-1}(y)$ and $g'(y)=2y\,e^{-\beta y^2/2}$, hence~$g(y)=-(2/\beta)\,e^{-\beta y^2/2}=-\gamma\,e^{-\beta y^2/2}$. The boundary term at $+\infty$ vanishes by Gaussian decay, and using $H_{m-1}'(y)=2(m-1)\,H_{m-2}(y)$,
\[
    \int_u^\infty 2y\,H_{m-1}(y)\,e^{-\beta y^2/2}\,\mathrm{d}y
    \;=\;\gamma\,H_{m-1}(u)\,e^{-\beta u^2/2}\;+\;2(m-1)\,\gamma\,J_{m-2}^\beta(u)\,.
\]
Combining with the Hermite recurrence $H_m(y)=2y\,H_{m-1}(y)-2(m-1)\,H_{m-2}(y)$,
\[
    J_m^\beta(u) \;=\; \int_u^\infty 2y\,H_{m-1}(y)\,e^{-\beta y^2/2}\,\mathrm{d}y\;-\;2(m-1)\,J_{m-2}^\beta(u)\,,
\]
and $2(m-1)\gamma-2(m-1)=2(m-1)(\gamma-1)=2(m-1)\vartheta$, yielding~\eqref{eq:Jm_beta_rec}.

\medskip

\noindent
$\bullet$ For the boundary cases, $J_0^\beta(u)=\int_u^\infty e^{-\beta y^2/2}\,\mathrm{d}y=\sqrt{2\pi/\beta}\,\bar\Phi(u\sqrt\beta)$ by direct change of variable, and~$J_1^\beta(u)=\int_u^\infty 2y\,e^{-\beta y^2/2}\,\mathrm{d}y=(2/\beta)\,e^{-\beta u^2/2}=\gamma\,e^{-\beta u^2/2}$.

\medskip

\noindent
$\bullet$ Unrolling~\eqref{eq:Jm_beta_rec} $\ell$ times and using the telescoping product $(m-1)(m-3)\cdots(m-2\ell+1)=(m-1)!!/(m-2\ell-1)!!$,
\[
    J_m^\beta(u)\;=\;\gamma\,e^{-\beta u^2/2}\!\!\sum_{j=0}^{\ell-1}(2\vartheta)^j\,\frac{(m-1)!!}{(m-2j-1)!!}\,H_{m-2j-1}(u)\;+\;(2\vartheta)^\ell\,\frac{(m-1)!!}{(m-2\ell-1)!!}\,J_{m-2\ell}^\beta(u)\,.
\]
For $m$ odd, take $\ell=(m-1)/2$: $m-2\ell-1=0$, hence $(m-2\ell-1)!!=0!!=1$, the residual factor $J_{m-2\ell}^\beta(u)=J_1^\beta(u)=\gamma\,e^{-\beta u^2/2}$, and the residual term coincides with the $j=(m-1)/2$ summand of the main sum (since $H_{m-2\ell-1}(u)=H_0(u)=1$ and $\gamma e^{-\beta u^2/2}=J_1^\beta(u)$); absorbing it into the main sum gives $R_m^\beta=0$. For $m$ even, take $\ell=m/2$: $m-2\ell-1=-1$, $(m-2\ell-1)!!=(-1)!!=1$, the residual factor is $J_0^\beta(u)$, and $R_m^\beta(u)=(2\vartheta)^{m/2}\,(m-1)!!\,J_0^\beta(u)$.

\medskip

\end{proof}

\begin{corollary}[Squared-Hermite Gaussian tail]
\label{cor:Kj_beta}
For every $\beta>0$, every integer $j\ge 0$, and every $u\in\R$, the squared-Hermite tail
\begin{subequations}\label{eq:Kj_beta}
\begin{equation}\label{eq:Kj_beta_def}
    K_j^\beta(u)\;:=\;\int_u^\infty H_j(y)^2\,e^{-\beta y^2/2}\,\mathrm{d}y
\end{equation}
admits the exact closed form
\begin{equation}\label{eq:Kj_beta_form}
    K_j^\beta(u)\;=\;\sum_{p=0}^j 2^p\,p!\,\binom{j}{p}^{\!2}\,J_{2j-2p}^\beta(u)\,,
\end{equation}
\end{subequations}
where each $J_{2j-2p}^\beta$ is given in closed form by Lemma~\ref{lem:Jm_beta}.
\end{corollary}

\begin{proof}[Proof of Corollary~\ref{cor:Kj_beta}]
The classical linearization of the product of two physicist Hermite polynomials reads
\[
    H_m(y)\,H_n(y)\;=\;\sum_{p=0}^{\min(m,n)} 2^p\,p!\,\binom{m}{p}\binom{n}{p}\,H_{m+n-2p}(y)\,.
\]
Specializing to $m=n=j$ yields $H_j(y)^2=\sum_{p=0}^j 2^p\,p!\,\binom{j}{p}^{\!2}\,H_{2j-2p}(y)$. Multiplying by $e^{-\beta y^2/2}$, integrating from $u$ to $\infty$, and exchanging the finite sum with the integral (justified termwise since each $H_m\,e^{-\beta y^2/2}$ is absolutely integrable on $[u,\infty)$), one obtains~\eqref{eq:Kj_beta_form}.
\end{proof}

\section*{Acknowledgments}
Yohann De Castro acknowledges the support of the Institut Universitaire de France (IUF).

\section*{Statements and Declarations}
\paragraph{Competing interests.} The authors have no competing interests to declare that are relevant to the content of this article.
\paragraph{Data availability.} The scripts that produce the figures and the numerical values reported in this article are available at \verifrepo. No other data were generated.

\section*{Disclosure of AI use}
In line with the Leiden Declaration on Artificial Intelligence and Mathematics
\citep{leiden2026declaration}, we disclose the use of AI agents (Claude Code, Anthropic) in the
preparation of this paper. These agents were used to audit the statements, the proofs and the
numerical tables, through repeated adversarial verification passes recorded in a change-log, and
to help with editorial tasks: consistency of notation, bibliography formatting, and typesetting.
The authors produced and checked every result and every proof themselves and take sole
responsibility for the correctness of the paper and for its citations.

\end{document}